# EXISTENCE OF MULTI-SOLITARY WAVES WITH LOGARITHMIC RELATIVE DISTANCES FOR THE NLS EQUATION

NGUYỄN TIẾN VINH

ABSTRACT. We construct in this paper global (for $t \geq 0$) and bounded solutions $u(t)$ for the nonlinear Schrödinger equation

$$i\partial_t u + \Delta u + |u|^{p-1}u = 0, \qquad t \in \mathbb{R}, \ x \in \mathbb{R}^d$$

in mass sub-critical cases ($1 < p < 1 + \frac{4}{d}$) and mass super-critical cases ($1 + \frac{4}{d} < p < \frac{d+2}{d-2}$) such that $u(t)$ decomposes asymptotically into two solitary waves with logarithmic distance

$$\left\| u(t) - e^{i\gamma(t)} \sum_{k=1}^{2} Q(. - x_k(t)) \right\|_{H^1} \to 0$$

and

$$|x_1(t) - x_2(t)| \sim 2\log t, \qquad \text{as } t \to +\infty.$$

The logarithmic distance is related to strong interactions between solitary waves. In the integrable case ($d = 1$ and $p = 3$) the existence of such solutions has been shown in [14].

## 1. INTRODUCTION

We consider the nonlinear Schrödinger equation in $\mathbb{R}^d$, for any $d \geq 1$:

$$\begin{cases} i\partial_t u = -\Delta u - |u|^{p-1}u, & (t,x) \in [0,T) \times \mathbb{R}^d \\ u(0,x) = u_0, & u_0 \in H^1 : \mathbb{R}^d \to \mathbb{C}. \end{cases} \tag{NLS}$$

It is well-known (see e.g. [2], [9]) that the equation (NLS) is locally well-posed in $H^1(\mathbb{R}^d)$ for $1 < p < \frac{d+2}{d-2}$: for any $u_0 \in H^1(\mathbb{R}^d)$, there exist $T^* > 0$ and a unique maximal solution $u \in \mathcal{C}([0, T^*), H^1(\mathbb{R}^d))$ of (NLS). Moreover, the following blow up criterion holds

$$T^* < +\infty \quad \text{implies} \quad \lim_{t \uparrow T^*} \|\nabla u(t)\|_{L^2} = +\infty. \tag{1.1}$$

Recall that the solution $u$ satisfies the following three conservation laws:

- Mass :

$$\int_{\mathbb{R}^d} |u(t,x)|^2 dx = \int_{\mathbb{R}^d} |u_0(x)|^2 \tag{1.2}$$

- Energy :

$$E(u(t)) = \frac{1}{2} \int_{\mathbb{R}^d} |\nabla u(t,x)|^2 - \frac{1}{p+1} \int_{\mathbb{R}^d} |u(t,x)|^{p+1} dx = E(u_0) \tag{1.3}$$

- Momentum:

$$M(u(t)) = \text{Im} \int_{\mathbb{R}^d} \nabla u(t,x)\bar{u}(t,x)dx = M(u_0) \tag{1.4}$$

for all $t \in [0, T^*)$. Recall also that (NLS) admits the following symmetries: the transformation of initial data implies the corresponding transformation of solution:

- Scaling: $\lambda > 0, \lambda^{\frac{2}{p-1}} u_0(\lambda x) \mapsto \lambda^{\frac{2}{p-1}} u(\lambda^2 t, \lambda x)$;
- Space translation: $x_0 \in \mathbb{R}^d, u_0(x + x_0) \mapsto u(t, x + x_0)$;
- Time translation: $t_0 \in \mathbb{R}, u_{t_0}(x) \mapsto u(t + t_0, x)$;
- Space rotation: $A \in SO(d), u_0(A \cdot x_0) \mapsto u(t, A \cdot x_0)$;
- Phase: $\gamma \in \mathbb{R}, u_0(x)e^{i\gamma} \mapsto u(t, x)e^{i\gamma}$;
- Galilean: $\beta \in \mathbb{R}^d, u_0(x)e^{i\beta x} \mapsto u(t, x - \beta t)e^{i\frac{\beta}{2}(x - \frac{\beta}{2}t)}$.

As a consequence of (1.2), (1.3) and the Gagliardo-Nirenberg inequality (see [2]), all solutions of (NLS) are global in $L^2$ sub-critical cases $(1 < p < 1 + \frac{4}{d})$, whereas blow-up solutions exist in $L^2$ critical case $(p = 1 + \frac{4}{d})$ and $L^2$ super-critical cases $(1 + \frac{4}{d} < p < \frac{d+2}{d-2})$.

This article is concerned with the construction of special solutions of (NLS) involving solitary wave solutions (or solitons). We recall the solitary wave: for $\lambda_0 > 0$,

$$u(t, x) = e^{i\lambda_0^2 t} Q_{\lambda_0}(x) \quad \text{with} \quad Q_{\lambda_0}(x) = \lambda_0^{\frac{2}{p-1}} Q(\lambda_0 x)$$

where $Q$, the ground state, is the unique radial positive solution (up to symmetries) of

$$\Delta Q - Q + Q^p = 0, \quad Q > 0, \quad Q \in H^1(\mathbb{R}^d). \tag{1.5}$$

For more properties of $Q$, see for example [2] and [23]. Recall also that for $L^2$ sub-critical cases $(1 < p < 1 + \frac{4}{d})$, the solitary waves are stable ([26]) and for $L^2$ critical and $L^2$ super-critical cases $(1 + \frac{4}{d} \leq p < \frac{d+2}{d-2})$, the solitary waves are unstable ([10]).

### 1.1. **Main result.**

In this article, we prove the existence of multi-solitary wave solutions with logarithmic relative distances by exhibiting a general non free Galilean motion due to strong interactions between solitons.

**Main Theorem** (Multi-solitary waves with logarithmic distance). *Let $1 < p < 1 + \frac{4}{d}$ (mass sub-critical cases) or $1 + \frac{4}{d} < p < \frac{d+2}{d-2}$ (mass super-critical cases). Then there exists an $H^1$ solution $u(t)$ of* (NLS) *on $[0, +\infty)$ which decomposes asymptotically into two solitary waves*

$$\left\| u(t) - e^{i\gamma(t)} \sum_{k=1}^{2} Q(. - x_k(t)) \right\|_{H^1} \lesssim \frac{1}{t} \tag{1.6}$$

*and*

$$|\delta x(t)| := |x_1(t) - x_2(t)| = 2(1 + o(1)) \log t, \quad \text{as} \ \ t \to +\infty. \tag{1.7}$$

Note that we restrict ourselves to a "two-bubble" solution, but the same proof applies to any number $K \geq 2$ of solitons located on a regular polygon of size $\log t$ as in [19]. Our result holds for general (NLS) both in $L^2$ sub-critical and $L^2$ super-critical, moreover, by scaling, we can replace $Q$ by $Q_{\lambda_0}$ for any $\lambda_0 > 0$. For $L^2$ critical, we refer to [19] and remarks below. We observe also that in the result, solitons need to have the same sign, the same scaling and the same phase, in fact, the solution is symmetric by $\tau : x \mapsto -x$.

**Remark 1.** *For* (NLS), *multiple bubble solutions with weak interactions and asymptotically free Galilean motion have been constructed in various settings, both in stable and unstable* *contexts, see in particular [4, 15, 20]. As a typical illustration of weakly interacting dynamics, there exist multi-solitary wave solutions of* (NLS) *satisfying for large t,*

$$\left\| u(t) - \sum_{k=1}^{K} e^{-i\Gamma_k(t,x)} Q_{\lambda_k}(. - \nu_k t) \right\|_{H^1} \lesssim e^{-\gamma t}, \quad \gamma > 0, \tag{1.8}$$

*for any given set of parameters* $\{\nu_k, \lambda_k\}_k \in \mathbb{R}^d \times (0, \infty)$ *with the decoupling condition* $\nu_k \neq \nu_{k'}$ *if* $k \neq k'$. *In the integrable case* $(d = 1$ *and* $p = 3)$, *the existence of such solutions is well-known, see [27] for a derivation of their explicit expression. Moreover, these solutions are very special: they describe the perfect interaction between several solitary waves and are actually global K-pure solitons. In the non-integrable cases, little is known on the behavior of such solutions as* $t \to -\infty$.

*Concerning strong interactions, for* (NLS) *equation, we recall Theorem 1 in [19]: in* $L^2$ *critical two dimensional case, there exists a global (for* $t \geq 0$) *solution* $u(t)$ *that decomposes asymptotically into a sum of solitary waves*

$$\left\| u(t) - e^{i\gamma(t)} \sum_{k=1}^{K} \frac{1}{\lambda(t)} Q\left( \frac{. - x_k(t)}{\lambda(t)} \right) \right\|_{H^1} \to 0, \quad \lambda(t) = \frac{1 + o(1)}{\log t} \quad as \quad t \to +\infty, \tag{1.9}$$

*where the translation parameters* $x_k(t)$ *converge to the vertices of a K-sided regular polygon and the solution blows up in infinite time with the rate*

$$\|\nabla u(t)\|_{L^2} \sim |\log t| \quad as \quad t \to +\infty.$$

*We also refer to [12, 13, 17] for works where refined analysis of interactions between solitons is important for other nonlinear equations.*

**Remark 2.** *The existence of the regime described in Main Theorem is not a surprise. In the integrable case* $(d = 1, p = 3)$, *the existence of 2-soliton solutions with logarithmic distance, called double pole solutions, was reported in [14] and the dynamics of interacting pulses in several models was formally studied in [8]. However, to our knowledge, the present work may be the first general proof of existence of global and bounded multi-solitons with logarithmic distances. We expect such solutions to be unstable, even in* $L^2$ *sub-critical cases, since generic perturbation can give collision or on the contrary weak interaction. The appearance of the regime in Main Theorem is closely related to the equation*

$$\ddot{z}(t) = -e^{-2z(t)}$$

*where* $\log t$ *is a solution with initial conditions* $z(1) = 0, \dot{z}(1) = 1$. *From the theory of perturbation, for* $z(t) = \log t + \epsilon v_1 + ...$ *with initial conditions* $z(1) = \epsilon, \dot{z}(1) = 1$, *one has at the linear level*

$$\ddot{v}_1 = \frac{2v_1}{t^2}, \quad v_1(1) = 1, \quad \dot{v}_1(1) = 0,$$

*whose solution is* $\frac{1}{3}t^2 + \frac{2}{3}\frac{1}{t}$ *so we see that the* $\log t$ *solution is an unstable state as* $t \to +\infty$.

**Remark 3.** *For the* $L^2$ *critical case* $(p = 1 + \frac{4}{d})$, *interestingly enough, the existence of bounded multi-solitary wave solutions with logarithmic distances as* (1.6)–(1.7) *in Main Theorem is ruled out by the nonlinear instability related to degeneracy of the scaling direction: for such solutions, one would have*

$$\int_{\mathbb{R}^d} |x|^2 |u(t,x)|^2 dx \sim \log^2(t) \tag{1.10}$$

*which is in contradiction to the virial identity*

$$\frac{d^2}{dt^2} \int_{\mathbb{R}^d} |x|^2 |u|^2 = 16 E(u_0).$$

*In fact, the scaling instability directions in the critical case are excited by the nonlinear interactions which leads to the infinite time concentration displayed in (1.9) (see [19]).*

*In the $L^2$ super-critical cases, solitons have an exponential instability, however such instability can be controlled by a topological argument used in [4] (see Section 6).*

**Remark 4.** *The distance (1.7) between solitons in the Main Theorem can be described asymptotically as*

$$|\delta x(t)| = 2 \log t - \frac{d-1}{2} \log(\log t) - C + O(\log^{-\frac{1}{2}}(t)) \quad as \ t \to +\infty$$

*where $C > 0$ a constant depending only on $d$ and $p$ (see (3.21)).*

The article is organized as follows. Sections 2, 3 and 4 concern the proof of Main Theorem in $L^2$ sub-critical cases with $p > 2$. In Section 2, we consider an approximate solution (an ansatz solution) to (NLS) made of two symmetric "bubbles" and extract the formal evolution system of the geometrical parameters of the bubbles (scaling, position, phase). The key observation is that this system contains forcing terms due to the nonlinear interactions of the waves, and has a special solution corresponding at the main order to the regime of Main Theorem. In Section 3, we construct, using modulation, particular backwards solutions of (NLS) related to the special regime of Main Theorem and prove backward uniform estimates by energy method. In Section 4, we use compactness arguments on a suitable sequence of such backwards solutions to finish the proof. Sections 5 deals with the case $1 < p \leq 2$; in this case, there are some extra technical difficulties, even if the strategy of the proof is similar: the interaction becomes stronger, we have to add extra terms in the approximate solution and due to lost of regularity, we have to use some truncations. Finally, the algebraic computations in the proof for $L^2$ sub-critical cases are still valid in $L^2$ super-critical cases, Section 6 presents additional arguments and modifications needed for $L^2$ super-critical cases.

1.2. **Notation.** The $L^2$ scalar product of two complex valued functions $f, g \in L^2(\mathbb{R}^d)$ is denoted by

$$\langle f, g \rangle = \text{Re} \left( \int_{\mathbb{R}^d} f(x) \overline{g}(x) dx \right).$$

We denote by $Q(x) := q(|x|)$ the unique radial positive ground state of (NLS):

$$q'' + \frac{d-1}{r} q' - q + q^p = 0, \quad q'(0) = 0, \quad \lim_{r \to +\infty} q(r) = 0. \tag{1.11}$$

It is well-known and easily checked by ODE arguments that for some constant $c_Q > 0$,

$$\text{for all } r > 1, \quad \left| q(r) - c_Q r^{-\frac{d-1}{2}} e^{-r} \right| + \left| q'(r) + c_Q r^{-\frac{d-1}{2}} e^{-r} \right| \lesssim r^{-\frac{d-1}{2}-1} e^{-r}. \tag{1.12}$$

We set

$$I_Q = \int Q^p(x) e^{-x_1} dx, \quad x = (x_1, ..., x_d).$$

We denote by $\mathcal{Y}$ the set of smooth functions $f$ such that

$$\text{for all } p \in \mathbb{N}, \text{ there exists } q \in \mathbb{N}, \text{ s.t. for all } x \in \mathbb{R}^d, \quad |f^{(p)}(x)| \lesssim |x|^q e^{-|x|}. \tag{1.13}$$

Let $\Lambda$ be the generator of $L^2$-scaling corresponding to (NLS):

$$\Lambda f = \frac{2}{p-1} f + x \cdot \nabla f.$$

The linearization of (NLS) around $Q$ involves the following Schrödinger operators:

$$L_+ := -\Delta + 1 - pQ^{p-1}, \qquad L_- := -\Delta + 1 - Q^{p-1}.$$

From [25], recall the generalized null space relations in sub-critical and super-critical cases:

$$\begin{aligned}
& L_- Q = 0, \quad L_+(\Lambda Q) = -2Q, \\
& L_+(\nabla Q) = 0, \quad L_-(xQ) = -2\nabla Q.
\end{aligned} \tag{1.14}$$

We recall the coercivity property in $L^2$ sub-critical (see [15], [21], [25], [26]): there exists $\mu > 0$ such that for all $\eta \in H^1$,

$$\langle L_+ \operatorname{Re}\eta, \operatorname{Re}\eta \rangle + \langle L_- \operatorname{Im}\eta, \operatorname{Im}\eta \rangle \geq \mu \|\eta\|_{H^1}^2 - \frac{1}{\mu} \left( \langle \eta, Q \rangle^2 + |\langle \eta, xQ \rangle|^2 + \langle \eta, i\Lambda Q \rangle^2 \right). \tag{1.15}$$

In $L^2$ super-critical (but $H^1$ sub-critical), we do not have the same situation since the negative direction can not be controlled by the scaling parameter. We consider the operator

$$\mathcal{L}v = iL_+ v_1 - L_- v_2 \ \text{ with } v = v_1 + iv_2.$$

The spectrum $\sigma(\mathcal{L})$ of $\mathcal{L}$ satisfies

$$\sigma(\mathcal{L}) \cap \mathbb{R} = \{-e_0, 0, e_0\}.$$

It is easy to see that $iQ, \nabla Q$ are independent and belong to the kernel of $\mathcal{L}$. In [4], [6], [7], [11], it is proved that there exist two eigenfunctions $Y^\pm$ (normalized by $\|Y^\pm\|_{L^2} = 1$) associated to eigenvalues $\pm e_0$

$$\mathcal{L}(Y^\pm) = \pm e_0 Y^\pm \tag{1.16}$$

and $Y^+ = \overline{Y^-}$ belong to $\mathcal{Y}$, in other words, $\operatorname{Re} Y^+, \operatorname{Im} Y^+ \in \mathcal{Y}$. Moreover, there holds a property of positivity based on $Y^\pm$: there exists $\mu > 0$ such that for all $\eta \in H^1$,

$$\langle L_+ \operatorname{Re}\eta, \operatorname{Re}\eta \rangle + \langle L_- \operatorname{Im}\eta, \operatorname{Im}\eta \rangle \geq \mu \|\eta\|_{H^1}^2$$
$$- \frac{1}{\mu} \left( \langle \eta, iY^+ \rangle^2 + \langle \eta, iY^- \rangle^2 + |\langle \eta, xQ \rangle|^2 + \langle \eta, i\Lambda Q \rangle^2 \right). \tag{1.17}$$

## 2. Approximate solution for $p > 2$

2.1. **System of modulation equations.** Let $p > 2$. Consider a time dependent $\mathcal{C}^1$ function of parameters $\vec{q}$ of the form

$$\vec{q} = (\lambda, z, \gamma, v) \in (0, +\infty) \times \mathbb{R}^d \times \mathbb{R} \times \mathbb{R}^d,$$

with $|v| \ll 1$ and $|z| \gg 1$. We renormalize the flow by considering

$$u(t,x) = \frac{e^{i\gamma(s)}}{\lambda^{\frac{2}{p-1}}(s)} w(s,y), \quad dt = \lambda^2(s)ds, \quad y = \frac{x}{\lambda(s)}, \tag{2.1}$$

so that

$$i\partial_t u + \Delta u + |u|^{p-1} u = \frac{e^{i\gamma}}{\lambda^{2+\frac{2}{p-1}}} \left[ i\dot{w} + \Delta w - w + |w|^{p-1} w - i\frac{\dot{\lambda}}{\lambda}\Lambda w + (1 - \dot{\gamma})w \right] \tag{2.2}$$

($\dot{w}$ denotes derivation with respect to $s$). We introduce the following $\vec{q}$-modulated ground state solitary waves, for $k \in \{1,2\}$,

$$P_k(s,y) = e^{i\Gamma_k(s,y-z_k(s))}Q(y-z_k(s)) = e^{iv_k(s)(y-z_k(s))}Q(y-z_k(s)), \qquad (2.3)$$

where we set

$$v_1(s) = -v_2(s) = \frac{1}{2}v(s), \quad z_1(s) = -z_2(s) = \frac{1}{2}z(s), \quad \Gamma_k(s,y) = v_k(s) \cdot y, \qquad (2.4)$$

Let

$$\mathbf{P}(s,y) = \mathbf{P}(y;(z(s),v(s))) = \sum_{k=1}^{2} P_k(s,y). \qquad (2.5)$$

Then, $\mathbf{P}$ is an approximate solution of the rescaled equation in the following sense.

**Lemma 1** (Leading order approximate flow). *Let the vectors of modulation equations be*

$$\vec{m}_k = \begin{pmatrix} \frac{\dot{\lambda}}{\lambda} \\ \dot{z}_k - 2v_k + \frac{\dot{\lambda}}{\lambda}z_k \\ \dot{\gamma} - 1 + |v_k|^2 - \frac{\dot{\lambda}}{\lambda}(v_k \cdot z_k) - (v_k \cdot \dot{z}_k) \\ \dot{v}_k - \frac{\dot{\lambda}}{\lambda}v_k \end{pmatrix}, \qquad \vec{\mathrm{M}}V = \begin{pmatrix} -i\Lambda V \\ -i\nabla V \\ -V \\ -yV \end{pmatrix}. \qquad (2.6)$$

*Then the error $\mathcal{E}_{\mathbf{P}}$ to the re-normalized flow (2.2) at $\mathbf{P}$,*

$$\mathcal{E}_{\mathbf{P}} = i\dot{\mathbf{P}} + \Delta\mathbf{P} - \mathbf{P} + |\mathbf{P}|^{p-1}\mathbf{P} - i\frac{\dot{\lambda}}{\lambda}\Lambda\mathbf{P} + (1-\dot{\gamma})\mathbf{P} \qquad (2.7)$$

*decomposes as*

$$\mathcal{E}_{\mathbf{P}} = [e^{i\Gamma_1}\vec{m}_1 \cdot \vec{\mathrm{M}}Q](y-z_1(s)) + [e^{i\Gamma_2}\vec{m}_2 \cdot \vec{\mathrm{M}}Q](y-z_2(s)) + G \qquad (2.8)$$

*where the interaction term $G = |\mathbf{P}|^{p-1}\mathbf{P} - |P_1|^{p-1}P_1 - |P_2|^{p-1}P_2$ satisfies*

$$\|G\|_{L^\infty} \lesssim |z|^{-\frac{d-1}{2}}e^{-|z|}, \quad \|\nabla G\|_{L^\infty} \lesssim |z|^{-\frac{d-1}{2}}e^{-|z|}. \qquad (2.9)$$

*Proof of Lemma 1.* Firstly, we compute $\mathcal{E}_{P_k} = i\dot{P}_k + \Delta P_k - P_k + |P_k|^{p-1}P_k - i\frac{\dot{\lambda}}{\lambda}\Lambda P_k + (1-\dot{\gamma})P_k$. Let $y_{z_k} = y - z_k$, by computations

$$i\dot{P}_k = \left[ -(\dot{v}_k \cdot y_{z_k})Q(y_{z_k}) + (v_k \cdot \dot{z}_k)Q(y_{z_k}) - i\dot{z}_k \cdot \nabla Q(y_{z_k}) \right]e^{iv_k \cdot y_{z_k}}$$

$$\nabla P_k = \left[ \nabla Q(y_{z_k}) + iv_k Q(y_{z_k}) \right]e^{iv_k \cdot y_{z_k}}$$

$$\Delta P_k = \left[ \Delta Q(y_{z_k}) + 2iv_k \cdot \nabla Q(y_{z_k}) - v_k^2 Q(y_{z_k}) \right]e^{iv_k \cdot y_{z_k}}$$

$$\Lambda P_k = \left[ \frac{2}{p-1}Q(y_{z_k}) + y \cdot [\nabla Q(y_{z_k}) + iv_k Q(y_{z_k})] \right]e^{iv_k y_{z_k}}$$

$$= \left[ \Lambda Q(y_{z_k}) + iv_k \cdot y_{z_k}Q(y_{z_k}) + iv_k \cdot z_k Q(y_{z_k}) + z_k \cdot \nabla Q(y_{z_k}) \right]e^{iv_k \cdot y_{z_k}}.$$

Therefore, we get

$$\mathcal{E}_{P_k} = \left[ -i\frac{\dot{\lambda}}{\lambda}\Lambda Q(y_{z_k}) - i(\dot{z}_k - 2v_k + z_k\frac{\dot{\lambda}}{\lambda})\cdot\nabla Q(y_{z_k}) - (\dot{\gamma} - 1 - v_k\cdot\dot{z}_k + |v_k|^2 - v_k\cdot z_k\frac{\dot{\lambda}}{\lambda})Q(y_{z_k}) \right.$$

$$\left. - (\dot{v}_k - v_k\frac{\dot{\lambda}}{\lambda})\cdot y_{z_k}Q(y_{z_k}) + \Delta Q(y_{z_k}) - Q(y_{z_k}) + |Q(y_{z_k})|^{p-1}Q(y_{z_k}) \right] e^{i\Gamma_k(s, y - z_k)}.$$

Since $\Delta Q - Q + |Q|^{p-1}Q = 0$, we have

$$\mathcal{E}_{P_k} = [e^{i\Gamma_k}\vec{m}_k\cdot\vec{M}Q](y - z_k(s)). \tag{2.10}$$

Returning to the error of renormalized flow, we obtain

$$\mathcal{E}_{\mathbf{P}} = \mathcal{E}_{P_1} + \mathcal{E}_{P_2} + |\mathbf{P}|^{p-1}\mathbf{P} - \sum_{k=1}^{2} |P_k|^{p-1}P_k. \tag{2.11}$$

Next, we estimate the interaction term $G = |\mathbf{P}|^{p-1}\mathbf{P} - |P_1|^{p-1}P_1 - |P_2|^{p-1}P_2$. Clearly,

$$|G| \lesssim |P_1|^{p-1}|P_2| + |P_2|^{p-1}|P_1|.$$

We observe that for $z = z_1 - z_2$, by (1.12),

$$Q(y)Q(y-z) \lesssim (1 + |y|)^{-\frac{d-1}{2}}(1 + |y-z|)^{-\frac{d-1}{2}}e^{-|y|}e^{-|z|+|y|} \lesssim |z|^{-\frac{d-1}{2}}e^{-|z|} \tag{2.12}$$

which yields

$$|P_1|^{p-1}|P_2| \lesssim |P_1||P_2||P_1|^{p-2} \lesssim |z|^{-\frac{d-1}{2}}e^{-|z|}|P_1|^{p-2}.$$

Thus,

$$|G(s,y)| \lesssim |z|^{-\frac{d-1}{2}}e^{-|z|}\sum_{k=1}^{2} Q^{p-2}(y - z_k(s)) \tag{2.13}$$

and since $p > 2$, we get

$$||G||_{L^\infty} \lesssim |z|^{-\frac{d-1}{2}}e^{-|z|}. \tag{2.14}$$

Similarly, by (1.12),

$$||\nabla G||_{L^\infty} \lesssim |z|^{-\frac{d-1}{2}}e^{-|z|}.$$

$\square$

2.2. **Nonlinear forcing.** For the next parts of the article, we will need the first-order and the second-order approximations of $F(u) = |u|^{p-1}u$ where $u = a + ib$. We consider the expansion for $|u| \ll 1$

$$F(1 + u) = 1 + pa + ib + \frac{p(p-1)}{2}a^2 + \frac{p-1}{2}b^2 + (p-1)iab + O(|u|^k) \tag{2.15}$$

for any $2 < k \le 3$. From which, we can deduce formally

$$F'(\mathbf{P}).\epsilon = \frac{p+1}{2}|\mathbf{P}|^{p-1}\epsilon + \frac{p-1}{2}|\mathbf{P}|^{p-3}\mathbf{P}^2\bar{\epsilon} \tag{2.16}$$

and

$$\frac{\bar{\epsilon}.F''(\mathbf{P}).\epsilon}{2} = \frac{p-1}{2}\epsilon^2\bar{\mathbf{P}}|\mathbf{P}|^{p-3} + (p-1)|\epsilon|^2\mathbf{P}|\mathbf{P}|^{p-3} + (p-1)\left(\frac{p}{2} - \frac{3}{2}\right)\text{Re}\,(\epsilon\bar{\mathbf{P}})^2\mathbf{P}|\mathbf{P}|^{p-5}.$$

In the case $p > 2$, set

$$2^+ = \min(3, \frac{p+2}{2}).$$

Remark that $2^+ < 2^*$ when $p > 2$ (where $2^* = \frac{2d}{d-2}$ is the critical exponent of the Sobolev injection). Then, from (2.15), we have

$$F(\mathbf{P} + \epsilon) = F(\mathbf{P}) + F'(\mathbf{P}).\epsilon + O(|\epsilon|^p) + O\left(\left|\frac{\epsilon}{\mathbf{P}}\right|^2 |\mathbf{P}|^p\right) \tag{2.17}$$

and

$$F(\mathbf{P} + \epsilon) = F(\mathbf{P}) + F'(\mathbf{P}).\epsilon + \frac{\bar{\epsilon}.F''(\mathbf{P}).\epsilon}{2} + O(|\epsilon|^p) + O\left(\left|\frac{\epsilon}{\mathbf{P}}\right|^{2^+} |\mathbf{P}|^p\right) \tag{2.18}$$

(note that for $\left|\frac{\epsilon}{\mathbf{P}}\right| \gg 1$ we have $F(\mathbf{P} + \epsilon) \sim F(\epsilon)$).

**Lemma 2** (Nonlinear interaction estimates). *For $|z| \gg 1, |v| \ll 1$, let*

$$H(z) = p\left[\int_{y \cdot \frac{z}{|z|} > -\frac{|z|}{2}} Q^{p-1}(y)\nabla Q(y)Q(y+z)dy + \int_{y \cdot \frac{z}{|z|} < -\frac{|z|}{2}} Q^{p-1}(y+z)\nabla Q(y)Q(y)dy\right]. \tag{2.19}$$

*Then the following estimates hold:*

$$\left|\langle G, e^{i\Gamma_1(y - z_1(s))}\nabla Q(y - z_1(s))\rangle - H(z)\right| \lesssim (|v|^2|z|^2 + |v|^2)|z|^{-\frac{d-1}{2}}e^{-|z|} + |z|^{-\frac{3(d-1)}{4}}e^{-\frac{3}{2}|z|} \tag{2.20}$$

*and*

$$\left|H(z) - C_p \frac{z}{|z|}|z|^{-\frac{d-1}{2}}e^{-|z|}\right| \lesssim |z|^{-\frac{d-1}{2}-1}e^{-|z|} \tag{2.21}$$

*where $C_p > 0$.*

**Remark 5.** *The estimate (2.21) on the leading order of the core part $H(z)$ of the projection $\langle G, [e^{i\Gamma_1}\nabla Q](y - z_1(s))\rangle$ is valid not only in the case $p > 2$ but also in the case $1 < p \leq 2$.*

*Proof of Lemma 2.* **step 1** Nonlinear interaction estimates. We prove the estimate (2.21) and in this step we will have $p > 1$. Consider

$$H(z) = p\int_{y \cdot \frac{z}{|z|} < -\frac{|z|}{2}} Q^{p-1}(y+z)\nabla Q(y)Q(y)dy + p\int_{y \cdot \frac{z}{|z|} > -\frac{|z|}{2}} Q^{p-1}(y)\nabla Q(y)Q(y+z)dy.$$

Recall that

$$Q(y)Q(y+z) \lesssim |z|^{-\frac{d-1}{2}}e^{-|z|}$$

$$Q(y)|\nabla Q(y+z)| \lesssim |z|^{-\frac{d-1}{2}}e^{-|z|}$$

then with $p > 2$, we have

$$\left|\int_{y \cdot \frac{z}{|z|} < -\frac{|z|}{2}} Q^{p-1}(y+z)\nabla Q(y)Q(y)dy\right| \lesssim e^{-\min(p-1, \frac{3}{2})|z|}$$

and with $1 < p \leq 2$, from the decay property of $Q$, we have for $\delta = \frac{p-1}{2}$

$$\left|\int_{y \cdot \frac{z}{|z|} < -\frac{|z|}{2}} Q^{p-1}(y+z)\nabla Q(y)Q(y)dy\right| \lesssim e^{-(p-1)|z|}\left|Q\left(\frac{|z|}{2}\right)\right|^{3-p-\delta}\int Q^\delta(y)dy$$

$$\lesssim e^{-\frac{p+3}{4}|z|}.$$

We claim that

$$\left| \int_{y \cdot \frac{z}{|z|} > -\frac{|z|}{2}} Q^{p-1}(y)\nabla Q(y)Q(y+z)dy - c_Q|z|^{-\frac{d-1}{2}}e^{-|z|} \int Q^{p-1}(y)\nabla Q(y)e^{-y \cdot \frac{z}{|z|}}dy \right| \tag{2.22}$$

$$\lesssim |z|^{-1-\frac{d-1}{2}}e^{-|z|}.$$

Indeed, let $0 < \theta < 1$ such that $p\theta > 1$. For $|y| \geq \theta|z|$, we have:

$$\left| \int_{\substack{|y| \geq \theta|z| \\ y \cdot \frac{z}{|z|} > -\frac{|z|}{2}}} Q^{p-1}(y)\nabla Q(y)Q(y+z)dy \right| \lesssim e^{-p\theta|z|} \left| \int Q(y+z)dy \right|$$

$$\lesssim e^{-p\theta|z|}.$$

For $|y| < \theta|z|$, as $Q(x) = q(|x|)$ and $|q(r) - c_Q r^{-\frac{d-1}{2}}e^{-r}| \lesssim r^{-\frac{d-1}{2}-1}e^{-r}$, we have:

$$\left| Q(y+z) - c_Q|y+z|^{-\frac{d-1}{2}}e^{-|y+z|} \right| \lesssim |y+z|^{-1-\frac{d-1}{2}}e^{-|y+z|}$$

$$\leq |1-\theta||z|^{-1-\frac{d-1}{2}}e^{-|z|}e^{|y|}.$$

Thus we get:

$$\left| \int_{\substack{|y| < \theta|z| \\ y \cdot \frac{z}{|z|} > -\frac{|z|}{2}}} Q^{p-1}(y)\nabla Q(y)\nabla Q(y+z)dy - c_Q \int_{\substack{|y| < \theta|z| \\ y \cdot \frac{z}{|z|} > -\frac{|z|}{2}}} Q^{p-1}(y)\nabla Q(y)|y+z|^{-\frac{d-1}{2}}e^{-|y+z|}dy \right|$$

$$\lesssim |z|^{-1-\frac{d-1}{2}}e^{-|z|}$$

since $\int Q^{p-1}(y)|\nabla Q(y)|e^{|y|}dy < +\infty$. On the other hand, $|y| < |z|$ implies

$$\left| |y+z|^{-k} - |z|^{-k} \right| \lesssim |z|^{-1-k}|y|$$

for any $k > 0$ and

$$\left| \frac{y+z}{|y+z|} - \frac{z}{|z|} \right| \lesssim |z|^{-1}|y|.$$

Moreover

$$\left| |y+z| - |z| - y \cdot \frac{z}{|z|} \right| \lesssim |z|^{-1}|y|^2$$

then

$$\left| e^{-|y+z|} - e^{-|z|-y \cdot \frac{z}{|z|}} \right| \lesssim |z|^{-1}|y|^2 e^{-|z|}e^{|y|}.$$

Thus we obtain that

$$\left| |y+z|^{-\frac{d-1}{2}}e^{-|y+z|} - |z|^{-\frac{d-1}{2}}e^{-|z|-y \cdot \frac{z}{|z|}} \right| \lesssim (1+|y|^2)|z|^{-1-\frac{d-1}{2}}e^{-|z|}e^{|y|}.$$

Therefore we have

$$\left| \int_{\substack{|y| < \theta|z| \\ y \cdot \frac{z}{|z|} > -\frac{|z|}{2}}} Q^{p-1}(y)\nabla Q(y)|y+z|^{-\frac{d-1}{2}}e^{-|y+z|}dy - c_Q|z|^{-\frac{d-1}{2}}e^{-|z|} \int_{\substack{|y| < \theta|z| \\ y \cdot \frac{z}{|z|} > -\frac{|z|}{2}}} Q^{p-1}(y)\nabla Q(y)e^{-y \cdot \frac{z}{|z|}}dy \right|$$

$$\lesssim |z|^{-1-\frac{d-1}{2}}e^{-|z|}.$$

Next we observe that

$$|z|^{-\frac{d-1}{2}}e^{-|z|}\int_{\substack{|y|\geq\theta|z| \\ y\cdot\frac{z}{|z|}>-\frac{|z|}{2}}}Q^{p-1}(y)\nabla Q(y)e^{-y\cdot\frac{z}{|z|}}dy \lesssim e^{-p\theta|z|}$$

and by (1.12)

$$\left|\int_{y\cdot\frac{z}{|z|}<-\frac{|z|}{2}}Q^{p-1}(y)\nabla Q(y)e^{-y\cdot\frac{z}{|z|}}dy\right| \lesssim e^{-\frac{p-1}{4}|z|}$$

which finish the proof of (2.22). Finally, in order to obtain (2.21) with $C_p = c_Q I_Q$, we use integration by parts

$$p\int Q^{p-1}(y)\nabla Q(y)e^{-y\cdot\frac{z}{|z|}}dy = \frac{z}{|z|}\int Q^p(y)e^{-y\cdot\frac{z}{|z|}}dy$$

and remark from the parity of the integral that

$$\int Q^p(y)e^{-y\cdot\frac{z}{|z|}}dy = \int Q^p(y)e^{-y_1}dy = I_Q.$$

**step 2** Error bound. From (2.15), we have the following estimates : if $y\cdot\frac{z}{|z|}>0$ then $|P_1|>|P_2|$

$$\left||\mathbf{P}|^{p-1}\mathbf{P} - |P_1|^{p-1}P_1 - |P_2|^{p-1}P_2 - \frac{p+1}{2}|P_1|^{p-1}P_2 - \frac{p-1}{2}|P_1|^{p-3}P_1^2\overline{P_2}\right| \lesssim |P_2|^2|P_1|^{p-2} \tag{2.23}$$

and if $y\cdot\frac{z}{|z|}<0$ then $|P_2|>|P_1|$

$$\left||\mathbf{P}|^{p-1}\mathbf{P} - |P_1|^{p-1}P_1 - |P_2|^{p-1}P_2 - \frac{p+1}{2}|P_2|^{p-1}P_1 - \frac{p-1}{2}|P_2|^{p-3}P_2^2\overline{P_1}\right| \lesssim |P_1|^2|P_2|^{p-2}. \tag{2.24}$$

We combine (2.23)–(2.24) to obtain, for all $y$,

$$\left||\mathbf{P}|^{p-1}\mathbf{P} - |P_1|^{p-1}P_1 - |P_2|^{p-1}P_2 - \left[\frac{p+1}{2}|P_1|^{p-1}P_2 + \frac{p-1}{2}|P_1|^{p-3}P_1^2\overline{P_2}\right].\mathbb{1}_{y\cdot\frac{z}{|z|}>0}\right.$$
$$\left.-\left[\frac{p+1}{2}|P_2|^{p-1}P_1 + \frac{p-1}{2}|P_2|^{p-3}P_2^2\overline{P_1}\right].\mathbb{1}_{y\cdot\frac{z}{|z|}<0}\right| \lesssim \min(|P_1|^2,|P_2|^2)\max(|P_1|^{p-2},|P_2|^{p-2}). \tag{2.25}$$

**step 3** Projection estimates. Since $\min(|P_1|^2,|P_2|^2)\leq |P_2|^{\frac{3}{2}}|P_1|^{\frac{1}{2}}$ and $\max(|P_1|^{p-2},|P_2|^{p-2})\leq |P_1|^{p-2}+|P_2|^{p-2}$, we have

$$\int Q^{\frac{3}{2}}(y-z)|\nabla Q(y)|Q^{\frac{1}{2}}(y)(Q^{p-2}(y) + Q^{p-2}(y+z))dy$$
$$\lesssim |z|^{-\frac{3(d-1)}{4}}e^{-\frac{3}{2}|z|}\int(Q^{p-2}(y) + Q^{p-2}(y+z))dy \lesssim |z|^{-\frac{3(d-1)}{4}}e^{-\frac{3}{2}|z|}$$

so we deduce from the error bound (2.25)

$$\left| \langle G, [e^{i\Gamma_1}\nabla Q](y - z_1(s)) \rangle - \left\langle \left[ \frac{p+1}{2}|P_1|^{p-1}P_2 + \frac{p-1}{2}|P_1|^{p-3}P_1^2\overline{P_2} \right] .\mathbb{1}_{y \cdot \frac{z}{|z|} > 0} \right.\right.$$
$$\left.\left. + \left[ \frac{p+1}{2}|P_2|^{p-1}P_1 + \frac{p-1}{2}|P_2|^{p-3}P_2^2\overline{P_1} \right] .\mathbb{1}_{y \cdot \frac{z}{|z|} < 0}, [e^{i\Gamma_1}\nabla Q](y - z_1(s)) \right\rangle \right| \lesssim |z|^{-\frac{3(d-1)}{4}}e^{-\frac{3}{2}|z|}.$$
(2.26)

Using a change of variables, we have

$$\langle |P_1|^{p-1}P_2\mathbb{1}_{y \cdot \frac{z}{|z|} > 0}, [e^{i\Gamma_1}\nabla Q](y - z_1(s)) \rangle$$
$$= \text{Re} \int_{y \cdot \frac{z}{|z|} > -\frac{|z|}{2}} Q^{p-1}(y)\nabla Q(y)Q(y - z_2 + z_1)e^{iv_2\cdot(y - z_2 + z_1) - iv_1\cdot y}dy$$
$$= \int_{y \cdot \frac{z}{|z|} > -\frac{|z|}{2}} Q^{p-1}(y)\nabla Q(y)Q(y + z)\cos(v_2\cdot(y + z) - v_1\cdot y)dy$$

with $z(s) = z_1(s) - z_2(s)$. Note that

$$|\cos(v_2\cdot(y + z) - v_1\cdot y) - 1| \lesssim |v|^2|z|^2 + |v|^2|y|^2$$

as the same method to prove (2.22), we get

$$\left| \langle |P_1|^{p-1}P_2.\mathbb{1}_{y \cdot \frac{z}{|z|} > 0}, [e^{i\Gamma_1}\nabla Q](y - z_1(s)) \rangle - \int_{y \cdot \frac{z}{|z|} > -\frac{|z|}{2}} Q^{p-1}(y)\nabla Q(y)Q(y + z)dy \right|$$
$$\lesssim (|v|^2|z|^2 + |v|^2)|z|^{-\frac{d-1}{2}}e^{-|z|}. \quad (2.27)$$

Similarly, for the other projections, we have

$$\left| \langle |P_1|^{p-3}P_1^2\overline{P_2}.\mathbb{1}_{y \cdot \frac{z}{|z|} > 0}, [e^{i\Gamma_1}\nabla Q](y - z_1(s)) \rangle - \int_{y \cdot \frac{z}{|z|} > -\frac{|z|}{2}} Q^{p-1}(y)\nabla Q(y)Q(y + z)dy \right|$$
$$\lesssim (|v|^2|z|^2 + |v|^2)|z|^{-\frac{d-1}{2}}e^{-|z|} \quad (2.28)$$

$$\left| \langle |P_1|^{p-3}P_1^2\overline{P_2}.\mathbb{1}_{y \cdot \frac{z}{|z|} < 0}, [e^{i\Gamma_1}\nabla Q](y - z_1(s)) \rangle - \int_{y \cdot \frac{z}{|z|} < -\frac{|z|}{2}} Q^{p-1}(y + z)\nabla Q(y)Q(y)dy \right|$$
$$\lesssim (|v|^2|z|^2 + |v|^2)|z|^{-\frac{d-1}{2}}e^{-|z|} \quad (2.29)$$

and finally

$$\langle |P_2|^{p-1}P_1.\mathbb{1}_{y \cdot \frac{z}{|z|} < 0}, [e^{i\Gamma_1}\nabla Q](y - z_1(s)) \rangle = \text{Re} \int_{y \cdot \frac{z}{|z|} < 0} Q^{p-1}(y - z_2(s))Q(y - z_1(s))\nabla Q(y - z_1(s))dy$$
$$= \int_{y \cdot \frac{z}{|z|} < -\frac{|z|}{2}} Q^{p-1}(y + z)\nabla Q(y)Q(y)dy. \quad (2.30)$$

From (2.26)–(2.30), we obtain the desired result (2.20).

$$\square$$

2.3. **Formal resolution and estimates of leading order.** From Lemma 1, we derive a simplified modulation system with forcing term and we determine one of its approximate solution that is relevant for the regime of Main Theorem. Formally, we have the following bounds (making this rigorous will be the goal of the bootstrap estimates in Sect. 3.2)

$$|\vec{m}_1| \lesssim |z|^{-\frac{d-1}{2}} e^{-|z|}, \tag{2.31}$$

from which we derive a simplified system ($\vec{m}_k$ is defined in (2.6)):

$$|\frac{\dot\lambda}{\lambda}| + |\dot z - 2v + \frac{\dot\lambda}{\lambda}z| \lesssim |z|^{-\frac{d-1}{2}} e^{-|z|}. \tag{2.32}$$

Furthermore, since we expect the interaction to be strong enough such that it will affect the main order of the modulation equations so by projecting $\mathcal{E}_{\mathbf{P}}$ onto the direction $e^{i\Gamma_1}\nabla Q(y - z_1(s))$, we obtain formally that

$$c_2 \dot v_1 \approx -\langle G, e^{i\Gamma_1}\nabla Q(y - z_1(s))\rangle \approx -H(z)$$

with $c_2 = \langle -yQ, \nabla Q \rangle > 0$. This remark suggests us to fix

$$\dot v = -\frac{2p}{c_2}\Big[\int Q^{p-1}(y)\nabla Q(y)Q(y+z)dy + \int Q^{p-1}(y+z)\nabla Q(y)Q(y)dy\Big] = -\frac{2}{c_2}H(z) \tag{2.33}$$

so $v(s)$ is completely determined by $z(s)$ and initial data $v^{in}$. In consequence, there are only three free parameters left $(\lambda, z, \gamma)$ corresponding to the scaling, translation and phase parameters which we will modulate to obtain orthogonality conditions (as shown below in Lemma 3). We use (2.21) to estimate the main order of $\dot v$

$$\left|\dot v + c\frac{z}{|z|}|z|^{-\frac{d-1}{2}} e^{-|z|}\right| \lesssim |z|^{-\frac{d-1}{2}-1} e^{-|z|} \tag{2.34}$$

with

$$c = \frac{2C_p}{c_2} = \frac{2c_Q I_Q}{c_2} > 0. \tag{2.35}$$

It can be checked that for some real functions $z_{mod}(s)$, $\lambda_{mod}(s)$, $v_{mod}(s)$ such that

$$\lambda_{mod}^{-1}(s) = 1, \quad v_{mod}(s) = s^{-1}, \quad z_{mod}^{-\frac{d-1}{2}} e^{-z_{mod}} = \frac{s^{-2}}{c} \tag{2.36}$$

then we have

$$z_{mod}(s) \sim 2\log(s), \qquad \dot v_{mod}(s) = -c z_{mod}^{-\frac{d-1}{2}}(s) e^{-z_{mod}(s)},$$
$$|\dot z_{mod}(s) - 2v_{mod}(s)| \lesssim s^{-1}\log^{-1}(s), \quad |\dot v_{mod}(s)| \lesssim s^{-2}. \tag{2.37}$$

Indeed, obviously $\dot v_{mod}(s) = -s^{-2} = -c z_{mod}^{-\frac{d-1}{2}}(s) e^{-z_{mod}(s)}$ and by differentiating the equation of $z_{mod}$, we get

$$-\dot z_{mod} z_{mod}^{-\frac{d-1}{2}} e^{-z_{mod}} - \frac{d-1}{2}\dot z_{mod} z_{mod}^{-\frac{d-1}{2}-1} e^{-z_{mod}} = -2\frac{s^{-3}}{c}$$

(in the case $d-1 = 0$, $-\dot z_{mod} e^{-z_{mod}} = -2\frac{s^{-3}}{c}$) so $|\dot z_{mod} - 2s^{-1}| \lesssim s^{-1}\log^{-1}(s)$ thus we can deduce $|\dot z_{mod}(s) - 2v_{mod}(s)| \lesssim s^{-1}\log^{-1}(s)$. The above estimates suggest that (2.36) is close to the first order asymptotics as $s \to +\infty$ for some particular solutions of (2.32) and matches the regime in Main Theorem.

## 3. Modulation and backward uniform estimates

Let $(\lambda^{in}, z^{in}, v^{in}) \in (0, +\infty) \times (0, +\infty) \times \mathbb{R}$ to be chosen with $|z^{in}| \gg 1$, $|v^{in}| \ll 1$, $T_{mod} > 0$ and $(\vec{e_1}, ..., \vec{e_d})$ standard basis of $\mathbb{R}^d$. Recall that in this section $p > 2$. Let $u(t, x)$ be the backward solution of (NLS) with initial data

$$u(T_{mod}, x) = \frac{1}{(\lambda^{in})^{\frac{2}{p-1}}} \mathbf{P}^{in} \left( \frac{x}{\lambda^{in}} \right) \quad \text{where} \quad \mathbf{P}^{in}(y) = \mathbf{P}(y; (z^{in}\vec{e_1}, v^{in})) \tag{3.1}$$

on some time interval including $T_{mod}$. Note that the equation (NLS) is invariant by rotation. In particular, if a solution of (NLS) is invariant by the symmetry $\tau : x \mapsto -x$ at some time, then it is invariant by the symmetry at any time.

### 3.1. Decomposition of $u(t)$.
We will state a standard modulation result with the same idea as in Lemma 3 of [15] or Lemma 2 of [22]. The choice of the special orthogonality conditions (3.5) is related to the generalized null space of the linearized equation around $Q$ in (1.14) and to the coercivity property (1.15) in sub-critical cases. See the proof of Lemma 5 for a technical justification of these choices. For $s^{in} \gg 1$ fixed.

**Lemma 3** (Modulation of the approximate solution). *Let $u(t, x)$ a solution invariant by $\tau$ on an interval $[T, T_{mod}]$ satisfying $u(T_{mod}, x) \in H^2(\mathbb{R}^d)$ and*

$$\left\| e^{-i\gamma^{in}} (\lambda^{in})^{\frac{2}{p-1}} u(T_{mod}, \lambda^{in} y) - \mathbf{P}(y; (z^{in}\vec{e_1}, v^{in})) \right\|_{H^1} \ll 1$$

*for $\mathbf{P}(s, y) = \mathbf{P}(y; (z(s), v(s)))$ as defined in (2.5). Then there exist a unique $\mathcal{C}^1$ function on an open interval $I \ni s^{in}$*

$$\vec{q}(s) = (\lambda, z, \gamma, v) : I \to (0, +\infty) \times \mathbb{R}^d \times \mathbb{R} \times \mathbb{R}^d,$$

*with $\vec{q}(s^{in}) = (\lambda^{in}, z^{in}\vec{e_1}, \gamma^{in}, v^{in})$ and a rescaling time function*

$$t(s) = T_{mod} - \int_s^{s^{in}} \lambda^2(\tau) d\tau \tag{3.2}$$

*such that $u(t, x)$ decomposes as follows*

$$u(t(s), x) = \frac{e^{i\gamma(s)}}{\lambda^{\frac{2}{p-1}}(s)} (\mathbf{P} + \epsilon)(s, y), \quad y = \frac{x}{\lambda(s)} \tag{3.3}$$

*where by setting*

$$\epsilon(s, y) = \left[ e^{i\Gamma_1} \eta_1 \right](s, y - z_1), \qquad \Gamma_k(s, y) = v_k(s) \cdot y, \tag{3.4}$$

*if initially $\langle \eta_1(s^{in}), Q \rangle = \langle \eta_1(s^{in}), yQ \rangle = \langle \eta_1(s^{in}), i\Lambda Q \rangle = 0$, the decomposition satisfies orthogonality conditions*

$$\langle \eta_1(s), Q \rangle = \langle \eta_1(s), yQ \rangle = \langle \eta_1(s), i\Lambda Q \rangle = 0 \tag{3.5}$$

*and the extra relation*

$$\dot{v}(s) = -\frac{2}{c_2} H(z(s)). \tag{3.6}$$

*Moreover, $\epsilon$ is also invariant by the symmetry $\tau$.*

*Proof of Lemma 3.* **step 1** Orthogonality conditions. Remark that we can go easily from the rescaled time $s$ to $t$ and conversely

$$s = s(t) = s^{in} - \int_t^{T_{mod}} \frac{d\tau}{\lambda^2(\tau)} \tag{3.7}$$

with $T_{mod} = t(s^{in})$. Denote

$$\mathbf{P}(s,y) = \left[ e^{i\Gamma_1} \mathbf{P}_1 \right](s, y - z_1), \qquad \mathcal{E}_{\mathbf{P}}(s,y) = \left[ e^{i\Gamma_1} \mathcal{E}_{\mathbf{P}_1} \right](s, y - z_1)$$

$$G(s,y) = [e^{i\Gamma_1} G_1](s, y - z_1)$$

where $G = |\mathbf{P}|^{p-1}\mathbf{P} - |P_1|^{p-1}P_1 - |P_2|^{p-1}P_2$. Let $w = \mathbf{P} + \epsilon$ as in (2.1). It follows from the equation of $w$ (2.2) and the equation of $\mathbf{P}$ (2.7) that

$$i\dot{\epsilon} + \Delta\epsilon - \epsilon + \left( |\mathbf{P} + \epsilon|^{p-1}(\mathbf{P} + \epsilon) - |\mathbf{P}|^{p-1}\mathbf{P} \right) - i\frac{\dot{\lambda}}{\lambda}\Lambda\epsilon + (1 - \dot{\gamma})\epsilon + \mathcal{E}_{\mathbf{P}} = 0. \tag{3.8}$$

We rewrite the equation of $\epsilon$ into the following equation for $\eta_1$ (see also the proof of Lemma 1)

$$i\dot{\eta}_1 + \Delta\eta_1 - \eta_1 + (|\mathbf{P}_1 + \eta_1|^{p-1}(\mathbf{P}_1 + \eta_1) - |\mathbf{P}_1|^{p-1}\mathbf{P}_1) + \vec{m}_1 \cdot \vec{\mathrm{M}}\eta_1 + \mathcal{E}_{\mathbf{P}_1} = 0. \tag{3.9}$$

Thus, for $A(y), B(y) \in \mathcal{Y}$, we get

$$\frac{d}{ds}\langle \eta_1, A + iB \rangle = -\langle \Delta\eta_1 - \eta_1 + (|\mathbf{P}_1 + \eta_1|^{p-1}(\mathbf{P}_1 + \eta_1) - |\mathbf{P}_1|^{p-1}\mathbf{P}_1)$$
$$+ \vec{m}_1 \cdot \vec{\mathrm{M}}\eta_1 + \mathcal{E}_{\mathbf{P}_1}, iA - B \rangle.$$

Choose $A = Q, B = 0$ and $A = yQ, B = 0$ and $A = 0, B = \Lambda Q$ then the conditions

$$\frac{d}{ds}\langle \eta_1(s), Q \rangle = \frac{d}{ds}\langle \eta_1(s), yQ \rangle = \frac{d}{ds}\langle \eta_1(s), i\Lambda Q \rangle = 0$$

are equivalent to

$$\begin{cases} \left\langle \Delta\eta_1 - \eta_1 + (|\mathbf{P}_1 + \eta_1|^{p-1}(\mathbf{P}_1 + \eta_1) - |\mathbf{P}_1|^{p-1}\mathbf{P}_1) + \vec{m}_1 \cdot \vec{\mathrm{M}}\eta_1 + \mathcal{E}_{\mathbf{P}_1}, iQ \right\rangle = 0 \\ \left\langle \Delta\eta_1 - \eta_1 + (|\mathbf{P}_1 + \eta_1|^{p-1}(\mathbf{P}_1 + \eta_1) - |\mathbf{P}_1|^{p-1}\mathbf{P}_1) + \vec{m}_1 \cdot \vec{\mathrm{M}}\eta_1 + \mathcal{E}_{\mathbf{P}_1}, iyQ \right\rangle = 0 \\ \left\langle \Delta\eta_1 - \eta_1 + (|\mathbf{P}_1 + \eta_1|^{p-1}(\mathbf{P}_1 + \eta_1) - |\mathbf{P}_1|^{p-1}\mathbf{P}_1) + \vec{m}_1 \cdot \vec{\mathrm{M}}\eta_1 + \mathcal{E}_{\mathbf{P}_1}, -\Lambda Q \right\rangle = 0. \end{cases}$$

We claim that the above system is equivalent to an autonomous system of ordinary differential equations on $(\theta(s), z(s), \gamma(s), v(s), t(s))$ where $\theta(s) = \ln(\lambda(s))$. Indeed, remark that

$$\epsilon(s,y) = e^{\frac{2}{p-1}\theta(s)} u(t(s), e^{\theta(s)}y) - \mathbf{P}(y; (z(s), v(s))) \tag{3.10}$$

and the expression of $\mathcal{E}_{\mathbf{P}_1}$ (from (2.7)–(2.8))

$$\mathcal{E}_{\mathbf{P}_1} = [\vec{m}_1 \cdot \vec{\mathrm{M}}Q](y) + [e^{i(\Gamma_2(y+z) - \Gamma_1(y))}\vec{m}_2 \cdot \vec{\mathrm{M}}Q](y+z) + G_1$$

then we get

$$\begin{cases} \langle \vec{m}_1 \cdot \vec{\mathrm{M}}Q, iQ \rangle + \langle e^{i(\Gamma_2(y+z) - \Gamma_1(y))}\vec{m}_2 \cdot \vec{\mathrm{M}}Q(y+z), iQ \rangle + \langle \vec{m}_1 \cdot \vec{\mathrm{M}}\eta_1, iQ \rangle = \mathcal{F}_1(\theta, z, \gamma, v, t) \\ \langle \vec{m}_1 \cdot \vec{\mathrm{M}}Q, iyQ \rangle + \langle e^{i(\Gamma_2(y+z) - \Gamma_1(y))}\vec{m}_2 \cdot \vec{\mathrm{M}}Q(y+z), iyQ \rangle + \langle \vec{m}_1 \cdot \vec{\mathrm{M}}\eta_1, iyQ \rangle = \mathcal{F}_2(\theta, z, \gamma, v, t) \\ \langle \vec{m}_1 \cdot \vec{\mathrm{M}}Q, -\Lambda Q \rangle + \langle e^{i(\Gamma_2(y+z) - \Gamma_1(y))}\vec{m}_2 \cdot \vec{\mathrm{M}}Q(y+z), -\Lambda Q \rangle + \langle \vec{m}_1 \cdot \vec{\mathrm{M}}\eta_1, -\Lambda Q \rangle = \mathcal{F}_3(\theta, z, \gamma, v, t) \end{cases} \tag{3.11}$$

with

$$\mathcal{F}_1(\theta, z, \gamma, v, t) = -\left\langle \Delta\eta_1 - \eta_1 + (|\mathbf{P}_1 + \eta_1|^{p-1}(\mathbf{P}_1 + \eta_1) - |\mathbf{P}_1|^{p-1}\mathbf{P}_1) + G_1, iQ \right\rangle$$

$$\mathcal{F}_2(\theta, z, \gamma, v, t) = -\left\langle \Delta\eta_1 - \eta_1 + (|\mathbf{P}_1 + \eta_1|^{p-1}(\mathbf{P}_1 + \eta_1) - |\mathbf{P}_1|^{p-1}\mathbf{P}_1) + G_1, iyQ \right\rangle$$

$$\mathcal{F}_3(\theta, z, \gamma, v, t) = -\left\langle \Delta\eta_1 - \eta_1 + (|\mathbf{P}_1 + \eta_1|^{p-1}(\mathbf{P}_1 + \eta_1) - |\mathbf{P}_1|^{p-1}\mathbf{P}_1) + G_1, -\Lambda Q \right\rangle.$$

Note that $\mathcal{F}_1, \mathcal{F}_2, \mathcal{F}_3$ are $\mathcal{C}^1$ functions. Indeed, if we replace $\eta_1$ by the expression (3.10) and its definition, it is clear that any term not containing $u$ is continuously differentiable. For terms concerning $u(t, x)$, by integration by parts and chain rule, we show how to prove that typical terms, integrals of the form

$$\frac{d}{dt} \mathrm{Re} \left( \int u(t, x) A(x) dx \right), \quad \frac{d}{dt} \mathrm{Re} \left( \int |u(t, x)|^{p-1} u(t, x) A(x) dx \right)$$

for $A(x)$ some complex functions such that $\mathrm{Re}\, A(x)$, $\mathrm{Im}\, A(x) \in \mathcal{Y}$, are continuous. We have

$$\frac{d}{dt} \mathrm{Re} \left( \int u(t, x) A(x) dx \right) = -\mathrm{Im} \left( \int u(t, x) \Delta A(x) dx \right) - \mathrm{Im} \left( \int |u(t, x)|^{p-1} u(t, x) A(x) dx \right) \tag{3.12}$$

and

$$\frac{d}{dt} \mathrm{Re} \left( \int |u(t, x)|^{p-1} u(t, x) A(x) dx \right) = p\, \mathrm{Re} \left( \int \partial_t u(t, x) |u(t, x)|^{p-1} A(x) dx \right) =$$

$$- p\, \mathrm{Im} \left( \int \Delta u(t, x) |u(t, x)|^{p-1} A(x) dx \right) - p\, \mathrm{Im} \left( \int |u(t, x)|^{2p-2} u(t, x) A(x) dx \right). \tag{3.13}$$

Recall the persistence of $H^2$ regularity for (NLS) equation (see Theorem 5.3.1 in [2]), since $u(T_{mod}, x) \in H^2(\mathbb{R}^d)$ then $u \in \mathcal{C}^1([0, T_{mod}], L^2(\mathbb{R}^d)) \bigcap \mathcal{C}([0, T_{mod}], H^2(\mathbb{R}^d))$. By Sobolev's injection ($\frac{d+6}{d-2} < \frac{2d}{d-4}$), we have $u \in \mathcal{C}([0, T_{mod}], L^{2p-1}(\mathbb{R}^d))$ thus the right-hand sides of (3.12), (3.13) are well-defined and continuous. Therefore, in particular, since initially

$$\langle \eta_1(s^{in}), Q \rangle = \langle \eta_1(s^{in}), yQ \rangle = \langle \eta_1(s^{in}), i\Lambda Q \rangle = 0,$$

the decomposition $(\vec{q}, \epsilon)$ will satisfy (3.5) if (3.11) holds.

**step 2** System of ODEs. We look a decomposition $(\vec{q}, \epsilon)$ of $u(t)$ and a rescaling time $t(s)$ such that

$$\begin{cases} \langle \vec{m}_1 \cdot \vec{\mathrm{M}} Q, iQ \rangle + \langle e^{i(\Gamma_2(y+z) - \Gamma_1(y))} \vec{m}_2 \cdot \vec{\mathrm{M}} Q(y+z), iQ \rangle + \langle \vec{m}_1 \cdot \vec{\mathrm{M}} \eta_1, iQ \rangle = \mathcal{F}_1(\theta, z, \gamma, v, t) \\ \langle \vec{m}_1 \cdot \vec{\mathrm{M}} Q, iyQ \rangle + \langle e^{i(\Gamma_2(y+z) - \Gamma_1(y))} \vec{m}_2 \cdot \vec{\mathrm{M}} Q(y+z), iyQ \rangle + \langle \vec{m}_1 \cdot \vec{\mathrm{M}} \eta_1, iyQ \rangle = \mathcal{F}_3(\theta, z, \gamma, v, t) \\ \langle \vec{m}_1 \cdot \vec{\mathrm{M}} Q, -\Lambda Q \rangle + \langle e^{i(\Gamma_2(y+z) - \Gamma_1(y))} \vec{m}_2 \cdot \vec{\mathrm{M}} Q(y+z), -\Lambda Q \rangle + \langle \vec{m}_1 \cdot \vec{\mathrm{M}} \eta_1, -\Lambda Q \rangle = \mathcal{F}_2(\theta, z, \gamma, v, t) \\ \dot{v} = -\frac{2}{c_2} H(z) \\ \dot{t}(s) = \lambda^2(s). \end{cases} \tag{3.14}$$

On the one hand, we calculate

$$\langle \vec{m}_1 \cdot \vec{\mathrm{M}} Q, iQ \rangle = (\tfrac{\dot{\lambda}}{\lambda}) \langle -i\Lambda Q, iQ \rangle = -c_1(\tfrac{\dot{\lambda}}{\lambda})$$

$$\langle \vec{m}_1 \cdot \vec{\mathrm{M}} Q, iyQ \rangle = (\dot{z} - 2v + \tfrac{\dot{\lambda}}{\lambda} z) \langle -i\nabla Q, iyQ \rangle = c_2(\dot{z} - 2v + \tfrac{\dot{\lambda}}{\lambda} z)$$

$$\langle \vec{m}_1 \cdot \vec{\mathrm{M}} Q, -\Lambda Q \rangle = c_1(\dot{\gamma} - 1 + |v|^2 - \tfrac{\dot{\lambda}}{\lambda}(v \cdot z) - (v \cdot \dot{z}))$$

with $c_1 = \langle \Lambda Q, Q \rangle, c_2 = \langle -\nabla Q, yQ \rangle$ non-zero. On the other hand, there exist a matrix $\mathcal{M}(\theta, z, \gamma, v, t) = (m_{ij})_{5 \times 5}$ and $\vec{\mathcal{G}}(\theta, z, \gamma, v, t)$ such that

$$
\begin{pmatrix}
\langle e^{i(\Gamma_2(y+z) - \Gamma_1(y))} \vec{m}_2 \cdot \vec{\mathbb{M}} Q(y+z), iQ \rangle + \langle \vec{m}_1 \cdot \vec{\mathbb{M}} \eta_1, iQ \rangle \\
\langle e^{i(\Gamma_2(y+z) - \Gamma_1(y))} \vec{m}_2 \cdot \vec{\mathbb{M}} Q(y+z), iyQ \rangle + \langle \vec{m}_1 \cdot \vec{\mathbb{M}} \eta_1, iyQ \rangle \\
\langle e^{i(\Gamma_2(y+z) - \Gamma_1(y))} \vec{m}_2 \cdot \vec{\mathbb{M}} Q(y+z), -\Lambda Q \rangle + \langle \vec{m}_1 \cdot \vec{\mathbb{M}} \eta_1, -\Lambda Q \rangle \\
0 \\
0
\end{pmatrix} = (\dot{\theta}, \dot{z}, \dot{\gamma}, \dot{v}, \dot{t}) \mathcal{M}(\theta, z, \gamma, v, t)
$$

$$
+ \vec{\mathcal{G}}(\theta, z, \gamma, v, t) \quad (3.15)
$$

where all entries of $\mathcal{M}(\theta, z, \gamma, v, t)$ are small $|m_{ij}| \ll 1$ as $z^{in} \gg 1$ and $||\epsilon(s^{in})||_{H^1} \ll 1$ (from hypothesis). Then the system (3.14) can be rewritten as an autonomous system

$$
(\dot{\theta}, \dot{z}, \dot{\gamma}, \dot{v}, \dot{t}) \mathcal{A}(\theta, z, \gamma, v, t) + (\dot{\theta}, \dot{z}, \dot{\gamma}, \dot{v}, \dot{t}) \mathcal{M}(\theta, z, \gamma, v, t) = \vec{\mathcal{H}}(\theta, z, \gamma, v, t) \quad (3.16)
$$

where

$$
\vec{\mathcal{H}}(\theta, z, \gamma, v, t) = \begin{pmatrix}
\mathcal{F}_1(\theta, z, \gamma, v, t) \\
\mathcal{F}_2(\theta, z, \gamma, v, t) + 2c_2 v \\
\mathcal{F}_3(\theta, z, \gamma, v, t) + c_1 - c_1 |v|^2 \\
-\frac{2}{c_2} H(z) \\
e^{2\theta}
\end{pmatrix} - \vec{\mathcal{G}}(\theta, z, \gamma, v, t)
$$

and the matrix $\mathcal{A}$ is given by

$$
\mathcal{A} = \begin{pmatrix}
-c_1 & c_2 z & c_1(v \cdot z) & 0 & 0 \\
0 & c_2 & c_1 v & 0 & 0 \\
0 & 0 & c_1 & 0 & 0 \\
0 & 0 & 0 & 1 & 0 \\
0 & 0 & 0 & 0 & 1
\end{pmatrix}.
$$

Therefore the perturbed matrix $(\mathcal{A} + \mathcal{M})(\theta, z, \gamma, v, t)$ is invertible $(\det \mathcal{A} = -c_1^2 c_2 < 0)$. As same as the way to deal with $\mathcal{F}$, one can check that $\mathcal{M}, \vec{\mathcal{G}}$ are continuously differentiable thus so are entries of $(\mathcal{A} + \mathcal{M})^{-1}$ and $\vec{\mathcal{H}}$. Therefore,

$$
\mathcal{R}(\theta, z, \gamma, v, t) = [(\mathcal{A} + \mathcal{M})^{-1} \cdot \vec{\mathcal{H}}](\theta, z, \gamma, v, t)
$$

satisfies the hypothesis of Cauchy-Lipschitz theorem and the system of ODEs

$$
(\dot{\theta}, \dot{z}, \dot{\gamma}, \dot{v}, \dot{t}) = \mathcal{R}(\theta, z, \gamma, v, t) \quad (3.17)
$$

admits a unique solution $(\theta(s), z(s), \gamma(s), v(s), t(s))$ to the initial value problem. We obtain the decomposition $(\lambda(s), z(s), \gamma(s), v(s))$ of $u(t)$ and the renormalization of time $t(s)$.

$\square$

Observe from (3.1) that the initial data

$$
w(s^{in}) = \mathbf{P}^{in}(y; (z^{in} \vec{e_1}, v^{in})), \quad \lambda(s^{in}) = \lambda^{in}, \quad \gamma(s^{in}) = 0,
$$
$$
z(s^{in}) = z^{in} \vec{e_1}, \quad v(s^{in}) = v^{in}, \quad \epsilon(s^{in}) \equiv 0 \quad (3.18)
$$

and $u(T_{mod}, x)$ satisfy the hypothesis of Lemma 3.

**Proposition 4** (Uniform backwards estimates for $p > 2$)**.** *There exists $s_0 > 10$ satisfying the following condition: for all $s^{in} > s_0$, there is a choice of initial parameters $(\lambda^{in}, z^{in}, v^{in})$ with*

$$\left| c^{-\frac{1}{2}} (z^{in})^{\frac{d-1}{4}} e^{\frac{1}{2} z^{in}} - s^{in} \right| < s^{in} \log^{-\frac{1}{2}}(s^{in}), \quad z^{in} > 0,$$

$$\lambda^{in} = 1, \quad v^{in} = c^{\frac{1}{2}} (z^{in})^{-\frac{d-1}{4}} e^{-\frac{1}{2} z^{in}} \cdot \vec{e_1},$$ (3.19)

*such that the solution $u$ of (NLS) corresponding to (3.1) exists. Moreover, the decomposition of $u$ given by Lemma 3 on the rescaled interval of time $[s_0, s^{in}]$*

$$u(s, x) = \frac{e^{i\gamma(s)}}{\lambda^{\frac{2}{p-1}}(s)} (\mathbf{P} + \epsilon)(s, y), \quad y = \frac{x}{\lambda(s)}, \quad dt = \lambda^2(s) ds$$

*verifies the uniform estimates for all $s \in [s_0, s^{in}]$*

$$| |z(s)| - 2\log(s)| \lesssim \log(\log(s)), \quad \left| \lambda^{-1}(s) - 1 \right| \lesssim s^{-1},$$

$$|v(s)| \lesssim s^{-1}, \quad \|\epsilon(s)\|_{H^1} \lesssim s^{-1}, \quad \left| |z(s)|^{\frac{d-1}{2}} e^{|z(s)|} - cs^2 \right| \lesssim s^2 \log^{-\frac{1}{2}}(s).$$ (3.20)

**Remark 6.** *The key point in Proposition 4 is that $s_0$ and the constants in (3.20) are independent of $s^{in}$ as $s^{in} \to +\infty$. Observe that the estimates (3.20) match the discussion in Sect. 2.3. The decomposition in Lemma 3 is only local but the estimates in (3.20) guarantee the global existence of the decomposition. The choice of $v^{in}$ is direct while the choice of $z^{in}$ is based on a contradiction argument and a topological constraint.*

The next subsection is devoted to the proof of Proposition 4 containing several technical steps. The proof relies on a bootstrap argument, integration of the differential system of geometrical parameters and energy estimates.

### 3.2. **Proof of Proposition 4.**

3.2.1. *Bootstrap bounds.* We shall consider the following bootstrap estimates

$$\left| c^{-\frac{1}{2}} |z|^{\frac{d-1}{4}} e^{\frac{1}{2}|z|} - s \right| \le s \log^{-\frac{1}{2}}(s),$$

$$\|\epsilon(s)\|_{H^1} \le C^* s^{-1}$$ (3.21)

with $C^* > 1$ to be chosen large enough. Note that the estimate on $z$ and the estimate (2.34) of $\dot{v}$ imply that, for $s$ large

$$\left| |z| - 2\log(s) \right| \lesssim \log(\log(s)), \quad \left| |\dot{v}| - s^{-2} \right| \lesssim s^{-2} \log^{-1}(s), \quad \left| |v| - s^{-1} \right| \lesssim s^{-1} \log^{-1}(s) \quad (3.22)$$

where the last inequality is obtained by integrating the second one with the choice of initial data $v^{in}$ in (3.19). Next, we define

$$s^* = \inf\{\tau \in [s_0, s^{in}]; \ (3.21) \text{ holds on } [\tau, s^{in}]\}.$$ (3.23)

3.2.2. *Control of the modulation equations.*

**Lemma 5** (Pointwise control of the modulation equations and the error)**.** *The following estimates hold on $[s^*, s^{in}]$.*

$$|\vec{m}_1(s)| \lesssim (C^*)^2 s^{-2}.$$ (3.24)

$$|\langle \eta_1(s), i\nabla Q \rangle| \lesssim (C^*)^2 s^{-2},$$ (3.25)

$$|\dot{z} - 2v| \lesssim s^{-1} \log^{-1}(s).$$ (3.26)

*Moreover, for all $s \in [s^*, s^{in}]$, for all $y \in \mathbb{R}^2$,*

$$|\mathcal{E}_{\mathbf{P}}(s, y)| \lesssim (C^*)^2 s^{-2} \sum_{k=1}^{2} Q(y - z_k(s)) + |G(s, y)|. \tag{3.27}$$

*Proof of Lemma 5.* Since $\epsilon(s^{in}) \equiv 0$, we may define

$$s^{**} = \inf\{s \in [s^*, s^{in}]; \ |\langle \eta_1(\tau), i\nabla Q\rangle| \leq C^{**}\tau^{-2} \text{ holds on } [s, s^{in}]\},$$

for some constant $C^{**} > 0$ to be chosen large enough. We work on the interval $[s^{**}, s^{in}]$. Recall equation for $\eta_1$ as below

$$i\dot{\eta}_1 + \Delta\eta_1 - \eta_1 + (|\mathbf{P}_1 + \eta_1|^{p-1}(\mathbf{P}_1 + \eta_1) - |\mathbf{P}_1|^{p-1}\mathbf{P}_1) + \vec{m}_1 \cdot \vec{\mathrm{M}}\eta_1 + \mathcal{E}_{\mathbf{P}_1} = 0.$$

Let $A(y)$ and $B(y)$ be two real-valued functions in $\mathcal{Y}$. We claim the following estimate on $[s^{**}, s^{in}]$

$$\left| \frac{d}{ds}\langle \eta_1, A + iB\rangle - \left[\langle \eta_1, iL_-A - L_+B\rangle - \langle \vec{m}_1 \cdot \vec{\mathrm{M}}Q, iA - B\rangle\right] \right| \lesssim (C^*)^2 s^{-2} + s^{-1}|\vec{m}_1|. \tag{3.28}$$

We compute from (3.9),

$$\frac{d}{ds}\langle \eta_1, A + iB\rangle = \langle \dot{\eta}_1, A + iB\rangle = \langle i\dot{\eta}_1, iA - B\rangle$$

$$= \langle -\Delta\eta_1 + \eta_1 - (\frac{p+1}{2}Q^{p-1}\eta_1 + \frac{p-1}{2}Q^{p-1}\overline{\eta}_1), iA - B\rangle$$

$$- \langle |\mathbf{P}_1 + \eta_1|^{p-1}(\mathbf{P}_1 + \eta_1) - |\mathbf{P}_1|^{p-1}\mathbf{P}_1 - \frac{p+1}{2}Q^{p-1}\eta_1 - \frac{p-1}{2}Q^{p-1}\overline{\eta}_1, iA - B\rangle$$

$$- \langle \vec{m}_1 \cdot \vec{\mathrm{M}}\eta_1, iA - B\rangle - \langle \mathcal{E}_{\mathbf{P}_1}, iA - B\rangle.$$

First, since $A$ and $B$ are real-valued, we have

$$\langle -\Delta\eta_1 + \eta_1 - (\frac{p+1}{2}Q^{p-1}\eta_1 + \frac{p-1}{2}Q^{p-1}\overline{\eta}_1), iA - B\rangle = \langle \eta_1, iL_-A - L_+B\rangle.$$

Second, recall the expression of $\mathbf{P}_1$

$$\mathbf{P}_1 = Q(y) + e^{i(\Gamma_2(y-(z_2-z_1))-\Gamma_1(y))}Q(y - (z_2 - z_1)).$$

By the expansion in (2.17), we can deduce the first order and the error of

$$|\mathbf{P}_1 + \eta_1|^{p-1}(\mathbf{P}_1 + \eta_1) - |\mathbf{P}_1|^{p-1}\mathbf{P}_1 - \frac{p+1}{2}Q^{p-1}\eta_1 - \frac{p-1}{2}Q^{p-1}\overline{\eta}_1$$

$$= \frac{p+1}{2}(|\mathbf{P}_1|^{p-1} - Q^{p-1})\eta_1 + \frac{p-1}{2}(|\mathbf{P}_1|^{p-3}\mathbf{P}_1^2 - Q^{p-1})\overline{\eta}_1 + O\big(\big|\frac{\eta_1}{\mathbf{P}_1}\big|^2 |\mathbf{P}_1|^p\big) + O(|\eta_1|^p).$$

Therefore, by (3.21)–(3.22) for some $q > 0$,

$$|\langle(|\mathbf{P}_1|^{p-1} - Q^{p-1})\eta_1, (iA - B)\rangle| + |\langle(|\mathbf{P}_1|^{p-3}\mathbf{P}_1^2 - Q^{p-1})\overline{\eta}_1, (iA - B)\rangle|$$

$$\lesssim |z|^q e^{-|z|}\|\eta_1\|_{L^2} \lesssim C^* s^{-3} \log^q(s).$$

Using Cauchy-Schwartz and Gagliardo-Nirenberg inequalities (as $p > 2$)

$$\left\langle \big|\frac{\eta_1}{\mathbf{P}_1}\big|^2 |\mathbf{P}_1|^p, (iA - B)\right\rangle \lesssim \|\epsilon\|_{L^2}^2 \lesssim (C^*)^2 s^{-2},$$

$$\langle |\eta_1|^p, (iA - B)\rangle \lesssim \|\epsilon\|_{H^1}^p \lesssim (C^*)^2 s^{-2}.$$

Therefore

$$\left|\langle|\mathbf{P}_1+\eta_1|^{p-1}(\mathbf{P}_1+\eta_1)-|\mathbf{P}_1|^{p-1}\mathbf{P}_1-\frac{p+1}{2}Q^{p-1}\eta_1-\frac{p-1}{2}Q^{p-1}\overline{\eta}_1,iA-B\rangle\right|\lesssim(C^*)^2s^{-2}. \tag{3.29}$$

Next, using (3.21)–(3.22), we obtain

$$|\langle\vec{m}_1\cdot\vec{\mathrm{M}}\eta_1,iA-B\rangle|\lesssim C^*s^{-1}|\vec{m}_1(s)|.$$

Finally, we need to prove following estimate

$$\left|\langle\mathcal{E}_{\mathbf{P}_1},iA-B\rangle-\langle\vec{m}_1\cdot\vec{\mathrm{M}}Q,iA-B\rangle\right|\lesssim s^{-2}+s^{-1}|\vec{m}_1|. \tag{3.30}$$

Indeed, recall that we have

$$\mathcal{E}_{\mathbf{P}_1}=[\vec{m}_1\cdot\vec{\mathrm{M}}Q](y)+[e^{i(\Gamma_2(y-(z_2-z_1))-\Gamma_1(y))}\vec{m}_2\cdot\vec{\mathrm{M}}Q](y-(z_2-z_1))+G_1.$$

From (2.14) and (3.21)–(3.22),

$$|\langle G_1,iA-B\rangle|\lesssim\|G\|_{L^\infty}\lesssim|z|^{-\frac{d-1}{2}}e^{-|z|}\lesssim s^{-2}.$$

Since $A,B\in\mathcal{Y}$, we have

$$|\langle e^{i(\Gamma_2(y-(z_2-z_1))-\Gamma_1(y))}(\vec{m}_2\cdot\vec{\mathrm{M}}Q(.-(z_2-z_1))),iA-B\rangle|\lesssim s^{-1}|\vec{m}_1|,$$

so the proof of (3.30) is complete.

We now use (3.28) to control the modulation vector $\vec{m}_1$. Note that $\eta_1$ satisfies the orthogonality conditions (3.5).

$\underline{\langle\eta_1,Q\rangle=0.}$ Let $A=Q$ and $B=0$. Since $L_-Q=0$ and $\langle\vec{m}_1\cdot\vec{\mathrm{M}}Q,iQ\rangle=-c_1(\frac{\dot{\lambda}}{\lambda})$, we obtain

$$\left|\frac{\dot{\lambda}}{\lambda}\right|\lesssim(C^*)^2s^{-2}+s^{-1}|\vec{m}_1|. \tag{3.31}$$

$\underline{\langle\eta_1,i\Lambda Q\rangle=0.}$ Let $A=0$ and $B=\Lambda Q$. Since $L_+(\Lambda Q)=-2Q$, $\langle\eta_1,Q\rangle=0$ and $\langle\vec{m}_1\cdot\vec{\mathrm{M}}Q,-\Lambda Q\rangle=c_1(\dot{\gamma}-1+|v|^2-\frac{\dot{\lambda}}{\lambda}(v\cdot z)-(v\cdot\dot{z}))$, we obtain

$$\left|\dot{\gamma}-1+|v|^2-\frac{\dot{\lambda}}{\lambda}(v\cdot z)-(v\cdot\dot{z})\right|\lesssim(C^*)^2s^{-2}+s^{-1}|\vec{m}_1^a|. \tag{3.32}$$

$\underline{\langle\eta_1,yQ\rangle=0.}$ Let $A=yQ$ and $B=0$. Since $L_-(yQ)=-2\nabla Q$, $|\langle\eta_1,i\nabla Q\rangle|\lesssim C^{**}s^{-2}$ and $\langle\vec{m}_1\cdot\vec{\mathrm{M}}Q,iyQ\rangle=c_2(\dot{z}-2v+\frac{\dot{\lambda}}{\lambda}z)$, we obtain

$$\left|\dot{z}-2v+\frac{\dot{\lambda}}{\lambda}z\right|\lesssim C^{**}s^{-2}+(C^*)^2s^{-2}+s^{-1}|\vec{m}_1|. \tag{3.33}$$

By (3.22) and (3.31),

$$\left|\dot{v}-\frac{\dot{\lambda}}{\lambda}v\right|\lesssim|\dot{v}|+\left|\frac{\dot{\lambda}}{\lambda}\right||v|\lesssim s^{-2}. \tag{3.34}$$

Combining (3.31)–(3.34), we have proved, for all $s\in[s^{**},s^{in}]$,

$$|\vec{m}_1(s)|\lesssim C^{**}s^{-2}+(C^*)^2s^{-2}. \tag{3.35}$$

Now we turn to the control of momentum $M(\mathbf{P})$:

$$\frac{1}{2}\frac{d}{ds}M(\mathbf{P}) = \text{Im}\int \dot{\mathbf{P}}\nabla\overline{\mathbf{P}}dy = \langle i\dot{\mathbf{P}}, \nabla\mathbf{P}\rangle.$$

From the equation of $\mathbf{P}$ (2.7), we get

$$\frac{1}{2}\frac{d}{ds}M(\mathbf{P}) = \langle \mathcal{E}_{\mathbf{P}} - \Delta\mathbf{P} + \mathbf{P} - |\mathbf{P}|^{p-1}\mathbf{P} + i\frac{\dot{\lambda}}{\lambda}\Lambda\mathbf{P} - (1 - \dot{\gamma})\mathbf{P}, \nabla\mathbf{P}\rangle$$

$$= \langle \mathcal{E}_{\mathbf{P}}, \nabla\mathbf{P}\rangle = \sum_{k=1}^{2}\langle \mathcal{E}_{\mathbf{P}}, \nabla P_k\rangle.$$

First, we note that

$$\nabla P_1 = [e^{i\Gamma_1}(\nabla Q + iv_1 Q)](y - z_1(s))$$

and $\mathcal{E}_{\mathbf{P}} = [e^{i\Gamma_1}\vec{m}_1 \cdot \vec{\mathrm{M}}Q](y - z_1(s)) + [e^{i\Gamma_2}\vec{m}_2 \cdot \vec{\mathrm{M}}Q](y - z_2(s)) + G$. We will control each term of $\langle \mathcal{E}_{\mathbf{P}}, \nabla P_1\rangle$ as follows

$$\langle [e^{i\Gamma_1}\vec{m}_1 \cdot \vec{\mathrm{M}}Q](y - z_1(s)) + G, [e^{i\Gamma_1}\nabla Q](y - z_1(s))\rangle$$

$$= c_2\left(\dot{v}_1 - \frac{\dot{\lambda}}{\lambda}v_1\right) + \langle G, e^{i\Gamma_1}\nabla Q(y - z_1(s))\rangle.$$

From (2.20) and the choice of $v$ in (3.6), we get

$$\left|c_2\dot{v}_1 + \langle G, e^{i\Gamma_1}\nabla Q(y - z_1(s))\rangle - c_2\frac{\dot{\lambda}}{\lambda}v_1\right| \lesssim (|v|^2|z|^2 + |v|^2)|z|^{-\frac{1}{2}}e^{-|z|} + |z|^{-\frac{3(d-1)}{4}}e^{-\frac{3}{2}|z|} + |v|\left|\frac{\dot{\lambda}}{\lambda}\right|$$

$$\lesssim (C^*)^2 s^{-3},$$

then from (3.31), we obtain

$$|\langle [e^{i\Gamma_1}\vec{m}_1 \cdot \vec{\mathrm{M}}Q](y - z_1(s)) + G, i[e^{i\Gamma_1}v_1 Q](y - z_1(s))\rangle| \lesssim |v|\left(\left|\frac{\dot{\lambda}}{\lambda}\right| + |G|\right)$$

$$\lesssim (C^*)^2 s^{-3}.$$

Finally, we have

$$\langle \vec{m}_2 \cdot \vec{\mathrm{M}}Q, e^{i\Gamma_1(y - (z_1 - z_2)) - i\Gamma_2(y)}(\nabla Q + ivQ)(y - (z_1 - z_2))\rangle \lesssim |\vec{m}_1||z|^{-\frac{3}{8}(d-1)}e^{-\frac{3}{4}|z|}$$

$$\lesssim s^{-3}\log^{-2}(s).$$

By symmetry, we obtain the same estimate for $\langle \mathcal{E}_{\mathbf{P}}, \nabla P_2\rangle$ and get

$$\left|\frac{d}{ds}M(\mathbf{P})\right| \lesssim (C^*)^2 s^{-3}$$

thus $|M(\mathbf{P}) - M(\mathbf{P}^{in})| \lesssim (C^*)^2 s^{-2}$. Note that the constant on the right-hand side does not depend on $C^{**}$. We consider

$$\langle \epsilon(s), i\nabla\mathbf{P}\rangle = \frac{1}{2}(M(u) - M(\mathbf{P}) - \text{Im}\int \nabla\epsilon\,\bar{\epsilon})$$

then from the conservation of momentum $M(u) = M(u^{in}) = M(\mathbf{P}^{in})$, we obtain

$$|\langle \epsilon(s), i\nabla\mathbf{P}\rangle| \lesssim (C^*)^2 s^{-2}.$$

Note that by the symmetry, $\langle \epsilon(s), i\nabla P_1 \rangle = \langle \epsilon(s), i\nabla P_2 \rangle = \frac{1}{2}\langle \epsilon(s), i\nabla \mathbf{P} \rangle$ and

$$\langle \epsilon(s), i\nabla P_1 \rangle = \langle \eta_1, i\nabla Q \rangle + \langle \eta_1, -vQ \rangle = \langle \eta_1, i\nabla Q \rangle + O(s^{-3}),$$

this information implies that $|\langle \eta_1, i\nabla Q \rangle| \lesssim (C^*)^2 s^{-2}$ so if we take $C^{**}$ big enough such that $\frac{C^{**}}{2} \gtrsim (C^*)^2$ then $s^{**} = s^*$. Those estimates (3.26) and (3.27) are direct consequences of (3.22), (3.24) and (3.31). $\qquad\square$

### 3.2.3. *Energy functional.* Consider the nonlinear energy functional for $\epsilon$

$$\mathbf{H}(s,\epsilon) = \frac{1}{2} \int \left( |\nabla \epsilon|^2 + |\epsilon|^2 - \frac{2}{p+1} \left( |\mathbf{P} + \epsilon|^{p+1} - |\mathbf{P}|^{p+1} - (p+1)|\mathbf{P}|^{p-1} \operatorname{Re}(\epsilon \overline{\mathbf{P}}) \right) \right).$$

Pick a smooth function $\chi : [0, +\infty) \to [0, \infty)$, non increasing, with $\chi \equiv 1$ on $[0, \frac{1}{10}]$, $\chi \equiv 0$ on $[\frac{1}{8}, +\infty)$. We define the localized momentum:

$$\mathbf{J} = \sum_k J_k, \quad J_k(s,\epsilon) = v_k \cdot \operatorname{Im} \int (\nabla \epsilon \, \bar\epsilon) \chi_k, \quad \chi_k(s,y) = \chi \left( \log^{-1}(s)|y - z_k(s)| \right).$$

Finally, set

$$\mathbf{W}(s,\epsilon) = \mathbf{H}(s,\epsilon) - \mathbf{J}(s,\epsilon).$$

The functional $\mathbf{W}$ is coercive in $\epsilon$ at the main order and it is an almost conserved quantity for the problem (see [24] for a similar functional).

**Proposition 6** (Coercivity and time control of the energy functional). *For all $s \in [s^*, s^{in}]$,*

$$\mathbf{W}(s,\epsilon(s)) \gtrsim \|\epsilon(s)\|_{H^1}^2, \tag{3.36}$$

*and*

$$\left| \frac{d}{ds}[\mathbf{W}(s,\epsilon(s))] \right| \lesssim s^{-2}\|\epsilon(s)\|_{H^1}. \tag{3.37}$$

*Proof of Proposition 6.* **step 1** Coercivity. The proof of the coercivity (3.36) is a standard consequence of the coercivity property (1.15) around one solitary wave with the orthogonality properties (3.5), (3.25), and an elementary localization argument. We refer to the proof of Lemma 4.1 in Appendix B of [18] for a similar proof.

**step 2** Variation of the energy. We estimate the time variation of the functional $\mathbf{H}$ and claim that for all $s \in [s^*, s^{in}]$,

$$\left| \frac{d}{ds}[\mathbf{H}(s,\epsilon(s))] - \sum_{k=1}^2 \dot{z}_k \cdot \left\langle \nabla P_k, \frac{\bar\epsilon.F''(\mathbf{P}).\epsilon}{2} \right\rangle \right| \lesssim s^{-2}\|\epsilon(s)\|_{H^1} + s^{-2}\|\epsilon\|_{H^1}^2. \tag{3.38}$$

The time derivative of $s \mapsto H(s,\epsilon(s))$ splits into two parts

$$\frac{d}{ds}[\mathbf{H}(s,\epsilon(s))] = D_s\mathbf{H}(s,\epsilon(s)) + \langle D_\epsilon\mathbf{H}(s,\epsilon(s)), \dot{\epsilon}_s \rangle,$$

where $D_s$ denotes differentiation of $\mathbf{H}$ with respect to $s$ and $D_\epsilon$ denotes differentiation of $\mathbf{H}$ with respect to $\epsilon$. Firstly we compute:

$$
\begin{aligned}
D_s\mathbf{H} = & -\operatorname{Re}\int[\dot{\mathbf{P}}(|\mathbf{P}+\epsilon|^{p-1}\overline{(\mathbf{P}+\epsilon)}) - |\mathbf{P}|^{p-1}\overline{\mathbf{P}}) \\
& -\frac{p-1}{2}|\mathbf{P}|^{p-3}(\dot{\mathbf{P}}\overline{\mathbf{P}} + \dot{\overline{\mathbf{P}}}\mathbf{P})\operatorname{Re}(\overline{\epsilon}\mathbf{P}) - |\mathbf{P}|^{p-1}\overline{\epsilon}\dot{\mathbf{P}}](y)dy \\
= & -\operatorname{Re}\int[\dot{\mathbf{P}}(|\mathbf{P}+\epsilon|^{p-1}\overline{(\mathbf{P}+\epsilon)}) - |\mathbf{P}|^{p-1}\overline{\mathbf{P}}) \\
& -\frac{p-1}{2}|\mathbf{P}|^{p-3}\frac{\epsilon\overline{\mathbf{P}}^2\dot{\mathbf{P}} + \overline{\epsilon}|\mathbf{P}|^2\dot{\mathbf{P}} + \epsilon|\mathbf{P}|^2\dot{\overline{\mathbf{P}}} + \overline{\epsilon}\mathbf{P}^2\dot{\overline{\mathbf{P}}}}{2} - |\mathbf{P}|^{p-1}\overline{\epsilon}\dot{\mathbf{P}}](y)dy \\
= & -\langle\dot{\mathbf{P}}, |\mathbf{P}+\epsilon|^{p-1}(\mathbf{P}+\epsilon) - |\mathbf{P}|^{p-1}\mathbf{P} - \frac{p+1}{2}\epsilon|\mathbf{P}|^{p-1} - \frac{p-1}{2}\overline{\epsilon}\mathbf{P}^2|\mathbf{P}|^{p-3}\rangle.
\end{aligned}
$$

We observe that $\dot{P}_k = -\dot{z}_k\cdot\nabla P_k + i\dot{v}_k\cdot(y-z_k)P_k$. Denote

$$
K = |\mathbf{P}+\epsilon|^{p-1}(\mathbf{P}+\epsilon) - |\mathbf{P}|^{p-1}\mathbf{P} - \frac{p+1}{2}\epsilon|\mathbf{P}|^{p-1} - \frac{p-1}{2}\overline{\epsilon}\mathbf{P}^2|\mathbf{P}|^{p-3}
$$

then by (2.16), $K = |\mathbf{P}+\epsilon|^{p-1}(\mathbf{P}+\epsilon) - |\mathbf{P}|^{p-1}\mathbf{P} - F'(\mathbf{P}).\epsilon$, we deduce from (2.17) that

$$
|K| \lesssim |\epsilon|^2|\mathbf{P}|^{p-2} + |\epsilon|^p
$$

so we obtain

$$
|\langle i\dot{v}_k\cdot(y-z_k)P_k, K\rangle| \lesssim (\|\epsilon\|_{H^1}^2 + \|\epsilon\|_{H^1}^p)|\dot{v}| \lesssim s^{-2}\|\epsilon\|_{H^1}^2.
$$

Next we look more precisely at $K$

$$
K = \frac{\overline{\epsilon}.F''(\mathbf{P}).\epsilon}{2} + O\Big(\Big|\frac{\epsilon}{\mathbf{P}}\Big|^{2^+}|\mathbf{P}|^p\Big) + O(|\epsilon|^p)
$$

as $|\dot{z}_k| \lesssim s^{-1}$ and $p-2^+ > 0$, we have

$$
\Big|\Big\langle -\dot{z}_k\cdot\nabla P_k, \Big|\frac{\epsilon}{\mathbf{P}}\Big|^{2^+}|\mathbf{P}|^p\Big\rangle\Big| \lesssim s^{-1}\|\epsilon\|_{H^1}^{2^+}
$$

and

$$
|\langle -\dot{z}_k\cdot\nabla P_k, |\epsilon|^p\rangle| \lesssim s^{-1}\|\epsilon\|_{H^1}^p.
$$

Combining these computations, we get

$$
D_s\mathbf{H}(s,\epsilon) = \sum_{k=1}^2\langle\dot{z}_k\cdot\nabla P_k, \frac{\overline{\epsilon}.F''(\mathbf{P}).\epsilon}{2}\rangle + O(s^{-1}\|\epsilon\|_{H^1}^{2^+}) + O(s^{-2}\|\epsilon\|_{H^1}^2) + O(s^{-1}\|\epsilon\|_{H^1}^p). \quad (3.39)
$$

Secondly we consider

$$
D_\epsilon\mathbf{H}(s,\epsilon) = -\Delta\epsilon + \epsilon - \big(|\mathbf{P}+\epsilon|^{p-1}(\mathbf{P}+\epsilon) - |\mathbf{P}|^{p-1}\mathbf{P}\big)
$$

and note that the equation (3.8) of $\epsilon$ can be rewritten as

$$
i\dot{\epsilon} - D_\epsilon\mathbf{H}(s,\epsilon) - i\frac{\dot{\lambda}}{\lambda}\Lambda\epsilon + (1-\dot{\gamma})\epsilon + \mathcal{E}_\mathbf{P} = 0
$$

so that

$$\langle D_\epsilon \mathbf{H}(s,\epsilon), \dot\epsilon \rangle = \langle i D_\epsilon \mathbf{H}(s,\epsilon), i\dot\epsilon \rangle$$
$$= \frac{\dot\lambda}{\lambda} \langle D_\epsilon \mathbf{H}(s,\epsilon), \Lambda\epsilon \rangle - (1-\dot\gamma)\langle D_\epsilon \mathbf{H}(s,\epsilon), \epsilon \rangle - \langle i D_\epsilon \mathbf{H}(s,\epsilon), \mathcal{E}_{\mathbf{P}} \rangle.$$

On the other hand, from (3.24) and (3.21)–(3.22), we have

$$\left| \frac{\dot\lambda}{\lambda} \langle D_\epsilon \mathbf{H}(s,\epsilon), \Lambda\epsilon \rangle \right| \lesssim \left| \frac{\dot\lambda}{\lambda} \right| \left( \|\epsilon\|_{H^1}^2 + \|\epsilon\|_{H^1}^{p+1} \right) \lesssim (C^*)^2 s^{-2} \|\epsilon\|_{H^1}^2,$$

$$|(1-\dot\gamma)\langle i D_\epsilon \mathbf{H}(s,\epsilon), \epsilon \rangle| \lesssim |1-\dot\gamma| \left( \|\epsilon\|_{H^1}^2 + \|\epsilon\|_{H^1}^{p+1} \right) \lesssim (C^*)^2 s^{-2} \|\epsilon\|_{H^1}^2.$$

For the last term, we rewrite

$$\langle i D_\epsilon \mathbf{H}(s,\epsilon), \mathcal{E}_{\mathbf{P}} \rangle = \langle -i\Delta\epsilon + i\epsilon - i \left( |\mathbf{P}+\epsilon|^{p-1}(\mathbf{P}+\epsilon) - |\mathbf{P}|^{p-1}\mathbf{P} \right),$$
$$[e^{i\Gamma_1}\vec{m}_1 \cdot \vec{\mathrm{M}}Q](y-z_1(s)) + [e^{i\Gamma_2}\vec{m}_2 \cdot \vec{\mathrm{M}}Q](y-z_2(s)) + G \rangle.$$

Recall that with $\eta_1 = \eta_1^1 + i\eta_1^2$ for $\eta_1^1, \eta_1^2$ real, from the expression of operators $L_+$ and $L_-$

$$I_1 = \langle -i\Delta\epsilon + i\epsilon - i \left( |\mathbf{P}+\epsilon|^{p-1}(\mathbf{P}+\epsilon) - |\mathbf{P}|^{p-1}\mathbf{P} \right), [e^{i\Gamma_1}\vec{m}_1 \cdot \vec{\mathrm{M}}Q](y-z_1(s)) \rangle$$
$$= \langle -i\Delta\eta_1 + i\eta_1 - i(|\mathbf{P}_1+\eta_1|^{p-1}(\mathbf{P}_1+\eta_1) - |\mathbf{P}_1|^{p-1}\mathbf{P}_1), \vec{m}_1 \cdot \vec{\mathrm{M}}Q \rangle$$
$$= \langle iL_+\eta_1^1 - L_-\eta_1^2, \vec{m}_1 \cdot \vec{\mathrm{M}}Q \rangle$$
$$\quad - \left\langle i \left( |\mathbf{P}_1+\eta_1|^{p-1}(\mathbf{P}_1+\eta_1) - |\mathbf{P}_1|^{p-1}\mathbf{P}_1 - \frac{p+1}{2}Q^{p-1}\eta_1 - \frac{p-1}{2}Q^{p-1}\overline{\eta}_1 \right), \vec{m}_1 \cdot \vec{\mathrm{M}}Q \right\rangle$$
$$= -\frac{\dot\lambda}{\lambda} \langle \eta_1, -2Q \rangle + (\dot{v} - \frac{\dot\lambda}{\lambda}v)\langle \eta_1, -2i\nabla Q \rangle$$
$$\quad - \left\langle i \left( |\mathbf{P}_1+\eta_1|^{p-1}(\mathbf{P}_1+\eta_1) - |\mathbf{P}_1|^{p-1}\mathbf{P}_1 - \frac{p+1}{2}Q^{p-1}\eta_1 - \frac{p-1}{2}Q^{p-1}\overline{\eta}_1 \right), \vec{m}_1 \cdot \vec{\mathrm{M}}Q \right\rangle.$$

By orthogonality of $\eta_1$ (3.5), (3.25) and the estimate (3.24), (3.29), we get

$$|I_1| = O((C^*)^2 s^{-4}) + O((C^*)^4 s^{-4}).$$

By symmetry, we have the same estimate for $I_2$. Finally, from (2.13) and (3.21), we have $\|G\|_{H^1} \lesssim s^{-2}$ so using integration by parts and Cauchy-Schwarz inequality

$$|\langle -i\Delta\epsilon + i\epsilon - i \left( |\mathbf{P}+\epsilon|^{p-1}(\mathbf{P}+\epsilon) - |\mathbf{P}|^{p-1}\mathbf{P} \right), G \rangle| \lesssim s^{-2} \|\epsilon\|_{H^1}. \tag{3.40}$$

The collection of above estimates finishes the proof of (3.38).

**step 3** Variation of the localized momentum. We now claim: for all $s \in [s^*, s^{in}]$,

$$\left| \frac{d}{ds}[\mathbf{J}(s,\epsilon(s))] - \sum_{k=1}^2 2v_k \cdot \langle \nabla P_k, \frac{\overline{\epsilon}.F''(\mathbf{P}).\epsilon}{2} \rangle \right| \lesssim s^{-1}\log^{-1}(s)\|\epsilon(s)\|_{H^1}^2 + s^{-\frac{5}{2}}\|\epsilon(s)\|_{H^1}. \tag{3.41}$$

Indeed, we compute, for any $k$,

$$\frac{d}{ds}[J_k(s,\epsilon(s))] = \dot{v}_k \cdot \text{Im} \int (\nabla\epsilon\,\overline{\epsilon})\chi_k + v_k \cdot \text{Im} \int (\nabla\epsilon\,\overline{\epsilon})\dot{\chi}_k + v_k \langle i\dot\epsilon, 2\chi_k\nabla\epsilon + \epsilon\nabla\chi_k \rangle.$$

By (3.21) and (3.22), we have

$$\left| \dot{v}_k \cdot \text{Im} \int (\nabla\epsilon\,\overline{\epsilon})\chi_k \right| \lesssim s^{-2}\|\epsilon\|_{H^1}^2.$$

Note that by direct computations, (3.24) and (3.21)–(3.22),

$$|\dot{\chi}_k| \lesssim (s^{-1}\log^{-1}(s)|y - z_k| + |\dot{z}_k|)\log^{-1}(s)|\chi'(\log^{-1}(s)(y - z_k(s)))| \lesssim s^{-1}\log^{-1}(s)$$

and so, by (3.21)–(3.22),

$$\left| v_k \cdot \operatorname{Im} \int (\nabla\epsilon\,\bar{\epsilon})\dot{\chi}_k \right| \lesssim s^{-2}\log^{-1}(s)\|\epsilon\|_{H^1}^2.$$

Now, we use the equation (3.8) of $\epsilon$ to estimate $v_k\langle i\dot{\epsilon}, 2\chi_k\nabla\epsilon + \epsilon\nabla\chi_k\rangle$. By integration by parts, we check the following

$$\langle \Delta\epsilon, 2\chi_k\nabla\epsilon + \epsilon\nabla\chi_k\rangle = -\langle\nabla\epsilon\cdot\nabla\chi_k, \nabla\epsilon\rangle + \frac{1}{2}\int|\epsilon|^2\nabla(\Delta\chi_k).$$

We have

$$|v_k\langle\nabla\epsilon\cdot\nabla\chi_k, \nabla\epsilon\rangle| \lesssim s^{-1}\log^{-1}(s)\|\epsilon\|_{H^1}^2$$

and as $|\nabla(\Delta\chi_k)| \lesssim \log^{-3}(s)$ we obtain

$$\left| v_k \cdot \int|\epsilon|^2\nabla(\Delta\chi_k) \right| \lesssim s^{-1}\log^{-3}(s)\|\epsilon\|_{H^1}^2.$$

In conclusion for term $\Delta\epsilon$ in the equation of $\epsilon$, we get

$$|v_k\langle\Delta\epsilon, 2\chi_k\nabla\epsilon + \epsilon\nabla\chi_k\rangle| \lesssim s^{-1}\log^{-1}(s)\|\epsilon\|_{H^1}^2.$$

For the term $\epsilon$, we simply verify by integration by parts that

$$\langle\epsilon, 2\chi_k\nabla\epsilon + \epsilon\nabla\chi_k\rangle = 0.$$

We also have that

$$\left| v_k\langle\frac{\dot{\lambda}}{\lambda}i\Lambda\epsilon, 2\chi_k\nabla\epsilon + \epsilon\nabla\chi_k\rangle \right| \lesssim |v_k|\left|\frac{\dot{\lambda}}{\lambda}\right|\|\epsilon\|_{H^1}^2 \lesssim (C^*)^2 s^{-3}\|\epsilon\|_{H^1}^2,$$

$$\left| v_k\langle\mathcal{E}_{\mathbf{P}}, 2\chi_k\nabla\epsilon + \epsilon\nabla\chi_k\rangle \right| \lesssim |v_k|\,|\mathcal{E}_{\mathbf{P}}|\log^{\frac{d}{2}}(s)\|\epsilon\|_{H^1} \lesssim s^{-\frac{5}{2}}\|\epsilon\|_{H^1}.$$

where the last estimate, we use Cauchy-Schwarz inequality and the fact that the support of $\chi_k$ is contained in $\{y|\,|y - z_k(s)| \leq \frac{1}{8}\log(s)\}$. Now we only have to deal with the term $v_k\langle|\mathbf{P} + \epsilon|^{p-1}(\mathbf{P} + \epsilon) - |\mathbf{P}|^{p-1}\mathbf{P}, 2\chi_k\nabla\epsilon + \epsilon\nabla\chi_k\rangle$. On the one hand, by (2.17), we consider

$$|\mathbf{P} + \epsilon|^{p-1}(\mathbf{P} + \epsilon) - |\mathbf{P}|^{p-1}\mathbf{P} = F'(\mathbf{P}).\epsilon + O(|\epsilon|^p) + O\left(\left|\frac{\epsilon}{\mathbf{P}}\right|^2|\mathbf{P}|^p\right)$$

and using Gagliardo-Nirenberg inequality (note that if $p > 2$ then $3 < 2^*$)

$$|v_k|\left|\left\langle\left|\frac{\epsilon}{\mathbf{P}}\right|^2|\mathbf{P}|^p, 2\chi_k\nabla\epsilon + \epsilon\nabla\chi_k\right\rangle\right| \lesssim s^{-1}\|\epsilon\|_{H^1}^3,$$

$$|v_k|\,|\langle|\epsilon|^p, 2\chi_k\nabla\epsilon + \epsilon\nabla\chi_k\rangle| \lesssim s^{-1}\|\epsilon\|_{H^1}^{p+1} \lesssim s^{-1}\|\epsilon\|_{H^1}^3.$$

On the other hand

$$|v_k\langle F'(\mathbf{P}).\epsilon, \epsilon\nabla\chi_k\rangle| \lesssim |v_k|\,|\nabla\chi_k|\,\|\epsilon\|_{H^1}^2 \lesssim s^{-1}\log^{-1}(s)\|\epsilon\|_{H^1}^2.$$

Finally by integration by parts, we get

$$v_k\langle F'(\mathbf{P}).\epsilon, \chi_k\nabla\epsilon\rangle = -\frac{1}{2}v_k\langle\nabla\mathbf{P}\chi_k, \bar{\epsilon}.F''(\mathbf{P}).\epsilon\rangle - \frac{1}{2}v_k\langle F'(\mathbf{P}).\epsilon, \epsilon\nabla\chi_k\rangle,$$

therefore the collection of above bounds gives

$$\frac{d}{ds}[\mathbf{J}(s,\epsilon(s))] = \sum_{k=1}^{2}\langle 2v_k \cdot \nabla\mathbf{P}\chi_k, \frac{\bar{\epsilon}.F''(\mathbf{P}).\epsilon}{2}\rangle + O(s^{-1}\log^{-1}(s)\|\epsilon(s)\|_{H^1}^2)$$
$$+ O(s^{-\frac{5}{2}}\|\epsilon(s)\|_{H^1}) + O(s^{-1}\|\epsilon(s)\|_{H^1}^3). \quad (3.42)$$

We finish the proof of (3.41) by showing the following estimate

$$|v_k \cdot \langle \nabla\mathbf{P}\chi_k, \bar{\epsilon}.F''(\mathbf{P}).\epsilon\rangle - v_k \cdot \langle \nabla P_k, \bar{\epsilon}.F''(\mathbf{P}).\epsilon\rangle| \lesssim |v_k|\left[\sum_{j\neq k}\left|\int_{|y-z_k(s)|<\frac{1}{8}\log s}(\bar{\epsilon}.F''(\mathbf{P}).\epsilon)\nabla P_j\right|\right.$$
$$\left.+ \left|\int_{|y-z_k(s)|>\frac{1}{10}\log s}(\bar{\epsilon}.F''(\mathbf{P}).\epsilon)\nabla P_k\right|\right] \lesssim s^{-\frac{1}{10}-1}\|\epsilon\|_{H^1}^2, \quad (3.43)$$

here we use (1.12).

**step 4** Conclusion. Recall that, by (3.26), $|\dot{z}_k - 2v_k| \lesssim s^{-1}\log^{-1}(s)$ so

$$\left|(\dot{z}_k - 2v_k)\cdot\langle\nabla P_k, \frac{\bar{\epsilon}.F''(\mathbf{P}).\epsilon}{2}\rangle\right| \lesssim s^{-1}\log^{-1}(s)\|\epsilon\|_{H^1}^2,$$

and (3.37) now follows from (3.38), (3.41). This concludes the proof of Proposition 6. $\square$

3.2.4. **End of the bootstrap argument.** We close the bootstrap estimates (3.21).

**step 1** Closing the estimate in $\epsilon$. By (3.37) in Proposition 6 and then (3.21)–(3.22), we have

$$\left|\frac{d}{ds}[\mathbf{W}(s,\epsilon(s))]\right| \lesssim s^{-2}\|\epsilon\|_{H^1} \lesssim C^*s^{-3}.$$

Thus, by integration on $[s, s^{in}]$ for any $s \in [s^*, s^{in}]$, using $\epsilon(s^{in}) = 0$ (see (3.18)), we obtain

$$|\mathbf{W}(s,\epsilon(s))| \lesssim C^*s^{-2}.$$

By (3.36) in Proposition 6, we get

$$\|\epsilon(s)\|_{H^1}^2 \leq C_0 C^* s^{-2}.$$

Therefore, for $C^*$ large enough such that $C_0 C^* \leq \frac{(C^*)^2}{4}$, we have $\|\epsilon\|_{H^1} \leq \frac{C^*}{2}s^{-1}$ which strictly improves the estimate on $\|\epsilon\|_{H^1}$ in (3.21).

**step 2** Closing the parameter $z$. Now, we need to finish the bootstrap argument for $z(s)$. Note that

$$\left|\dot{v} + c\frac{z}{|z|}|z|^{-\frac{d-1}{2}}e^{-|z|}\right| \lesssim s^{-2}\log^{-1}(s)$$
$$|\dot{z} - 2v| \lesssim s^{-1}\log^{-1}(s)$$

thus we deduce

$$\left|\dot{v}\cdot\frac{z}{|z|} + c|z|^{-\frac{d-1}{2}}e^{-|z|}\right| \lesssim s^{-2}\log^{-1}(s)$$
$$\left|\dot{z}\cdot\frac{z}{|z|} - 2v\cdot\frac{z}{|z|}\right| \lesssim s^{-1}\log^{-1}(s).$$

We get

$$\left|2\left(v\cdot\frac{z}{|z|}\right)\left(\dot{v}\cdot\frac{z}{|z|}\right) + c\,\dot{z}\cdot\frac{z}{|z|}|z|^{-\frac{d-1}{2}}e^{-|z|}\right| \lesssim s^{-3}\log^{-1}(s)$$

since $|v| \lesssim s^{-1}, |\dot{v}| \lesssim s^{-2}$. Therefore, by the explicit choice of initial data

$$v(s^{in}) = \sqrt{c}(z^{in})^{-\frac{d-1}{4}} e^{-\frac{1}{2}z^{in}} \vec{e}_1, \quad z(s^{in}) = z^{in} \vec{e}_1,$$

we integrate on $[s, s^{in}]$ for any $s \in [s^*, s^{in}]$, if $d - 1 > 0$

$$\left| \left( v \cdot \frac{z}{|z|} \right)^2 - c|z|^{-\frac{d-1}{2}} e^{-|z|} \right| \lesssim s^{-2} \log^{-1}(s) + \int_s^{s^{in}} |\dot{z}||z|^{-\frac{d-1}{2}-1} e^{-|z|} \lesssim s^{-2} \log^{-1}(s),$$

if $d - 1 = 0$, $\left| 2(v \cdot \frac{z}{|z|})(\dot{v} \cdot \frac{z}{|z|}) + c \dot{z} \cdot \frac{z}{|z|} e^{-|z|} \right| \lesssim s^{-3} \log^{-1}(s)$ implies also $\left| (v \cdot \frac{z}{|z|})^2 - ce^{-|z|} \right| \lesssim s^{-2} \log^{-1}(s)$. In both cases, combining with (3.26), we get

$$\left| (v \cdot \frac{z}{|z|}) - \sqrt{c}|z|^{-\frac{d-1}{4}} e^{-\frac{1}{2}|z|} \right| + \left| (\dot{z} \cdot \frac{z}{|z|}) - 2(v \cdot \frac{z}{|z|}) \right| \lesssim s^{-1} \log^{-1}(s)$$

so $\left| (\dot{z} \cdot \frac{z}{|z|}) - 2\sqrt{c}|z|^{-\frac{d-1}{4}} e^{-\frac{1}{2}|z|} \right| \lesssim s^{-1} \log^{-1}(s)$. Next, note that if $d - 1 > 0$

$$\frac{d}{ds}(|z|^{\frac{d-1}{4}} e^{\frac{1}{2}|z|}) = \frac{1}{2}\dot{z} \cdot \frac{z}{|z|}|z|^{\frac{d-1}{4}} e^{\frac{1}{2}|z|} + \frac{d-1}{4}\dot{z} \cdot \frac{z}{|z|}|z|^{\frac{d-1}{4}-1} e^{\frac{1}{2}|z|}$$

and if $d - 1 = 0$

$$\frac{d}{ds}(e^{\frac{1}{2}|z|}) = \frac{1}{2}\dot{z} \cdot \frac{z}{|z|} e^{\frac{1}{2}|z|}$$

thus

$$\left| \frac{d}{ds}\left( |z|^{\frac{d-1}{4}} e^{\frac{1}{2}|z|} \right) - c^{\frac{1}{2}} \right| \lesssim \log^{-1}(s) + \frac{d-1}{4}|\dot{z}||z|^{\frac{d-1}{4}-1} e^{\frac{1}{2}|z|} \lesssim \log^{-1}(s) \quad (3.44)$$

here we use $|z| \lesssim \log^{-1}(s)$ and $|\dot{z}| \lesssim s^{-1}$. Next, we need to adjust the initial choice of $z^{in}$ through a topological argument (see [4] for a similar argument). We define $\zeta$ and $\xi$ the following two functions on $[s^*, s^{in}]$

$$\zeta(s) = c^{-\frac{1}{2}}|z|^{\frac{d-1}{4}} e^{\frac{1}{2}|z|}, \quad \xi(s) = (\zeta(s) - s)^2 s^{-2} \log(s). \quad (3.45)$$

Then, (3.44) writes

$$|\dot{\zeta}(s) - 1| \lesssim \log^{-1}(s). \quad (3.46)$$

According to (3.21), our objective is to prove that there exists a suitable choice of

$$\zeta(s^{in}) = \zeta^{in} \in [s^{in} - s^{in} \log^{-\frac{1}{2}}(s^{in}), s^{in} + s^{in} \log^{-\frac{1}{2}}(s^{in})],$$

so that $s^* = s_0$. Assume for the sake of contradiction that for all $\zeta^{\sharp} \in [-1, 1]$, the choice

$$\zeta^{in} = s^{in} + \zeta^{\sharp} s^{in} \log^{-\frac{1}{2}}(s^{in})$$

leads to $s^* = s^*(\zeta^{\sharp}) \in (s_0, s^{in})$. Since all estimates in (3.21) except the one on $z(s)$ have been strictly improved on $[s^*, s^{in}]$, it follows from $s^*(\zeta^{\sharp}) \in (s_0, s^{in})$ and continuity that

$$|\zeta(s^*(\zeta^{\sharp})) - s^*| = s^* \log^{-\frac{1}{2}} s^* \quad \text{i.e.} \quad \zeta(s^*(\zeta^{\sharp})) = s^* \pm s^* \log^{-\frac{1}{2}} s^*.$$

We need a transversality condition to reach a contradiction. We compute:

$$\dot{\xi}(s) = 2(\zeta(s) - s)(\dot{\zeta}(s) - 1)s^{-2}\log(s) - (\zeta(s) - s)^2(2s^{-3}\log(s) - s^{-3}). \quad (3.47)$$

At $s = s^*$, this gives

$$|\dot{\xi}(s^*) + 2(s^*)^{-1}| \lesssim (s^*)^{-1} \log^{-\frac{1}{2}}(s^*).$$

Thus, for $s_0$ large enough,

$$\dot{\xi}(s^*) < -(s^*)^{-1}. \quad (3.48)$$

A consequence of the transversality property (3.48) is the continuity of the function $\zeta^\sharp \in [-1, 1] \mapsto s^*(\zeta^\sharp)$. Indeed, let $\epsilon > 0$ then there exists $\delta > 0$ such that $\xi(s^*(\zeta^\sharp) - \epsilon) > 1 + \delta$ and $\xi(s^*(\zeta^\sharp) + \epsilon) < 1 - \delta$. Moreover, by definition of $s^*(\zeta^\sharp)$ (choosing $\delta$ small enough) for all $s \in [s^*(\zeta^\sharp) + \epsilon, s^{in}]$ we have $\xi(s) < 1 - \delta$. But from the continuity of the flow, there exists $\iota > 0$ such that for all $|\tilde\zeta^\sharp - \zeta^\sharp| < \iota$

$$\forall s \in [s^*(\zeta^\sharp) - \epsilon, s^{in}], \quad |\tilde\xi(s) - \xi(s)| \leq \delta/2$$

so we obtain that $s^*(\zeta^\sharp) - \epsilon \leq s^*(\tilde\zeta^\sharp) \leq s^*(\xi^\sharp) + \epsilon$ and the continuity of $s^*(\zeta^\sharp)$ as expected. Thus we deduce the continuity of the function $\Phi$ defined by

$$\Phi : \zeta^\sharp \in [-1, 1] \mapsto (\zeta(s^*) - s^*)(s^*)^{-1} \log^{\frac{1}{2}}(s^*) \in \{-1, 1\}.$$

Moreover, for $\zeta^\sharp = -1$ and $\zeta^\sharp = 1$, in these two cases $\xi(s^{in}) = 1$, from (3.47) we have that $\dot\xi(s^{in}) < 0$ thus $s^* = s^{in}$. Therefore, $\Phi(-1) = -1$ and $\Phi(1) = 1$, but this is a contradiction with the continuity.

In conclusion, there exists at least a choice of

$$\zeta(s^{in}) = \zeta^{in} \in (s^{in} - s^{in} \log^{-\frac{1}{2}}(s^{in}), s^{in} + s^{in} \log^{-\frac{1}{2}}(s^{in}))$$

such that $s^* = s_0$. This concludes our bootstrap argument.

**step 3** Estimate on the parameter $\lambda$. From (3.24), we obtain

$$\left| \frac{\dot\lambda}{\lambda} \right| \lesssim s^{-2}.$$

By integration on $[s, s^{in}]$, for any $s \in [s_0, s^{in}]$, using the value $\lambda(s^{in}) = \lambda^{in} = 1$ (see (3.19)), we have

$$|\log(\lambda(s))| \lesssim s^{-1},$$

and thus

$$|\lambda(s) - 1| \lesssim s^{-1}$$

or in other words

$$\left| \lambda^{-1}(s) - 1 \right| \lesssim s^{-1}. \tag{3.49}$$

$$\square$$

## 4. Compactness arguments

### 4.1. Construction of a sequence of backwards solutions.

**Lemma 7.** *There exist $t_0 > 1$ and a sequence of solutions $u_n \in \mathcal{C}([t_0, T_n], H^1)$ of* (NLS), *where*

$$T_n \to +\infty \quad \text{as } n \to +\infty, \tag{4.1}$$

*satisfying the following estimates, for all $t \in [t_0, T_n]$,*

$$|\,|z_n(t)| - 2\log t| \lesssim \log(\log t), \quad |\lambda_n^{-1}(t) - 1| \lesssim t^{-1},$$

$$|v_n(t)| \lesssim t^{-1}, \quad \|\epsilon_n(t)\|_{H^1} \lesssim t^{-1}, \quad \left| |z_n(t)|^{\frac{d-1}{2}} e^{|z_n(t)|} - ct^2 \right| \lesssim t^2 \log^{-\frac{1}{2}}(t), \tag{4.2}$$

*where $(\lambda_n, z_n, \gamma_n, v_n)$ are the parameters of the decomposition of $u_n$, i.e.*

$$u_n(t, x) = \frac{e^{i\gamma_n(t)}}{\lambda_n^{\frac{2}{p-1}}(t)} \left( \sum_{k=1}^2 \left[ e^{i\Gamma_{k,n}} Q \right] \left( \frac{x}{\lambda_n(t)} + \frac{(-1)^k}{2} z_n(t) \right) + \epsilon_n \left( t, \frac{x}{\lambda_n(t)} \right) \right), \tag{4.3}$$

with $\Gamma_{k,n}(t,x) = \frac{(-1)^{k+1}}{2}v_n(t) \cdot \frac{x}{\lambda_n(t)}$.

*Proof of Lemma 7.* Applying Proposition 4 with $s^{in} = n$ for any large $n$, there exists a solution $u_n(t)$ of (NLS) defined on the time interval $[0, T_n]$ where

$$T_n = \int_{s_0}^{n} \lambda_n^2(s)ds.$$

and whose decomposition satisfies the uniform estimates (3.20). First, we see that $T_n \to +\infty$ as $n \to +\infty$ which follows directly from the estimate on $\lambda_n(s)$. From the definition of the rescaled time $s$ (see (3.2)), for any $s \in [s_0, n]$, we have

$$t(s) = \int_{s_0}^{s} \lambda_n^2(s')ds' \quad \text{where} \quad |\lambda_n^2(s) - 1| \lesssim s^{-1}.$$

Fix $t_0 = \bar{s}_0$ with $\bar{s}_0 > s_0$ large enough independent of $n$ such that for all $s$ with $n \geq s > \bar{s}_0$

$$\frac{1}{2}s \leq \int_{s_0}^{s} \lambda_n^2(s')ds' = s\left(1 + O(s^{-1})\right) \leq \frac{3}{2}s$$

then, for all $t \in [t_0, T_n]$

$$t(s) = s\left(1 + O(s^{-1})\right) \geq \frac{1}{2}s$$

and

$$s = t\left(1 + O(t^{-1})\right).$$

Thus, we get from (3.20)

$$\left|\,|z_n(s)| - 2\log(s)\right| \lesssim \log(\log(s)) \Leftrightarrow \left|\,|z_n(s(t))| - 2\log(t)\right| \lesssim \log(\log(t))$$

$$\left|\lambda_n^{-1}(s) - 1\right| \lesssim s^{-1} \Leftrightarrow \left|\lambda_n^{-1}(s(t)) - 1\right| \lesssim t^{-1}$$

$$\|\epsilon_n(s)\|_{H^1} \lesssim s^{-1} \Leftrightarrow \|\epsilon_n(s(t))\|_{H^1} \lesssim t^{-1}$$

$$|v_n(s)| \lesssim s^{-1} \Leftrightarrow |v_n(s(t))| \lesssim t^{-1}.$$

□

**4.2. Compactness argument.** Next, we claim a strong compactness result in $L^2(\mathbb{R}^d)$.

**Lemma 8.** *There exist $u_0 \in H^1(\mathbb{R}^d)$ and a sub-sequence, still denoted $u_n$, such that*

$$u_n(t_0) \rightharpoonup u_0 \text{ weakly in } H^1(\mathbb{R}^d)$$

$$u_n(t_0) \to u_0 \text{ in } H^\sigma(\mathbb{R}^d), \text{ for } 0 \leq \sigma < 1$$

*as $n \to \infty$.*

*Proof of Lemma 8.* By interpolation, it is enough to prove that the sub-sequence $u_n(t_0) \xrightarrow{L^2} u_0$ as $n \to \infty$. First, we claim the following: $\forall \delta_1 > 0, \delta_1 \ll 1, \exists n_0 \gg 1, \exists K_1 = K_1(\delta_1) > 0$ such that $\forall n \geq n_0$

$$\int_{|x| > K_1} |u_n(t_0, x)|^2 dx < \delta_1. \tag{4.4}$$

Indeed, denote $x_n(t) = z_n(t)\lambda_n(t)$ and

$$\tilde{R}_n(t,x) = e^{i\gamma_n(t)}\sum_{k=1}^{2}\left[e^{i\Gamma_{k,n}}Q_{\lambda_n^{-1}(t)}\right]\left(x + \frac{(-1)^k}{2}x_n(t)\right)$$

$$R_n(t, x) = e^{i\gamma_n(t)} \sum_{k=1}^{2} Q\left(x + \frac{(-1)^k}{2} x_n(t)\right)$$

then we have

$$||u_n(t) - R_n(t)||_{H^1} \leq ||\epsilon_n(t)||_{H^1} + 2||\tilde{R}_n(t) - R_n(t)||_{H^1}$$
$$\lesssim ||\epsilon_n(t)||_{H^1} + |\lambda_n^{-1}(t) - 1| + |v_n(t)| \lesssim t^{-1}. \quad (4.5)$$

We get a direct consequence of the above estimate

$$||u_n(t)||_{H^1} < C \quad (4.6)$$

for all $t \in [t_0, T_n]$ since $||R_n(t)||_{H^1} \leq 2||Q||_{H^1}$. Furthermore, for fixed $\delta_1$, there exists $t_1 > t_0$ such that

$$||u_n(t_1) - R_n(t_1)||_{H^1} \lesssim (t_1)^{-1} < \sqrt{\delta_1}$$

for $n$ large enough that $T_n > t_1$; in others words, we have

$$\int |u_n(t_1, x) - R_n(t_1, x)|^2 dx < \delta_1.$$

Besides, $|x_n(t_1) - 2\log(t_1)| \lesssim \log(\log t_1)$ then for $K_2 \gg 1$ large enough we have

$$\int_{|x| > K_2} |R_n(t_1, x)|^2 dx < \delta_1.$$

Consider now a $\mathcal{C}^1$ cut-off function $g : \mathbb{R} \to [0, 1]$ such that : $g \equiv 0$ on $(-\infty, 1]$, $0 < g' < 2$ on $(1, 2)$ and $g \equiv 1$ on $[2, +\infty)$. Since $||u_n(t)||_{H^1} < C$ bounded in $H^1$ independently of $n$ and $t \in [t_0, T_n]$, we can choose $\gamma_1 > 0$ independent of $n$ such that

$$\gamma_1 \geq \frac{2}{\delta_1}(t_1 - t_0)C^2.$$

We have by direct calculations, for $t \in [t_0, T_n]$

$$\left| \frac{d}{dt} \int |u_n(t, x)|^2 g\left(\frac{|x| - K_2}{\gamma_1}\right) dx \right| = \left| \frac{1}{\gamma_1} \operatorname{Im} \int u\left(\nabla \bar{u} \cdot \frac{x}{|x|}\right) g'\left(\frac{|x| - K_2}{\gamma_1}\right) \right|$$
$$\leq \frac{2}{\gamma_1} \sup_{T_n \geq t \geq t_0} ||u_n(t)||_{H^1}^2 \leq \frac{\delta_1}{t_1 - t_0}.$$

By integration from $t_0$ to $t_1$

$$\int |u_n(t_0, x)|^2 g\left(\frac{|x| - K_2}{\gamma_1}\right) dx - \int |u_n(t_1, x)|^2 g\left(\frac{|x| - K_2}{\gamma_1}\right) dx$$
$$\leq \int_{t_0}^{t_1} \left| \frac{d}{dt} \int |u_n(t, x)|^2 g\left(\frac{|x| - K_2}{\gamma_1}\right) dx \right| \leq \delta_1.$$

From the properties of $g$ we conclude:

$$\int_{|x| > 2\gamma_1 + K_2} |u_n(t_0, x)|^2 dx \leq \int |u_n(t_0, x)|^2 g\left(\frac{|x| - K_2}{\gamma_1}\right) dx$$
$$\leq \int |u_n(t_1, x)|^2 g\left(\frac{|x| - K_2}{\gamma_1}\right) dx + \delta_1 \leq \int_{|x| > K_2} |u_n(t_1, x)|^2 dx + \delta_1 \leq 5\delta_1.$$

Thus (4.4) is proved. As $||u_n(t_0)||_{H^1} < C$, there exists a subsequence of $(u_n)$ (still denoted by $(u_n)$) and $u_0 \in H^1$ such that

$$u_n(t_0) \rightharpoonup u_0 \quad \text{weakly in } H^1(\mathbb{R}^2),$$

$$u_n(t_0) \to u_0 \quad \text{in } L^2_{loc}(\mathbb{R}^d), \text{ as } n \to +\infty$$

and by (4.4), we conclude that $u_n(t_0) \xrightarrow{L^2} u_0$ as required. $\qquad \square$

Let us finish the proof of the Main Theorem in sub-critical cases with $p > 2$. We consider $u$ the solution of (NLS) corresponding to $u(t_0) = u_0$. By continuous dependence of the solution upon the initial data (see [2] and [3]), for all $0 \leq \sigma < 1$, for all $t \in [t_0, +\infty)$,

$$u_n(t) \to u(t) \quad \text{in } H^\sigma(\mathbb{R}^d).$$

Moreover, the decomposition $(\vec{q}, \epsilon)$ of $u$ satisfies, for all $t \geq t_0$,

$$\vec{q}_n(t) \to \vec{q}(t), \quad \epsilon_n(t) \to \epsilon(t) \text{ in } H^\sigma, \quad \epsilon_n(t) \rightharpoonup \epsilon(t) \text{ in } H^1 \qquad (4.7)$$

(see e.g. [21], Claim p.598). In particular, for all $t \in [t_0, +\infty)$, $u(t)$ decomposes as

$$u(t,x) = \frac{e^{i\gamma(t)}}{\lambda^{\frac{2}{p-1}}(t)} \left( \sum_{k=1}^{2} \left[ e^{i\Gamma_k} Q \right] \left( \frac{x + \frac{(-1)^k}{2} \lambda(t) z(t)}{\lambda(t)} \right) + \epsilon \left( t, \frac{x}{\lambda(t)} \right) \right), \qquad (4.8)$$

where $\Gamma_k(t,y) = \frac{(-1)^{k+1}}{2} v(t) \cdot y$ and it follows from the uniform estimates (4.2) that

$$| |z(t)| - 2\log t | \lesssim \log(\log t), \quad \left| \lambda^{-1}(t) - 1 \right| \lesssim t^{-1},$$

$$|v(t)| \lesssim t^{-1}, \quad \|\epsilon(t)\|_{H^1} \lesssim t^{-1}, \quad \left| |z(t)|^{\frac{d-1}{2}} e^{|z(t)|} - ct^2 \right| \lesssim t^2 \log^{-\frac{1}{2}}(t). \qquad (4.9)$$

We obtain $|\delta x(t)| = |x_1(t) - x_2(t)| = \lambda(t)|z(t)| \to 2(1 + o(1))\log t$, more precisely

$$| |\delta x(t)| - 2\log(t) | \lesssim \log(\log(t))$$

and the following estimate

$$\left\| u(t) - e^{i\gamma(t)} \sum_{k=1}^{2} Q \left( x + \frac{(-1)^k}{2} \delta x(t) \right) \right\|_{H^1} \lesssim \|\epsilon(t)\|_{H^1} + |\lambda^{-1}(t) - 1| + |v(t)| \lesssim t^{-1}. \quad (4.10)$$

## 5. Sub-critical cases with $1 < p \leq 2$

In this section, we show the difficulties occurring and sketch the proof of Main Theorem in the case $1 < p \leq 2$. In this case, let

$$2^+ = \min(2^*, \frac{p+3}{2}).$$

Note that $p - 2^+ > -1$. From (2.15), we deduce the following Taylor expansions:

$$F(\mathbf{P} + \epsilon) = F(\mathbf{P}) + F'(\mathbf{P}).\epsilon + O(|\epsilon|^p) \qquad (5.1)$$

$$F(\mathbf{P} + \epsilon) = F(\mathbf{P}) + F'(\mathbf{P}).\epsilon + O\left( \left| \frac{\epsilon}{\mathbf{P}} \right|^2 |\mathbf{P}|^p \right) \qquad (5.2)$$

(since $|\epsilon| > \frac{|\mathbf{P}|}{2}$ then $|\epsilon|^p \lesssim |\frac{\epsilon}{\mathbf{P}}|^2 |\mathbf{P}|^p$ and $|\epsilon| \leq \frac{|\mathbf{P}|}{2}$ then $|\frac{\epsilon}{\mathbf{P}}|^2 |\mathbf{P}|^p \lesssim |\epsilon|^p$) and

$$F(\mathbf{P} + \epsilon) = F(\mathbf{P}) + F'(\mathbf{P}).\epsilon + \frac{\bar{\epsilon}.F''(\mathbf{P}).\epsilon}{2} + O\left( \left| \frac{\epsilon}{\mathbf{P}} \right|^{2^+} |\mathbf{P}|^p \right). \qquad (5.3)$$

In the following remark, we identify new problems compared with the case $p > 2$.

**Remark 7.** *Let us try to control the nonlinear interaction term*

$$G(y; (z(s), v(s))) = |\mathbf{P}|^{p-1}\mathbf{P} - |P_1|^{p-1}P_1 - |P_2|^{p-1}P_2.$$

*Since $|P_1| > |P_2|$ for $y \cdot \frac{z}{|z|} > 0$ and $|P_2| > |P_1|$ for $y \cdot \frac{z}{|z|} < 0$, one has by (2.15)*

$$|G(y; (z(s), v(s)))| = \left||P_1 + P_2|^{p-1}(P_1 + P_2) - |P_1|^{p-1}P_1 - |P_2|^{p-1}P_2\right| \atop \lesssim |P_1|^{p-1}|P_2|.\mathbb{1}_{y \cdot \frac{z}{|z|} > 0} + |P_2|^{p-1}|P_1|.\mathbb{1}_{y \cdot \frac{z}{|z|} < 0}. \tag{5.4}$$

*Using the asymptotic behavior of $Q$, on the half space $\{y \cdot \frac{z}{|z|} > 0\}$*

$$|P_1|^{p-1}|P_2|.\mathbb{1}_{y \cdot \frac{z}{|z|} > 0} \lesssim |P_1 P_2|^{p-1}|P_2|^{2-p}.\mathbb{1}_{y \cdot \frac{z}{|z|} > 0} \lesssim |z|^{-\frac{(p-1)(d-1)}{2}}e^{-(p-1)|z|}|P_2|^{2-p}.\mathbb{1}_{y \cdot \frac{z}{|z|} > 0}$$

$$\lesssim |z|^{-\frac{(p-1)(d-1)}{2}}e^{-(p-1)|z|}\left|\frac{z}{2}\right|^{-\frac{(2-p)(d-1)}{2}}e^{-\frac{2-p}{2}|z|} \lesssim |z|^{-\frac{d-1}{2}}e^{-\frac{p}{2}|z|}. \tag{5.5}$$

*By symmetry, we have the same estimate on the other half space $\{y \cdot \frac{z}{|z|} < 0\}$ and thus*

$$\|G\|_{L^\infty} \lesssim |z|^{-\frac{d-1}{2}}e^{-\frac{p}{2}|z|} \sim s^{-p} \tag{5.6}$$

*(to be compared with (2.9)). Now for the projection of interaction, we recall that its core part (as identified in the proof of Lemma 2 and in step 4 of Proposition 9) is given by*

$$H(z) = p\int_{y \cdot \frac{z}{|z|} > -\frac{|z|}{2}} Q^{p-1}(y)\nabla Q(y)Q(y + z)dy + p\int_{y \cdot \frac{z}{|z|} < -\frac{|z|}{2}} Q^{p-1}(y + z)\nabla Q(y)Q(y)dy$$

*and the following estimate of $H(z)$ is still valid for $1 < p \leq 2$ (see Lemma 2)*

$$\left|H(z) - c_Q I_Q \frac{z}{|z|}|z|^{-\frac{d-1}{2}}e^{-|z|}\right| \lesssim |z|^{-1-\frac{d-1}{2}}e^{-|z|}. \tag{5.7}$$

*In summary, the projection $\langle G, e^{i\Gamma_1}\nabla Q(y - z_1(s))\rangle$ and thus $\dot{v}$ are still of order $s^{-2}$, however the interaction $G$ is of order $s^{-p} \gg s^{-2}$ in $L^\infty$ norm. Therefore, there still exist some terms in the interaction that perturb our regime and prevent us to close the bootstrap arguments (for example (3.40)).*

In view of the above remark, we look for a refined approximate solution $\mathbf{P}$ of the form

$$\mathbf{P}(s, y) = \mathbf{P}(y; (z(s), v(s))) = \sum_{k=1}^{2} e^{iv_k(s)(y - z_k(s))}Q(y - z_k(s)) + W(y; (z(s), v(s))) \atop = \sum_{k=1}^{2} P_k(s, y) + W(y; (z(s), v(s))), \tag{5.8}$$

where $W(y; (z(s), v(s)))$ to be determined.

**Proposition 9** (Expansion of the refined approximate solution). *There exists a series of $(J + 1)$ functions $R_j(y; (z(s), v(s)))$ which are invariant by $\tau$ such that by setting*

$$W(y; (z(s), v(s))) = \sum_{j=0}^{J} R_j(y; (z(s), v(s))),$$

*the error $\mathcal{E}_{\mathbf{P}}$ defined as in (2.7) admits the decomposition*

$$\mathcal{E}_{\mathbf{P}} = [e^{i\Gamma_1}\vec{m}_1 \cdot \vec{\mathbf{M}}Q](y - z_1(s)) + [e^{i\Gamma_2}\vec{m}_2 \cdot \vec{\mathbf{M}}Q](y - z_2(s)) + G_0, \tag{5.9}$$

*where under the bootstrap assumptions (3.21) and the pointwise control of the modulation equation (3.24)–(3.26)*

$$|z| \lesssim \log(s), \quad |\dot{z}| \lesssim s^{-1}, \quad |v| \lesssim s^{-1}, \quad |\dot{v}| \lesssim s^{-2}, \quad \left|\frac{\dot{\lambda}}{\lambda}\right| \lesssim (C^*)^2 s^{-2}, \quad |\dot{\gamma} - 1| \lesssim (C^*)^2 s^{-2},$$

*the corrected interaction term $G_0$ satisfies*

$$\|G_0\|_{L^2} \lesssim s^{-2}, \quad \|\nabla G_0\|_{L^2} \lesssim s^{-2}. \tag{5.10}$$

*Moreover, $G_0$ is symmetric and*

$$\left|\langle G_0, e^{i\Gamma_1(y-z_1(s))}\nabla Q(y - z_1(s))\rangle - C_p \frac{z}{|z|}|z|^{-\frac{d-1}{2}}e^{-|z|}\right| \lesssim s^{-2}\log^{-1}(s) \tag{5.11}$$

*with $C_p > 0$.*

**Remark 8.** *In fact, before the pointwise control of the modulation equations in Lemma 5, we bound $\|G_0\|_{L^2}, \|\nabla G_0\|_{L^2}$ by $z, v$ and $s^{-p}|\vec{m}_1|$ then once we have the control on $\vec{m}_1$, we will obtain (5.10).*

*Proof of Proposition 9.* **step 1** Properties of the Helmholtz operators. We recall well-known properties of $(-\Delta + 1)u_s(y) = f_s(y)$ in $\mathbb{R}^d$. The operator $(-\Delta + 1)^{-1}$ is continuous from $L^2$ to $H^1$, in particular

$$\|u\|_{H^1} \leq \|f\|_{L^2}.$$

It is self-adjoint

$$\langle u, (-\Delta + 1)g\rangle = \langle(-\Delta + 1)u, g\rangle = \langle f, g\rangle, \tag{5.12}$$

invariant by $\tau$ and $(-\Delta + 1)\dot{u}_s(y) = \dot{f}_s(y)$ ($\dot{f}$ denotes the derivative with respect to time $s$). Moreover, by theory of elliptic equation (see e.g [1]), we have an explicit kernel representation $E_d$ for $(-\Delta + 1)^{-1}$ as follows

$$E_d(x) = -(2\pi)^{-\frac{d}{2}}\left(\frac{1}{|x|}\right)^{\frac{d}{2}-1}\mathcal{K}_{\frac{d}{2}-1}(|x|)$$

$$u(x) = \int_{\mathbb{R}^d} E_d(x - y)f(y)dy \tag{5.13}$$

where $\mathcal{K}_\alpha$ is modified Bessel functions of second kind which is decreasing exponentially when $|x| \to +\infty$. This is a convolution of type $L^1 \star L^\infty$ so we deduce that

$$\|u\|_{L^\infty} \lesssim \|f\|_{L^\infty}. \tag{5.14}$$

Next, we claim the exponential decay property: assume that a regular function $f$ is exponentially decreasing in the direction $e_j$, $e^{\delta|y_j|}|f(y)| \leq C$ with $0 < \delta < 1$, then so is the solution $u$ of $(-\Delta + 1)^{-1}$.

Indeed, we consider

$$e^{\delta|x_j|}|u(x)| = e^{\delta|x_j|}\left|\int_{\mathbb{R}^d} E_d(x - y)f(y)dy\right|$$

$$\lesssim C\left|\int_{\mathbb{R}^d}\left(\frac{1}{|x-y|}\right)^{\frac{d}{2}-1}e^{-|x-y|}e^{\delta(|x_j|-|y_j|)}dy\right| \lesssim C\left\|\left(\frac{1}{|x|}\right)^{\frac{d}{2}-1}e^{-(1-\delta)|x|}\right\|_{L^1} \lesssim C.$$

**step 2** Iteration of $R_j$. We introduce a suitable smooth cut-off function that localizes the points whose distances to center of two solitons are smaller than $|z|$. Denote $\psi_0 : \mathbb{R} \to [0,1]$ such that

$$0 \leq \psi_0' \leq C, \quad \psi_0 \equiv 0 \text{ on } (-\infty, -1], \quad \psi_0 \equiv 1 \text{ on } [0, +\infty)$$

and

$$\psi(y; z(s)) = \psi_0\left(|z(s)| - \left|y + \frac{z(s)}{2}\right|\right) \psi_0\left(|z(s)| - \left|y - \frac{z(s)}{2}\right|\right).$$

Recall the definition of $G$

$$G(y; (z(s), v(s))) = |P_1 + P_2|^{p-1}(P_1 + P_2) - |P_1|^{p-1}P_1 - |P_2|^{p-1}P_2$$

and denote $\mathfrak{pr}_i$ the projection on the direction $\nabla Q$ around each soliton

$$\mathfrak{pr}_i(f) = \frac{\langle f(\cdot), \nabla Q(\cdot + \frac{(-1)^i}{2} z(s)) \rangle}{\|\nabla Q(\cdot + \frac{(-1)^i}{2} z(s))\|_{L^2}^2} \nabla Q(\cdot + \frac{(-1)^i}{2} z(s)).$$

Setting

$$A_0(y; (z(s), v(s))) = G(y; (z(s), v(s)))\psi(y; z(s)),$$

$$\tilde{A}_0 = A_0 - \mathfrak{pr}_1(A_0) - pr_2(A_0),$$

$$A_1 = |P_1 + P_2 + R_0|^{p-1}(P_1 + P_2 + R_0) - |P_1 + P_2|^{p-1}(P_1 + P_2),$$

$$\tilde{A}_1 = A_1 - \mathfrak{pr}_1(A_1) - pr_2(A_1)$$

and for $j \geq 2$

$$A_j = |P_1 + P_2 + \sum_{k=0}^{j-1} R_k|^{p-1}(P_1 + P_2 + \sum_{k=0}^{j-1} R_k) - |P_1 + P_2 + \sum_{k=0}^{j-2} R_k|^{p-1}(P_1 + P_2 + \sum_{k=0}^{j-2} R_k),$$

$$\tilde{A}_j = A_j - \mathfrak{pr}_1(A_j) - pr_2(A_j).$$

Observe that

$$\sum_{j=1}^{J} A_j = |P_1 + P_2 + \sum_{k=0}^{j-1} R_k|^{p-1}(P_1 + P_2 + \sum_{k=0}^{j-1} R_k) - |P_1 + P_2|^{p-1}(P_1 + P_2). \tag{5.15}$$

Then let

$$R_j(y; (z(s), v(s))) = (-\Delta + 1)^{-1} \tilde{A}_j.$$

We will show by induction on $j$ the following properties

- $R_j$ is almost orthogonal to $\nabla(Q^p)(\cdot \pm \frac{1}{2}z)$, i.e,

$$\langle R_j(\cdot), \nabla(Q^p)(\cdot \pm \frac{1}{2}z) \rangle \lesssim s^{-3}. \tag{5.16}$$

- The $L^\infty, H^1$ norm of $R_j$ satisfy

$$\|R_{j+1}\|_{L^\infty} \lesssim s^{-(p-1)}\|R_j\|_{L^\infty} \lesssim s^{-p},$$

$$\|R_{j+1}\|_{H^1} \lesssim s^{-(p-1-\kappa)}\|R_j\|_{H^1} \lesssim s^{-p}\log^{dp}(s)$$

with $0 < \kappa \ll 1$ to be determined (see (5.33), (5.34)).

- After a finite number $(J+1)$ of steps, the function $R_J$ satisfies the two following estimates: there is $\epsilon > 0$

$$|Q^{p-1}(y)R_J(y+\tfrac{z}{2})| + |Q^{p-1}(y)R_J(y+\tfrac{z}{2})| \lesssim e^{-\epsilon|y|}s^{-2} \tag{5.17}$$

$$\|R_J\|_{H^1}^p + s^{p(p-1)}\|R_J\|_{H^1} \ll s^{-2} \tag{5.18}$$

independently of $z, v$ ((5.18) means thats there exists $\delta > 0$ such that $\|R_J\|_{H^1}^p + s^{-p(p-1)}\|R_J\|_{H^1} \lesssim s^{-2-\delta}$).

Note that a direct consequence of the above estimates is

$$\|A_{J+1}\|_{L^2} = \left\| |P_1 + P_2 + \sum_{j=0}^{J} R_j|^{p-1}(P_1 + P_2 + \sum_{j=0}^{J} R_j) - |P_1 + P_2 + \sum_{j=0}^{J-1} R_j|^{p-1}(P_1 + P_2 + \sum_{j=0}^{J-1} R_j) \right\|_{L^2}$$

$$\lesssim \left\| |P_1 + P_2 + \sum_{j=0}^{J-1} R_j|^{p-1}|R_J| + |R_J|^p \right\|_{L^2} \lesssim \|Q^{p-1}(\cdot)R_J(\cdot+\tfrac{z}{2})\|_{L^2} + \|R_J\|_{L^2}^p + s^{p(p-1)}\|R_J\|_{L^2} \lesssim s^{-2} \tag{5.19}$$

since $\|R_j\|_{L^\infty} \lesssim \|R_0\|_{L^\infty} \lesssim s^{-p}, \forall j = \overline{1,J}$.

Let us begin with $R_0$. We have that

$$|G(y;(z(s),v(s)))| = \left| |P_1 + P_2|^{p-1}(P_1 + P_2) - |P_1|^{p-1}P_1 - |P_2|^{p-1}P_2 \right|$$
$$\lesssim |P_1|^{p-1}|P_2|.\mathbb{1}_{y\cdot\frac{z}{|z|}>0} + |P_2|^{p-1}|P_1|.\mathbb{1}_{y\cdot\frac{z}{|z|}<0}.$$

Consider

$$|P_1|^{p-1}|P_2|.\mathbb{1}_{y\cdot\frac{z}{|z|}>0} \lesssim e^{-(p-1)(z_1-y\cdot\frac{z}{|z|})}\left|\frac{z}{2}\right|^{-\frac{d-1}{2}}e^{-(y\cdot\frac{z}{|z|}-z_2)}$$
$$\lesssim |z|^{-\frac{d-1}{2}}e^{-\frac{p}{2}|z|}e^{-(2-p)y\cdot\frac{z}{|z|}} \lesssim s^{-p}e^{-(2-p)y\cdot\frac{z}{|z|}}|z|^{-\frac{(2-p)(d-1)}{2}}, \tag{5.20}$$

by symmetry, we also have the same estimate on $\{y\cdot\frac{z}{|z|}<0\}$. Thus, from definition of $\psi$, we get

$$\|e^{(2-p)|y\cdot\frac{z}{|z|}|}A_0(y;(z(s),v(s)))\|_{L^\infty} \lesssim s^{-p}|z|^{-\frac{(2-p)(d-1)}{2}} \lesssim s^{-p} \tag{5.21}$$

and

$$\|A_0(y;(z(s),v(s)))\|_{L^2} \lesssim s^{-p}\log^d(s). \tag{5.22}$$

The estimate (5.21) yields

$$|A_0(y+\tfrac{z}{2})| \lesssim e^{-(2-p)|y\cdot\frac{z}{|z|}+\frac{z}{|z|}|}s^{-p}|z|^{-\frac{(2-p)(d-1)}{2}}$$
$$\lesssim e^{(2-p)|y\cdot\frac{z}{|z|}|}e^{-(2-p)\frac{|z|}{2}}|z|^{-\frac{(2-p)(d-1)}{2}}s^{-p} \lesssim e^{(2-p)|y\cdot\frac{z}{|z|}|}s^{-2}$$

so it gives a control on projections of $A_0$

$$\left|\int_{\mathbb{R}^d} A_0(y+\tfrac{z}{2})\nabla Q(y)dy\right| \lesssim s^{-2}. \tag{5.23}$$

Therefore, from definition of $\tilde{A}_0$

$$\|e^{(2-p)|y\cdot\frac{z}{|z|}|}\tilde{A}_0\|_{L^\infty} \lesssim s^{-p}|z|^{-\frac{(2-p)(d-1)}{2}}, \quad \|\tilde{A}_0\|_{L^2} \lesssim s^{-p}\log^d(s).$$

From step 1, we can transfer these properties to $R_0(y; (z(s), v(s)))$

$$\|e^{(2-p)|y \cdot \frac{z}{|z|}|} R_0(y; (z(s), v(s)))\|_{L^\infty} \lesssim s^{-p}|z|^{-\frac{(2-p)(d-1)}{2}}, \tag{5.24}$$

$$\|R_0(y; (z(s), v(s)))\|_{H^1} \lesssim s^{-p}\log^d(s). \tag{5.25}$$

To show the almost orthogonality condition, we note that $(-\Delta + 1)\nabla Q = \nabla(Q^p)$ so from self-adjoint property (5.12) of $(-\Delta + 1)$, we have

$$\left|\langle R_0, \nabla(Q^p)(\cdot + \frac{z}{2})\rangle\right| = \left|\langle A_0 - \mathfrak{pr}_1(A_0) - \mathfrak{pr}_2(A_0), \nabla Q(\cdot + \frac{z}{2})\rangle\right|$$

$$= \left|\langle \mathfrak{pr}_1(A_0), \nabla Q(\cdot + \frac{z}{2})\rangle\right| \lesssim s^{-2}\langle \nabla Q(\cdot - \frac{z}{2}), \nabla Q(\cdot + \frac{z}{2})\rangle \lesssim s^{-3}.$$

If $\frac{3}{2} < p \le 2$, we see that $R_0$ satisfies already the conditions (5.17), (5.18) as

$$\|R_0\|_{H^1}^p \lesssim s^{-p^2}\log^{dp}(s) \le s^{-\frac{9}{4}}\log^{dp}(s) \ll s^{-2}$$

$$s^{p(p-1)}\|R_0\|_{H^1} \lesssim s^{-\frac{3}{4}}s^{-\frac{3}{2}} \ll s^{-2}$$

and $|R_0(y + \frac{z}{2})| \lesssim e^{-(2-p)|y \cdot \frac{z}{|z|} + \frac{z}{2}|}s^{-p}|z|^{-\frac{(2-p)(d-1)}{2}} \lesssim e^{(2-p)|y \cdot \frac{z}{|z|}|}s^{-2}$ so for $\epsilon = 2p - 3 > 0$

$$|Q^{p-1}(y)R_0(y + \frac{z}{2})| \lesssim e^{(2-p)|y \cdot \frac{z}{|z|}|}Q^{(2-p)}(y)s^{-2}Q^{(2p-3)}(y) \lesssim e^{-\epsilon|y|}s^{-2}.$$

Thus $J = 0$ and $W = R_0(y; (z(s), v(s)))$ in this case.

If $\frac{4}{3} < p \le \frac{3}{2}$, we consider $A_1(y; (z(s), v(s)))$, by (2.15), we obtain

$$\left||P_1 + P_2 + R_0|^{p-1}(P_1 + P_2 + R_0) - |P_1 + P_2|^{p-1}(P_1 + P_2)\right.$$

$$\left. - \frac{p+1}{2}|P_1 + P_2|^{p-1}R_0 - \frac{p-1}{2}|P_1 + P_2|^{p-3}(P_1 + P_2)^2\overline{R_0}\right| \lesssim |R_0|^p. \tag{5.26}$$

Next remark that for $1 < p \le 2$, $\left||P_1 + P_2|^{p-1} - |P_1|^{p-1} - |P_2|^{p-1}\right| \lesssim \min(|P_1|^{p-1}, |P_2|^{p-1})$ so the main part of $A_1 = |P_1 + P_2 + R_0|^{p-1}(P_1 + P_2 + R_0) - |P_1 + P_2|^{p-1}(P_1 + P_2)$ can be computed by

$$\left\|\frac{p+1}{2}|P_1 + P_2|^{p-1}R_0 + \frac{p-1}{2}|P_1 + P_2|^{p-3}(P_1 + P_2)^2\overline{R_0} - p|P_1 + P_2|^{p-1}R_0\right\|_{L^2}$$

$$\lesssim \left\|(|v|^2|y|^2 + |v|^2|z|^2)|R_0|(|P_1|^{p-1} + |P_2|^{p-1})\right\|_{L^2} \ll s^{-2} \tag{5.27}$$

$$\left||P_1 + P_2|^{p-1}R_0 - (|P_1|^{p-1} + |P_2|^{p-1})R_0\right| \lesssim \min(|P_1|^{p-1}, |P_2|^{p-1})|R_0| \tag{5.28}$$

here in (5.27) we use the bootstrap assumptions and the control of modulation equations. Let estimate $R_0(y)Q^{p-1}(y + \frac{z}{2})$, from the decreasing properties of $R_0$ (5.24), we have

$$|R_0(y)Q^{p-1}(y + \frac{z}{2})| \lesssim e^{-(2-p)|y \cdot \frac{z}{|z|}|}s^{-p}|z|^{-\frac{(2-p)(d-1)}{2}}e^{(p-1)|y \cdot \frac{z}{|z|}|}e^{-(p-1)\frac{|z|}{2}}$$

$$\lesssim e^{-(3-2p)|y \cdot \frac{z}{|z|}|}s^{-(2p-1)}|z|^{-\frac{(3-2p)(d-1)}{2}} \tag{5.29}$$

so for $\kappa \ll 1$ determined later in (5.33)

$$\|R_0(y)Q^{p-1}(y + \frac{z}{2})\|_{L^2} \lesssim s^{-(2p-1-\kappa)}. \tag{5.30}$$

The collection of above estimates gives a bound on norm $L^2$ and on the decay property of $A_1$

$$\|A_1\|_{L^2} \lesssim \|R_0\|_{L^2}^p + \|R_0|P_2|^{p-1}\|_{L^2} + \|R_0|P_2|^{p-1}\|_{L^2}$$
$$\lesssim s^{-p^2}\log^{dp}(s) + s^{-(2p-1-\kappa)} \leq s^{-(2p-1-\kappa)},$$

$$\|e^{(3-2p)|y \cdot \frac{z}{|z|}|}A_1\|_{L^\infty} \lesssim \|e^{(3-2p)|y \cdot \frac{z}{|z|}|}|R_0|^p\|_{L^\infty} + \|e^{(3-2p)|y \cdot \frac{z}{|z|}|}R_0|P_2|^{p-1}\|_{L^\infty} + \|e^{(3-2p)|y \cdot \frac{z}{|z|}|}R_0|P_2|^{p-1}\|_{L^\infty}$$
$$\lesssim s^{-p^2} + s^{-(2p-1)}|z|^{-\frac{(3-2p)(d-1)}{2}} \leq s^{-(2p-1)}|z|^{-\frac{(3-2p)(d-1)}{2}}$$

as the decay $e^{-(2-p)|y \cdot \frac{z}{|z|}|}$ of $R_0$ is faster than the one of $e^{-(3-2p)|y \cdot \frac{z}{|z|}|}$. Finally, we consider

$$\left| \langle A_1, \nabla Q(y + \frac{z}{2}) \rangle - p \left\langle Q^{p-1}(y - \frac{z}{2})R_0 + Q^{p-1}(y + \frac{z}{2})R_0, \nabla Q(y + \frac{z}{2}) \right\rangle \right|$$
$$\lesssim \left\langle |R_0|^p, \nabla Q(y + \frac{z}{2}) \right\rangle + \left\langle \min(|P_1|^{p-1}, |P_2|^{p-1})|R_0|, \nabla Q(y + \frac{z}{2}) \right\rangle$$
$$\lesssim \left\langle e^{-(2-p)p|y \cdot \frac{z}{|z|}|}s^{-p^2}|z|^{-\frac{(2-p)p(d-1)}{2}}e^{(2-p)|y \cdot \frac{z}{|z|}|}e^{-(2-p)p\frac{|z|}{2}}, Q^{1-(2-p)p}(y + \frac{z}{2}) \right\rangle$$
$$+ \left\langle s^{-(p-1)}e^{-(2-p)|y \cdot \frac{z}{|z|}|}s^{-p}|z|^{-\frac{(2-p)(d-1)}{2}}e^{(2-p)|y \cdot \frac{z}{|z|}|}e^{-(2-p)\frac{|z|}{2}}, Q^{1-(2-p)}(y + \frac{z}{2}) \right\rangle$$
$$\lesssim s^{-2p} + s^{-(p+1)} \ll s^{-2}.$$

We can deduce from the almost orthogonality (5.16) that

$$\langle A_1, \nabla Q(y \pm \frac{z}{2}) \rangle \ll s^{-2}, \tag{5.31}$$

in other words, we have

$$\|\mathfrak{pr}_i(A_1)\|_{L^2} \ll s^{-2}, \qquad i = 1, 2. \tag{5.32}$$

Therefore, we have the following estimates for $\tilde{A}_1 = A_1 - \mathfrak{pr}_1(A_1) - \mathfrak{pr}_2(A_2)$

$$\|\tilde{A}_1\|_{L^2} \lesssim s^{-(2p-1-\kappa)}, \quad \|e^{(3-2p)|y \cdot \frac{z}{|z|}|}\tilde{A}_1\|_{L^\infty} \lesssim s^{-(2p-1)}|z|^{-\frac{(3-2p)(d-1)}{2}}$$

and the analogue for $R_1$

$$\|R_1\|_{H^1} \lesssim s^{-(2p-1-\kappa)}, \quad \|e^{(3-2p)|y \cdot \frac{z}{|z|}|}R_1\|_{L^\infty} \lesssim s^{-(2p-1)}|z|^{-\frac{(3-2p)(d-1)}{2}}.$$

There exists $0 < \kappa \ll 1$ such that for all $p > \frac{4}{3}$

$$-(2p-1-\kappa)p < -2, \quad -(2p-1-\kappa) - p(p-1) < -2 \tag{5.33}$$

so $\|R_J\|_{H^1}^p + s^{p(p-1)}\|R_J\|_{H^1} \lesssim s^{-(2p-\kappa)p} + s^{-(2p-1-\kappa)-p(p-1)} \ll s^{-2}$ and for $\epsilon = 3p - 4 > 0$

$$|Q^{p-1}(y)R_1(y + \frac{z}{2})| \lesssim e^{-(3-2p)|y \cdot \frac{z}{|z|} + \frac{z}{2}|}s^{-(2p-1)}|z|^{-\frac{(3-2p)(d-1)}{2}}Q^{p-1}(y)$$
$$\leq e^{(3-2p)|y \cdot \frac{z}{|z|}|}Q^{(3-2p)}(y)s^{-2}Q^{(3p-4)}(y) \lesssim e^{-\epsilon|y|}s^{-2}.$$

The almost orthogonal property of $V_1$ is a direct consequence of $\langle \tilde{A}_1(\cdot \pm \frac{z}{2}), \nabla Q \rangle \lesssim s^{-3}$. Thus $J = 1$ and $W = R_0(y; (z(s), v(s))) + R_0(y; (z(s), v(s)))$ in this case.

If $\frac{J+3}{J+2} < p \leq \frac{J+2}{J+1}$, we proceed the same way and after $(J+1)$ steps, our process will finish with

$$W = \sum_{j=0}^{J} R_j(y; (z(s), v(s))),$$

$\epsilon = (J+2)p - (J+3) > 0$ and $0 < \kappa \ll 1$ such that for all $\frac{J+2}{J+1} < p \leq \frac{J+1}{J}$

$$-((J+1)p - J - \kappa)p < -2, \quad -((J+1)p - J - \kappa) - p(p-1) < -2. \quad (5.34)$$

**step 3** Estimate of $G_0$. Let $\mathbf{P} = P_1 + P_2 + W$ and put into the definition $\mathcal{E}_{\mathbf{P}}$, it follows from the computations in Lemma 1 that

$$\mathcal{E}_{\mathbf{P}} = [e^{i\Gamma_1}\vec{m}_1 \cdot \vec{\mathbb{M}}Q](y - z_1(s)) + [e^{i\Gamma_2}\vec{m}_2 \cdot \vec{\mathbb{M}}Q](y - z_2(s)) + |\mathbf{P}|^{p-1}\mathbf{P} - |P_1|^{p-1}P_1 - |P_1|^{p-1}P_1$$

$$+ \sum_{j=0}^{J}(\Delta - 1)R_j + \sum_{j=0}^{J}[i\dot{R}_j - i\frac{\dot{\lambda}}{\lambda}\Lambda R_j + (1 - \dot{\gamma})R_j]. \quad (5.35)$$

Note that

$$\sum_{j=1}^{J}(\Delta - 1)R_j = -\sum_{j=1}^{J}\tilde{A}_j = -\sum_{j=1}^{J}A_j + \sum_{j=1}^{J}[\mathfrak{pr}_1(A_j) + \mathfrak{pr}_2(A_j)]$$

thus following (5.9) and (5.15), we have the explicit expression of $G_0$

$$G_0 = |P_1 + P_2 + \sum_{j=0}^{J}R_j|^{p-1}(P_1 + P_2 + \sum_{j=0}^{J}R_j) - |P_1 + P_2 + \sum_{j=0}^{J-1}R_j|^{p-1}(P_1 + P_2 + \sum_{j=0}^{J-1}R_j)$$

$$+ \sum_{j=1}^{J}[A_j + (\Delta - 1)R_j] + |P_1 + P_2|^{p-1}(P_1 + P_2) - |P_1|^{p-1}P_1 - |P_2|^{p-1}P_2 - (\Delta - 1)R_0$$

$$+ \sum_{j=0}^{J}[i\dot{R}_j - i\frac{\dot{\lambda}}{\lambda}\Lambda R_j + (1 - \dot{\gamma})R_j]$$

$$= |P_1 + P_2 + \sum_{j=0}^{J}R_j|^{p-1}(P_1 + P_2 + \sum_{j=0}^{J}R_j) - |P_1 + P_2 + \sum_{j=0}^{J-1}R_j|^{p-1}(P_1 + P_2 + \sum_{j=0}^{J-1}R_j)$$

$$+ \sum_{j=1}^{J}[\mathfrak{pr}_1(A_j) + \mathfrak{pr}_2(A_j)] + G + (\Delta - 1)R_0 + \sum_{j=0}^{J}[i\dot{R}_j - i\frac{\dot{\lambda}}{\lambda}\Lambda R_j + (1 - \dot{\gamma})R_j]$$

$$= A_{J+1} + \sum_{j=1}^{J}[\mathfrak{pr}_1(A_j) + \mathfrak{pr}_2(A_j)] + \mathfrak{pr}_1(G\psi) + \mathfrak{pr}_2(G\psi) + G(1 - \psi) + \sum_{j=0}^{J}[i\dot{R}_j - i\frac{\dot{\lambda}}{\lambda}\Lambda R_j + (1 - \dot{\gamma})R_j].$$

We bound the first term by (5.19)

$$\|A_{J+1}\|_{L^2} \lesssim s^{-2}.$$

Next, from pointwise control of the modulation equations, we have $\left|\frac{\dot{\lambda}}{\lambda}\right|, |1 - \dot{\gamma}| \lesssim (C^*)^2 s^{-2}$ and $\|R_j\|_{H^1} < \|R_0\|_{H^1} \lesssim s^{-p}\log^d(s)$, therefore

$$\left\|\sum_{j=0}^{J} i\frac{\dot{\lambda}}{\lambda}\Lambda R_j - (1 - \dot{\gamma})R_j\right\|_{L^2} \ll s^{-2}. \quad (5.36)$$

We recall (5.23) that

$$\|\mathfrak{pr}_1(G\psi)\|_{L^2} + \|\mathfrak{pr}_2(G\psi)\|_{L^2} \lesssim s^{-2}$$

and similarly to (5.32), we have

$$\|\mathfrak{pr}_1(A_j)\|_{L^2} + \|\mathfrak{pr}_2(A_j)\|_{L^2} \ll s^{-2}, \qquad \forall j \geq 1.$$

The term

$$\|G(1-\psi)\|_{L^2} \lesssim |z|^{-\frac{d-1}{2}} e^{-|z|} (\|P_1\|_{L^2}^{p-1} + \|P_2\|_{L^2}^{p-1}) \lesssim s^{-2},$$

this is a consequence of the choice of localized cut-off function $\psi$ and the decay property of $Q$. For the last term, we have $\dot{R}_j = (-\Delta+1)^{-1}\dot{\tilde{A}}_j$, so

$$\|\dot{R}_j\|_{H^1} \leq \|\dot{\tilde{A}}_j\|_{L^2}.$$

We consider $R_0$ and $A_0$, proceeding as the way we control $G$ in (5.4), we have that $\dot{G}$ decays more rapidly because of extra terms $\dot{z}$ and $\dot{v}$. In fact, we have

$$|\dot{G}| \leq \left| (\dot{P}_1 + \dot{P}_2)|P_1 + P_2|^{p-1} - \dot{P}_1|P_1|^{p-1} - \dot{P}_2|P_2|^{p-1} \right|$$

$$+ \left| (\dot{P}_1 + \dot{P}_2)|P_1 + P_2|^{p-2}(P_1 + P_2) - \dot{P}_1|P_1|^{p-2}P_1 - \dot{P}_2|P_2|^{p-2}P_2 \right|$$

and

$$\dot{P}_k = \dot{z}_k \nabla P_k + i\dot{v}_k(y - z_k)P_k.$$

Then for $|P_1| > |P_2|$, we deduce from the asymptotic behavior of $Q, \nabla Q$ at infinity that

$$\left| (\nabla P_1 - \nabla P_2)|P_1 + P_2|^{p-1} - \nabla P_1|P_1|^{p-1} + \nabla P_2|P_2|^{p-1} \right|$$

$$= \left| \nabla P_1 |P_1|^{p-1} \left[ \left( 1 - \frac{\nabla P_2}{\nabla P_1} \right) \left| 1 + \frac{P_2}{P_1} \right|^{p-1} - 1 + \frac{\nabla P_2}{\nabla P_1} \left| \frac{P_2}{P_1} \right|^{p-1} \right] \right| \lesssim |P_1|^{p-1}|P_2|.\mathbb{1}_{y\cdot\frac{z}{|z|}>0}$$

and

$$\left| (\nabla P_1 - \nabla P_2)|P_1 + P_2|^{p-2}(P_1 + P_2) - \nabla P_1|P_1|^{p-2}P_1 + \nabla P_2|P_2|^{p-2}P_2 \right|$$

$$= \left| \nabla P_1 |P_1|^{p-2} P_1 \left[ \left( 1 - \frac{\nabla P_2}{\nabla P_1} \right) \left| 1 + \frac{P_2}{P_1} \right|^{p-2} \left( 1 + \frac{P_2}{P_1} \right) - 1 + \frac{\nabla P_2}{\nabla P_1} \left| \frac{P_2}{P_1} \right|^{p-2} \frac{P_2}{P_1} \right] \right|$$

$$\lesssim |P_1|^{p-1}|P_2|.\mathbb{1}_{y\cdot\frac{z}{|z|}>0}.$$

We do the same way in case $|P_2| > |P_1|$ and for function $(y - z_k)P_k$ thus we obtain from (5.20) that

$$|\dot{G}| \lesssim |\dot{z}|s^{-p} + |\dot{v}|s^{-p} \lesssim s^{-(p+1)}$$

so $\left\| \dot{A}_0 \right\|_{L^2} \lesssim \|\dot{G}\psi\|_{L^2} + |\dot{z}| \|G\nabla\psi\|_{L^2} \ll s^{-2}$. Next remark that for a function $f$

$$\left| \frac{d}{ds}\mathfrak{pr}_i(f) \right| \lesssim \left| \mathfrak{pr}_i(\dot{f}) \right| + |\dot{z}| |\mathfrak{pr}_i(f)|, \qquad i = 1, 2 \tag{5.37}$$

thus $\left\| \dot{\tilde{A}}_0 \right\|_{L^2} \ll s^{-2}$, by properties of $(-\Delta+1)^{-1}$, this implies $\|\dot{R}_0\|_{L^2} \ll s^{-2}$. We will prove by induction that

$$\|R_j\|_{L^2}, \forall j \geq 1.$$

For $A_j$ $(j \geq 1)$, we have

$$|\dot{A}_j| \lesssim \left| (\dot{P}_1 + \dot{P}_2 + \sum_{k=0}^{j-1} \dot{R}_k)|P_1 + P_2 + \sum_{k=0}^{j-1} R_k|^{p-1} - (\dot{P}_1 + \dot{P}_2 + \sum_{k=0}^{j-2} \dot{R}_k)|P_1 + P_2 + \sum_{k=0}^{j-2} R_k|^{p-1} \right|.$$

As $\|\dot{R}_k\|_{L^2} \ll s^{-2}$ for $0 \leq k < j$, it is sufficient to prove that

$$\left| (\dot{P}_1 + \dot{P}_2)|P_1 + P_2 + \sum_{k=0}^{j-1} R_k|^{p-1} - (\dot{P}_1 + \dot{P}_2)|P_1 + P_2 + \sum_{k=0}^{j-2} R_k|^{p-1} \right| \ll s^{-2}. \tag{5.38}$$

Let estimate

$$B_j = \left| (\nabla P_1 - \nabla P_2)|P_1 + P_2 + \sum_{k=0}^{j-1} R_k|^{p-1} - (\nabla P_1 - \nabla P_2)|P_1 + P_2 + \sum_{k=0}^{j-2} R_k|^{p-1} \right|.$$

We have three cases to consider, at a given point $x$, if $\max(|P_1|, |P_2|, |V_0|, ..., |V_{j-1}|) > \max(|P_1|, |P_2|)$ then

$$B_j \lesssim \sum_{k=0}^{j-1} |V_k|^p \lesssim s^{-p};$$

otherwise $\max(|P_1|, |P_2|, |V_0|, ..., |V_{j-1}|) = |P_1|$ then, by the first-order Taylor expansion

$$B_j = \left| \nabla P_1 |P_1|^{p-1} \left[ \frac{1 - \nabla P_2/\nabla P_1}{1 + P_2/P_1 + \sum_{k=0}^{j-1} R_k/P_1} \left( 1 + \frac{P_2}{P_1} \sum_{k=0}^{j-1} \frac{R_k}{P_1} \right) \left| 1 + \frac{P_2}{P_1} \sum_{k=0}^{j-1} \frac{R_k}{P_1} \right|^{p-1} \right. \right.$$

$$\left. \left. -(1 - \frac{\nabla P_2}{\nabla P_1}) \left| 1 + \frac{P_2}{P_1} + \sum_{k=0}^{j-2} \frac{R_k}{P_1} \right|^{p-1} \right] \right| \lesssim |P_1|^{p-1}|P_2|.\mathbb{1}_{y \cdot \frac{z}{|z|} > 0} + \sum_{k=0}^{j-1} |P_1|^{p-1}|R_k| \lesssim s^{-p},$$

and similarly for the case $\max(|P_1|, |P_2|, |V_0|, ..., |V_{j-1}|) = |P_2|$ so $B_j \lesssim s^{-p}$, from which we deduce (5.38). Recall the estimate for the derivative of a projection (5.37) so we get $\left\| \dot{A}_j \right\|_{L^2} \ll s^{-2}$. In conclusion, we have $\|G_0\|_{L^2} \lesssim s^{-2}$. Similarly, the same estimate holds for $\nabla G_0$, which finishes the proof of (5.10).

**step 4** Estimate of projection. From step 3, the terms whose norm $L^2$ is of order $s^{-2}$ are $A_{J+1}, \mathfrak{pr}_1(G\psi), \mathfrak{pr}_2(G\psi), G(1-\psi)$. As $|\langle \mathfrak{pr}_2(G\psi), e^{i\Gamma_1(y - z_1(s))} \nabla Q(y - z_1(s)) \rangle| \ll s^{-2}$ and similarly to (5.32), we can show $|\mathfrak{pr}_1(A_{J+1})| \ll s^{-2}$ thus

$$\langle G_0, e^{i\Gamma_1(y - z_1(s))} \nabla Q(y - z_1(s)) \rangle = \langle G, e^{i\Gamma_1(y - z_1(s))} \nabla Q(y - z_1(s)) \rangle + o(s^{-2}).$$

For $1 < p \leq 2$, we also have the analogous estimates of (2.23), (2.24)

$$\left| |\mathbf{P}|^{p-1}\mathbf{P} - |P_1|^{p-1}P_1 - |P_2|^{p-1}P_2 - \left[ \frac{p+1}{2}|P_1|^{p-1}P_2 + \frac{p-1}{2}|P_1|^{p-3}P_1^2\overline{P_2} \right].\mathbb{1}_{y \cdot \frac{z}{|z|} > 0} \right.$$

$$\left. - \left[ \frac{p+1}{2}|P_2|^{p-1}P_1 + \frac{p-1}{2}|P_2|^{p-3}P_2^2\overline{P_1} \right].\mathbb{1}_{y \cdot \frac{z}{|z|} < 0} \right| \lesssim \min(|P_1|^p, |P_2|^p). \tag{5.39}$$

We note that for $\delta = \frac{p-1}{2} > 0$

$$\left| \int_{y \cdot \frac{z}{|z|} > -\frac{|z|}{2}} Q^p(y+z)\nabla Q(y) dy \right| \lesssim |z|^{-\frac{d-1}{2}} e^{-|z|} Q^{(p-1)-\delta}\left(\frac{|z|}{2}\right) \int Q^\delta(y) dy$$

$$\lesssim s^{-(p+1-\delta)} \ll s^{-2}\log^{-1}(s),$$

$$\left| \int_{y \cdot \frac{z}{|z|} < -\frac{|z|}{2}} Q^p(y)\nabla Q(y) dy \right| \lesssim Q^{(p+1)-\delta}\left(\frac{|z|}{2}\right) \int Q^\delta(y) dy \lesssim s^{-(p+1-\delta)} \ll s^{-2}\log^{-1}(s).$$

We repeat the approach in step 3 of Lemma 2 and combine it with (5.7) to conclude that

$$\left| \langle G_0, e^{i\Gamma_1(y-z_1(s))}\nabla Q(y-z_1(s)) \rangle - C_p\frac{z}{|z|}|z|^{-\frac{d-1}{2}}e^{-|z|} \right| \lesssim |z|^{-1-\frac{d-1}{2}}e^{-|z|} \lesssim s^{-2}\log^{-1}(s)$$

as required.                                                                            □

The modulation part remains the same as for $p > 2$ (see Lemma 3) except the extra relation will be

$$\dot{v} = -\frac{2}{c_2}H_0(v,z) \tag{5.40}$$

where

$$H_0(v,z) = \left\langle G_0(y;(v(s),z(s))), e^{i\frac{v(s)}{2}(y-\frac{z(s)}{2}(s))}\nabla Q\left(y-\frac{z(s)}{2}\right) \right\rangle \tag{5.41}$$

$$= \langle G_0, e^{i\Gamma_1(y-z_1(s))}\nabla Q(y-z_1(s))\rangle.$$

Remark that by (5.11), the main order of $\dot{v}$ still remains

$$\left| \dot{v} + c\frac{z}{|z|}|z|^{-\frac{d-1}{2}}e^{-|z|} \right| \lesssim |z|^{-\frac{d-1}{2}-1}e^{-|z|}.$$

We claim the following analogue of Proposition 4 in the context $1 < p \leq 2$ for $L^2$ sub-critical.

**Proposition 10** (Uniform backwards estimates for $1 < p \leq 2$). *There exists $s_0 \gg 1$ satisfying the following condition: for all $s^{in} > s_0$, there is a choice of initial parameters $(\lambda^{in}, z^{in}, v^{in})$ such that the solution $u$ of (NLS) corresponding to (3.1) exists. Moreover, the decomposition of $u$ with extra relation (5.40) on the rescaled interval of time $[s_0, s^{in}]$*

$$u(s,x) = \frac{e^{i\gamma(s)}}{\lambda^{\frac{2}{p-1}}(s)}(\mathbf{P}+\epsilon)(s,y), \quad y = \frac{x}{\lambda(s)}, \quad dt = \lambda^2(s)ds$$

*verifies the uniform estimates for all $s \in [s_0, s^{in}]$*

$$||z(s)| - 2\log(s)| \lesssim \log(\log(s)), \quad |\lambda^{-1}(s) - 1| \lesssim s^{-1},$$

$$|v(s)| \lesssim s^{-1}, \quad \|\epsilon(s)\|_{H^1} \lesssim s^{-1}, \quad \left||z(s)|^{\frac{d-1}{2}}e^{|z(s)|} - cs^2\right| \lesssim s^2\log^{-\frac{1}{2}}(s). \tag{5.42}$$

*Proof of Proposition 10.* We only sketch the proof since it is very similar to Section 3.2, the main difference is the localization to avoid singularities due to the small power $p$ in Taylor expansions (5.1)–(5.3).

**step 1** Modulation equations. Consider

$$\frac{d}{ds}\langle \eta_1, A+iB\rangle = \langle \eta_1, iL_-A - L_+B\rangle - \langle \vec{m}_1 \cdot \vec{\mathbf{M}}\eta_1, iA - B\rangle - \langle \mathcal{E}_{\mathbf{P}_1}, iA - B\rangle$$

$$- \langle |\mathbf{P}_1 + \eta_1|^{p-1}(\mathbf{P}_1 + \eta_1) - |\mathbf{P}_1|^{p-1}\mathbf{P}_1 - \frac{p+1}{2}Q^{p-1}\eta_1 - \frac{p-1}{2}Q^{p-1}\overline{\eta}_1, iA - B\rangle$$

where the expression of $\mathbf{P}_1$ is given by

$$\mathbf{P}_1 = Q(y) + e^{i(\Gamma_2(y-(z_2-z_1))-\Gamma_1(y))}Q(y-(z_2-z_1)) + \sum_{j=0}^{J} e^{-i\Gamma_1(y)}R_j(y+z_1).$$

Let $C$ the set such that $\max(|R_0(y+z_1)|, ..., |R_J(y+z_1)|) \geq \frac{1}{J+2}Q(y)$ then for $y \in C$

$$|Q(y)| \lesssim \|R_i\|_{L^\infty} \leq s^{-p}, \qquad \text{for some } i \in \{0, ...J\}.$$

Since $|A|, |B| \lesssim |x|^q e^{-|x|}$, from the asymptotic behavior (1.12) of $Q$, over the set $C$, we have

$$|A| + |B| \lesssim s^{-p}\log^q s. \tag{5.43}$$

Next, denote

$$\Gamma(s,y) = \Gamma_2(y-(z_2-z_1)) - \Gamma_1(y) = -\frac{1}{2}iv \cdot (y+z) - \frac{1}{2}iv \cdot y, \tag{5.44}$$

from the estimates $||z| - 2\log(s)| \lesssim \log(\log(s))$ and $||v| - s^{-1}| \lesssim s^{-1}\log^{-1}(s)$, there exists a constant $c_0$ (independent of $s^{in}$) such that if $|y| \leq c_0 s$ then $|\Gamma(s,y)| \leq \frac{\pi}{2}$. Let $D = \{y \in \mathbb{R}^d, |y| > c_0 s\}$, we have for $y \in C^c \cap D^c$

$$\frac{1}{J+2}Q(y) \leq |\mathbf{P}_1(y)| \lesssim 1 \tag{5.45}$$

since $|R_0(y+z_1)|, ..., |R_J(y+z_1)| < \frac{1}{J+2}Q(y)$ and $\text{Re}[e^{i\Gamma}Q(y+z)] > 0$. And we have for $y \in C \cup D$, using $A, B \in \mathcal{Y}$ and (5.43),

$$|A(y)| + |B(y)| \lesssim \min(e^{-\frac{c_0}{2}s}, s^{-p}\log^q(s)) \lesssim s^{-1^+} \tag{5.46}$$

with $1^+ = \frac{p+1}{2}$. We denote

$$\varphi(s,y) = \mathbb{1}_{D^c}\mathbb{1}_{C^c}. \tag{5.47}$$

A consequence of (5.45) and (5.46) is that

$$|\mathbf{P}_1(y)|^{-m}Q(y)^n\varphi(s,y) \lesssim 1 \quad \text{for } n \geq m > 0 \tag{5.48}$$

and

$$(|A(y)| + |B(y)|)(1 - \varphi(s,y)) \lesssim s^{-1^+}. \tag{5.49}$$

By Cauchy-Schwarz and Gagliardo-Nirenberg inequalities

$$|\langle |\mathbf{P}_1 + \eta_1|^{p-1}(\mathbf{P}_1 + \eta_1) - |\mathbf{P}_1|^{p-1}\mathbf{P}_1 - \frac{p+1}{2}Q^{p-1}\eta_1 - \frac{p-1}{2}Q^{p-1}\overline{\eta}_1, (iA - B)(1 - \varphi(s,y))\rangle|$$

$$\lesssim \langle |\eta_1| + |\eta_1|^p, (iA+B)(1 - \varphi(s,\cdot))\rangle \lesssim s^{-1^+}(\|\eta_1\|_{H^1} + \|\eta_1\|_{H^1}^p) \lesssim C^* s^{-(1+1^+)}.$$

From the expansion in (5.1), we get

$$\left[ |\mathbf{P}_1 + \eta_1|^{p-1}(\mathbf{P}_1 + \eta_1) - |\mathbf{P}_1|^{p-1}\mathbf{P}_1 - \frac{p+1}{2}Q^{p-1}\eta_1 - \frac{p-1}{2}Q^{p-1}\overline{\eta}_1 \right]\varphi(s,y)$$

$$= \left[ \frac{p+1}{2}(|\mathbf{P}_1|^{p-1} - Q^{p-1})\eta_1 + \frac{p-1}{2}(|\mathbf{P}_1|^{p-3}\mathbf{P}_1^2 - Q^{p-1})\overline{\eta}_1 + O\left(\left|\frac{\eta_1}{\mathbf{P}_1}\right|^2|\mathbf{P}_1|^p\right) \right]\varphi(s,y).$$

We control the first two terms as before in the case $p > 2$

$$|\langle (|\mathbf{P}_1|^{p-1} - Q^{p-1})\eta_1, (iA - B)\varphi(s, \cdot)\rangle| + |\langle (|\mathbf{P}_1|^{p-3}\mathbf{P}_1^2 - Q^{p-1})\overline{\eta}_1, (iA - B)\varphi(s, \cdot)\rangle| \lesssim C^* s^{-(p+1)} \log^q(s)$$

and for the last term, we use (5.48) to remark that $|\mathbf{P}_1|^{p-2}|iA - B|\varphi(s, \cdot) \lesssim 1$ then deduce the inequality

$$\left\langle \left|\frac{\eta_1}{\overline{\mathbf{P}}_1}\right|^2 |\mathbf{P}_1|^p, (iA - B)\varphi(s, \cdot) \right\rangle \lesssim \|\epsilon\|_{L^2}^2 \lesssim (C^*)^2 s^{-2}.$$

To summarize, we have shown that

$$\left\langle |\mathbf{P}_1 + \eta_1|^{p-1}(\mathbf{P}_1 + \eta_1) - |\mathbf{P}_1|^{p-1}\mathbf{P}_1 - \frac{p+1}{2}Q^{p-1}\eta_1 - \frac{p-1}{2}Q^{p-1}\overline{\eta}_1, iA - B \right\rangle \lesssim s^{-2}. \tag{5.50}$$

Next, it is obvious that we still have as before

$$|\langle \vec{m}_1 \cdot \vec{\mathrm{M}}\eta_1, iA - B\rangle| \lesssim C^* s^{-1} |\vec{m}_1(s)|.$$

To prove the estimate

$$\left| \langle \mathcal{E}_{\mathbf{P}_1}, iA - B\rangle - \langle \vec{m}_1 \cdot \vec{\mathrm{M}}Q, iA - B\rangle \right| \lesssim s^{-2} + s^{-1}|\vec{m}_1|, \tag{5.51}$$

we recall $\mathcal{E}_{\mathbf{P}_1} = [\vec{m}_1 \cdot \vec{\mathrm{M}}Q](y) + [e^{i(\Gamma_2(y - (z_2 - z_1)) - \Gamma_1(y))}\vec{m}_2 \cdot \vec{\mathrm{M}}Q](y - (z_2 - z_1)) + e^{-i\Gamma_1(y)}G_0(y + z_1)$. From (5.10)

$$|\langle e^{-i\Gamma_1(y)}G_0(y + z_1), iA - B\rangle| \lesssim \|G_0\|_{L^2} \lesssim s^{-2}$$

and finally since $A, B \in \mathcal{Y}$, we have

$$|\langle e^{i(\Gamma_2(y - (z_2 - z_1)) - \Gamma_1(y))}(\vec{m}_2 \cdot \vec{\mathrm{M}}Q(\cdot - (z_2 - z_1))), iA - B\rangle| \lesssim s^{-1}|\vec{m}_1|,$$

which yields the estimate (3.28) in the case $1 < p \leq 2$. We project $\eta_1$ onto three null spaces of the linearized equation around $Q$ and obtain the almost orthogonality for the forth null space by the conservation of momentum thanks to the special choice of $\dot{v}$ in (5.41) (as in Section 3.2.2).

**step 2** Control the energy functional. We still consider the energy functional

$$\begin{aligned}
\mathbf{W}(s, \epsilon) &= \mathbf{H}(s, \epsilon) - \mathbf{J}(s, \epsilon) \\
&= \frac{1}{2} \int \left( |\nabla \epsilon|^2 + |\epsilon|^2 - \frac{2}{p+1} \left( |\mathbf{P} + \epsilon|^{p+1} - |\mathbf{P}|^{p+1} - (p+1)|\mathbf{P}|^{p-1}\mathrm{Re}(\epsilon\overline{\mathbf{P}}) \right) \right) \\
&\quad - \sum_{k=1}^{2} v_k \cdot \mathrm{Im} \int (\nabla \epsilon \, \overline{\epsilon})\chi_k
\end{aligned}$$

and remark that we still have the coercivity property

$$\mathbf{W}(s, \epsilon(s)) \gtrsim \|\epsilon(s)\|_{H^1}^2$$

(see for example [13], [17]). Define

$$\varphi_1(s, y) = \varphi(s, y - z_1(s)) \tag{5.52}$$

a function localized to the first soliton $\mathbf{P}_1$. Similarly, we can define an analogous function $\varphi_2(s, y)$ localized to the second soliton $\mathbf{P}_2$.

We claim an estimate on the derivative of $\mathbf{H}$ by $\dot{z}_k \cdot \langle \nabla P_k, \frac{\bar{\epsilon}.F''(\mathbf{P}).\epsilon}{2} \rangle$ but now localized by $\varphi_k$

$$\left| \frac{d}{ds}[\mathbf{H}(s, \epsilon(s))] - \sum_{k=1}^{2} \dot{z}_k \cdot \langle \varphi_k \nabla P_k, \frac{\bar{\epsilon}.F''(\mathbf{P}).\epsilon}{2} \rangle \right| \lesssim s^{-2} \|\epsilon(s)\|_{H^1} + s^{-2} \|\epsilon\|_{H^1}^2. \tag{5.53}$$

Recall that we have

$$\frac{d}{ds}[\mathbf{H}(s, \epsilon(s))] = D_s \mathbf{H}(s, \epsilon(s)) + \langle D_\epsilon \mathbf{H}(s, \epsilon(s)), \dot{\epsilon}_s \rangle,$$

and

$$D_s \mathbf{H} = \langle \dot{\mathbf{P}}, K \rangle, \quad \langle D_\epsilon \mathbf{H}(s, \epsilon), \dot{\epsilon} \rangle = \frac{\dot{\lambda}}{\lambda} \langle D_\epsilon \mathbf{H}(s, \epsilon), \Lambda \epsilon \rangle - (1 - \dot{\gamma}) \langle i D_\epsilon \mathbf{H}(s, \epsilon), \epsilon \rangle - \langle i D_\epsilon \mathbf{H}(s, \epsilon), \mathcal{E}_{\mathbf{P}} \rangle$$

with $K = |\mathbf{P} + \epsilon|^{p-1}(\mathbf{P} + \epsilon) - |\mathbf{P}|^{p-1}\mathbf{P} - \frac{p+1}{2}\epsilon|\mathbf{P}|^{p-1} - \frac{p-1}{2}\bar{\epsilon}\mathbf{P}^2|\mathbf{P}|^{p-3}$. We observe from (5.46) that for $\dot{P}_k = -\dot{z}_k \cdot \nabla P_k + i \dot{v}_k \cdot (y - z_k) P_k$, over the set $C \cup D$ , $|\dot{P}_k| \lesssim s^{-(1+1^+)}$ then

$$|\langle \dot{P}_k, K(1 - \varphi_k) \rangle| \lesssim s^{-(1+1^+)} \|\epsilon\|_{H^1}.$$

From (5.1), $|K| \lesssim |\epsilon|^2 |\mathbf{P}|^{p-2}$ so we obtain

$$|\langle i \dot{v}_k \cdot (y - z_k) P_k, K \varphi_k \rangle| \lesssim |\dot{v}| \, \|\epsilon\|_{H^1}^2 \lesssim s^{-2} \|\epsilon\|_{H^1}^2$$

since $\frac{Q(y - z_k)}{|\mathbf{P}|\varphi_k} \lesssim 1$ by (5.48). Next we look more precisely at $K$

$$K = \frac{\bar{\epsilon}.F''(\mathbf{P}).\epsilon}{2} + O\left(\left|\frac{\epsilon}{\mathbf{P}}\right|^{2^+} |\mathbf{P}|^p\right)$$

since $|\dot{z}_k| \lesssim s^{-1}$ and $p - 2^+ > -1$, we also have

$$\left| \left\langle -\dot{z}_k \cdot \nabla P_k, \left|\frac{\epsilon}{\mathbf{P}}\right|^{2^+} |\mathbf{P}|^p \varphi_k \right\rangle \right| \lesssim s^{-1} \|\epsilon\|_{H^1}^{2^+}.$$

We deal the first two terms of $\langle D_\epsilon \mathbf{H}(s, \epsilon), \dot{\epsilon} \rangle$ as in the case $p > 2$

$$\left| \frac{\dot{\lambda}}{\lambda} \langle D_\epsilon \mathbf{H}(s, \epsilon), \Lambda \epsilon \rangle \right| \lesssim \left| \frac{\dot{\lambda}}{\lambda} \right| \left( \|\epsilon\|_{H^1}^2 + \|\epsilon\|_{H^1}^{p+1} \right) \lesssim (C^*)^2 s^{-2} \|\epsilon\|_{H^1}^2,$$

$$|(1 - \dot{\gamma}) \langle i D_\epsilon \mathbf{H}(s, \epsilon), \epsilon \rangle| \lesssim |1 - \dot{\gamma}| \left( \|\epsilon\|_{H^1}^2 + \|\epsilon\|_{H^1}^{p+1} \right) \lesssim (C^*)^2 s^{-2} \|\epsilon\|_{H^1}^2.$$

Recall that for the last term we have

$$\langle i D_\epsilon \mathbf{H}(s, \epsilon), \mathcal{E}_{\mathbf{P}} \rangle = \langle -i\Delta\epsilon + i\epsilon - i\left(|\mathbf{P} + \epsilon|^{p-1}(\mathbf{P} + \epsilon) - |\mathbf{P}|^{p-1}\mathbf{P}\right),$$

$$[e^{i\Gamma_1}\vec{m}_1 \cdot \vec{M}Q](y - z_1(s)) + [e^{i\Gamma_2}\vec{m}_2 \cdot \vec{M}Q](y - z_2(s)) + G_0 \rangle$$

so from the properties of operators $L_+$ and $L_-$

$$I_1 = \langle -i\Delta\epsilon + i\epsilon - i\left(|\mathbf{P} + \epsilon|^{p-1}(\mathbf{P} + \epsilon) - |\mathbf{P}|^{p-1}\mathbf{P}\right), [e^{i\Gamma_1}\vec{m}_1 \cdot \vec{M}Q](y - z_1(s)) \rangle$$

$$= -\frac{\dot{\lambda}}{\lambda} \langle \eta_1, -2Q \rangle + (\dot{v} - \frac{\dot{\lambda}}{\lambda} v) \langle \eta_1, -2i\nabla Q \rangle$$

$$- \left\langle i\left(|\mathbf{P}_1 + \eta_1|^{p-1}(\mathbf{P}_1 + \eta_1) - |\mathbf{P}_1|^{p-1}\mathbf{P}_1 - \frac{p+1}{2}Q^{p-1}\eta_1 - \frac{p-1}{2}Q^{p-1}\bar{\eta}_1\right), \vec{m}_1 \cdot \vec{M}Q \right\rangle.$$

By the same way to prove (5.50), combining with the orthogonality of $\eta_1$ (3.5), (3.25) and the estimate of modulation equation (3.24), we get

$$|I_1| = O((C^*)^4 s^{-4}) + O((C^*)^2 s^{-4}).$$

Finally, using integration by parts and Cauchy-Schwarz inequality, from the bound for $H^1$ norm of $G_0$ (5.10), we obtain

$$|\langle -i\Delta\epsilon + i\epsilon - i\left(|\mathbf{P}+\epsilon|^{p-1}(\mathbf{P}+\epsilon) - |\mathbf{P}|^{p-1}\mathbf{P}\right), G_0\rangle| \lesssim s^{-2}\|\epsilon\|_{H^1}.$$

Combining these computations, the proof of (5.53) is finished. We now claim the estimate for localized momentum $J_k$: for all $s \in [s^*, s^{in}]$,

$$\left|\frac{d}{ds}[\mathbf{J}(s,\epsilon(s))] - \sum_{k=1}^{2} 2v_k \cdot \langle \varphi_k \nabla P_k, \frac{\bar{\epsilon}.F''(\mathbf{P}).\epsilon}{2}\rangle\right| \lesssim s^{-1}\log^{-1}(s)\|\epsilon(s)\|_{H^1}^2 + s^{-\frac{5}{2}}\|\epsilon(s)\|_{H^1}. \tag{5.54}$$

Recall that from the equation of $i\dot{\epsilon}$ (3.8), we have

$$\frac{d}{ds}[J_k(s,\epsilon(s))] = \dot{v}_k \cdot \text{Im} \int (\nabla\epsilon\,\bar{\epsilon})\chi_k + v_k \cdot \text{Im} \int (\nabla\epsilon\,\bar{\epsilon})\dot{\chi}_k$$

$$- v_k\langle \Delta\epsilon - \epsilon + \left(|\mathbf{P}+\epsilon|^{p-1}(\mathbf{P}+\epsilon) - |\mathbf{P}|^{p-1}\mathbf{P}\right) - i\frac{\dot{\lambda}}{\lambda}\Lambda\epsilon + (1-\dot{\gamma})\epsilon + \mathcal{E}_{\mathbf{P}}, 2\chi_k\nabla\epsilon + \epsilon\nabla\chi_k\rangle.$$

We proceed the same way as in Section 3.2.3 for $L^2$ sub-critical cases with $p > 2$, except for the term

$$v_k\left\langle|\mathbf{P}+\epsilon|^{p-1}(\mathbf{P}+\epsilon) - |\mathbf{P}|^{p-1}\mathbf{P}, 2\chi_k\nabla\epsilon + \epsilon\nabla\chi_k\right\rangle.$$

First, by (5.1)

$$|\mathbf{P}+\epsilon|^{p-1}(\mathbf{P}+\epsilon) - |\mathbf{P}|^{p-1}\mathbf{P} = F'(\mathbf{P}).\epsilon + O(|\epsilon|^p)$$

and then we have

$$|v_k||\langle|\epsilon|^p, 2\chi_k\nabla\epsilon + \epsilon\nabla\chi_k\rangle| \lesssim s^{-1}\|\epsilon\|_{H^1}^{p+1}.$$

Second, we consider

$$|v_k\langle F'(\mathbf{P}).\epsilon, \epsilon\nabla\chi_k\rangle| \lesssim |v_k|\,|\nabla\chi_k|\,\|\epsilon\|_{H^1}^2 \lesssim s^{-1}\log^{-1}(s)\|\epsilon\|_{H^1}^2.$$

Finally by integration by parts, we obtain

$$v_k\langle F'(\mathbf{P}).\epsilon, \chi_k\nabla\epsilon\rangle = -\frac{1}{2}v_k\langle \nabla\mathbf{P}\chi_k, \bar{\epsilon}.F''(\mathbf{P}).\epsilon\rangle - \frac{1}{2}v_k\langle F'(\mathbf{P}).\epsilon, \epsilon\nabla\chi_k\rangle.$$

These estimates yield

$$\frac{d}{ds}[\mathbf{J}(s,\epsilon(s))] = \sum_{k=1}^{2}\langle 2v_k \cdot \nabla\mathbf{P}\chi_k, \frac{\bar{\epsilon}.F''(\mathbf{P}).\epsilon}{2}\rangle + O(s^{-1}\log^{-1}(s)\|\epsilon(s)\|_{H^1}^2)$$

$$+ O(s^{-\frac{5}{2}}\|\epsilon(s)\|_{H^1}) + O(s^{-1}\|\epsilon(s)\|_{H^1}^{p+1}). \tag{5.55}$$

Since in the support of $\chi_k$, we have $|P_k| \gtrsim s^{-\frac{1}{8}} \geq \|V_j\|_{L^\infty}, \forall j = \overline{0, J}$ so $\varphi_k \equiv 1$ then

$$|v_k \cdot \langle \nabla\mathbf{P}\chi_k, \bar{\epsilon}.F''(\mathbf{P}).\epsilon\rangle - v_k \cdot \langle \varphi_k\nabla P_k, \bar{\epsilon}.F''(\mathbf{P}).\epsilon\rangle|$$

$$\lesssim |v_k|\left[\sum_{j\neq k}\left|\int_{|y-z_k(s)|<\frac{1}{8}\log s}(\bar{\epsilon}.F''(\mathbf{P}).\epsilon)\nabla P_j\right| + \left|\int_{|y-z_k(s)|>\frac{1}{10}\log s}\varphi_k(\bar{\epsilon}.F''(\mathbf{P}).\epsilon)\nabla P_k\right|\right]$$

$$\lesssim s^{-\frac{p-1}{10}-1}\|\epsilon\|_{H_1}^2$$

here we use the property (5.48) of $\varphi_k$ that $\varphi_k \neq 0$ implies $\left|\frac{\nabla P_i}{\mathbf{P}}\right| \lesssim 1$ and $\left|\frac{\nabla P_k}{\mathbf{P}}\right| \lesssim 1$. The proof of (5.54) is complete. Then we can deduce from the modulation equation $|\dot{z}_k - 2v_k| \lesssim s^{-1} \log^{-1}(s)$ that

$$\left|\frac{d}{ds}[\mathbf{W}(s, \epsilon(s))]\right| \lesssim s^{-2} \|\epsilon(s)\|_{H^1}.$$

The rest of the proof stays unchanged in comparison to the case $p > 2$ in Section 3.2.4. □

From the uniform backwards estimates in Proposition 10, since $\|R_j\|_{H^1} \ll s^{-1}$ for $j = \overline{0, J}$, we have that

$$\left\| u(t(s), x) - \frac{e^{i\gamma(s)}}{\lambda^{\frac{2}{p-1}}(s)} \sum_{k=1}^{2} \left[e^{i\Gamma_k} Q\right] \left(\frac{x}{\lambda(s)} + \frac{(-1)^k}{2} z(s)\right) \right\|_{H^1} \lesssim \|\epsilon(s)\|_{H^1} + \sum_{j=0}^{J} \|R_j(s)\|_{H^1}$$
$$\lesssim s^{-1}$$

then we proceed like in Section 4 to obtain the existence of a solution $u(t)$ satisfying the regime (1.6) in sub-critical cases with $1 < p \leq 2$

$$\left\| u(t) - e^{i\gamma(t)} \sum_{k=1}^{2} Q(. - x_k(t)) \right\|_{H^1} \lesssim \frac{1}{t}.$$

## 6. Super-critical cases

In this section, we will present the necessary modifications to prove the result in the $L^2$ super-critical cases $(1 + \frac{4}{d} < p < \frac{d+2}{d-2})$ (see [4]). For $k \in \{1, 2\}$, $z_1(s) = -z_2(s) = \frac{1}{2}z(s)$, $v_1(s) = -v_2(s) = \frac{1}{2}v(s)$, denote

$$Y_k^{\pm}(s, y) = e^{i\Gamma_k(s, y - z_k(s))} Y^{\pm}(y - z_k(s)) \tag{6.1}$$
$$Z_k(s, y) = e^{i\Gamma_k(s, y - z_k(s))} i\Lambda Q(y - z_k(s))$$
$$V_k(s, y) = e^{i\Gamma_k(s, y - z_k(s))} i\nabla Q(y - z_k(s))$$
$$W_k(s, y) = e^{i\Gamma_k(s, y - z_k(s))} (y - z_k(s)) Q(y - z_k(s)).$$

Let

$$\mathbf{Y}^{\pm}(s, y) = \mathbf{Y}^{\pm}(y; (z(s), v(s))) = \sum_{k=1}^{2} Y_k^{\pm}(s, y), \mathbf{Z}(s, y) = \mathbf{Z}(y; (z(s), v(s))) = \sum_{k=1}^{2} Z_k(s, y),$$

$$\mathbf{V}(s, y) = \mathbf{V}(y; (z(s), v(s))) = V_1(s, y) - V_2(s, y), \mathbf{W}(s, y) = \mathbf{W}(y; (z(s), v(s))) = W_1(s, y) - W_2(s, y).$$

We need some extra parameters to control the instability created by $Y^{\pm}$. Consider a solution of (NLS) with symmetric initial data like below: for $\mathfrak{b} = (\mathfrak{b}^+, \mathfrak{b}^-, \mathfrak{b}_1, \mathfrak{b}_2, \mathfrak{b}_3) \in \mathbb{R}^5$, $\|\mathfrak{b}\| \leq C(s^{in})^{-\frac{3}{2}}$ (the constant $C$ independent of $s^{in}$ and given in Lemma 11)

$$u(T_{mod}, x) = \frac{1}{(\lambda^{in})^{\frac{2}{p-1}}} w(s^{in}, y), \quad y = \frac{x}{\lambda^{in}} \tag{6.2}$$

with

$$w(s^{in}) = \mathbf{P}^{in}(y; (z^{in}\vec{e}_1, v^{in})) + \mathfrak{b}^+ i\mathbf{Y}^+(y; (z^{in}\vec{e}_1, v^{in})) + \mathfrak{b}^- i\mathbf{Y}^-(y; (z^{in}\vec{e}_1, v^{in}))$$
$$+ \mathfrak{b}_1 \mathbf{Z}(y; (z^{in}\vec{e}_1, v^{in})) + \mathfrak{b}_2 \mathbf{V}(y; (z^{in}\vec{e}_1, v^{in})) + \mathfrak{b}_3 \mathbf{W}(y; (z^{in}\vec{e}_1, v^{in})). \tag{6.3}$$

Then we get

$$\epsilon(s^{in}) = \mathfrak{b}^+ i\mathbf{Y}^+(y; (z^{in}\vec{e}_1, v^{in})) + \mathfrak{b}^- i\mathbf{Y}^-(y; (z^{in}\vec{e}_1, v^{in}))$$
$$+ \mathfrak{b}_1 \mathbf{Z}(y; (z^{in}\vec{e}_1, v^{in})) + \mathfrak{b}_2 \mathbf{V}(y; (z^{in}\vec{e}_1, v^{in})) + \mathfrak{b}_3 \mathbf{W}(y; (z^{in}\vec{e}_1, v^{in})).$$

**Lemma 11** (Modulated data in direction $Y^\pm$). *There exists $C > 0$ such that for all $s^{in} \geq s_0$ and for all $a^{in} \in [-(s^{in})^{-\frac{3}{2}}, (s^{in})^{-\frac{3}{2}}]$, there is a unique $\mathfrak{b}$ so that $||\mathfrak{b}|| \leq C|a^{in}|$ (C independent of $s^{in}$) and the initial data satisfies*

$$\langle \eta_1(s^{in}), iY^-\rangle = a^{in}, \quad \langle \eta_1(s^{in}), iY^+\rangle = \langle \eta_1(s^{in}), i\Lambda Q\rangle = \langle \eta_1(s^{in}), yQ\rangle = \langle \eta_1(s^{in}), i\nabla Q\rangle = 0 \tag{6.4}$$

*with $\eta_1$ defined as in (3.4).*

*Proof of Lemma 11.* Let

$$\mathfrak{c} = (\langle \eta_1(s^{in}), iY^+\rangle, \langle \eta_1(s^{in}), iY^-\rangle, \langle \eta_1(s^{in}), i\Lambda Q\rangle, \langle \eta_1(s^{in}), i\nabla Q\rangle, \langle \eta_1(s^{in}), yQ\rangle).$$

We consider the linear maps

$$\Psi : \mathbb{R}^5 \to H^1(\mathbb{R}^d) \qquad \Phi : H^1(\mathbb{R}^d) \to \mathbb{R}^5$$
$$\mathfrak{b} \mapsto \epsilon(s^{in}) \qquad\qquad \epsilon(s^{in}) \to \mathfrak{c}$$

and $\Omega = \Phi \circ \Psi : \mathbb{R}^5 \to \mathbb{R}^5$. We compute

$$\Psi(h) = (i\mathbf{Y}^+(y; (z^{in}\vec{e}_1, v^{in})), i\mathbf{Y}^-(y; (z^{in}\vec{e}_1, v^{in})), \mathbf{Z}(y; (z^{in}\vec{e}_1, v^{in})), \mathbf{V}(y; (z^{in}\vec{e}_1, v^{in})), \mathbf{W}(y; (z^{in}\vec{e}_1, v^{in})) \cdot h$$

$$\Phi(v) = \begin{pmatrix} \int v(y)[e^{-i\Gamma_1}\overline{iY^+}](y - \frac{1}{2}z^{in}\vec{e}_1)dy \\ \int v(y)[e^{-i\Gamma_1}\overline{iY^-}](y - \frac{1}{2}z^{in}\vec{e}_1)dy \\ \int v(y)[e^{-i\Gamma_1}\overline{i\Lambda Q}](y - \frac{1}{2}z^{in}\vec{e}_1)dy \\ \int v(y)[e^{-i\Gamma_1}\overline{i\nabla Q}](y - \frac{1}{2}z^{in}\vec{e}_1)dy \\ \int v(y)[e^{-i\Gamma_1}\overline{yQ}](y - \frac{1}{2}z^{in}\vec{e}_1)dy \end{pmatrix}$$

then we can deduce that for some complex functions $A(y), B(y) \in \mathbf{Y}$

$$\Omega = \Phi \circ \Psi = N + O\big(|\langle A(y + z^{in}\vec{e}_1), B(y)\rangle|\big) = N + O(e^{-|z^{in}|})$$

where

$$N = \begin{pmatrix} \langle iY^+, iY^+\rangle & \langle iY^-, iY^+\rangle & \langle i\Lambda Q, iY^+\rangle & \langle i\nabla Q, iY^+\rangle & \langle yQ, iY^+\rangle \\ \langle iY^+, iY^-\rangle & \langle iY^-, iY^-\rangle & \langle i\Lambda Q, iY^-\rangle & \langle i\nabla Q, iY^-\rangle & \langle yQ, iY^-\rangle \\ \langle iY^+, i\Lambda Q\rangle & \langle iY^-, i\Lambda Q\rangle & \langle i\Lambda Q, i\Lambda Q\rangle & \langle i\nabla Q, i\Lambda Q\rangle & \langle yQ, i\Lambda Q\rangle \\ \langle iY^+, i\nabla Q\rangle & \langle iY^-, i\nabla Q\rangle & \langle i\Lambda Q, i\nabla Q\rangle & \langle i\nabla Q, i\nabla Q\rangle & \langle yQ, i\nabla Q\rangle \\ \langle iY^+, yQ\rangle & \langle iY^-, yQ\rangle & \langle i\Lambda Q, yQ\rangle & \langle i\nabla Q, yQ\rangle & \langle yQ, yQ\rangle \end{pmatrix}$$

and $\Omega(0) = 0$. Remark that $N$ is the Gramian matrix of $iY^+, iY^-, i\Lambda Q, i\nabla Q, yQ$ which are linearly independent since if for some $m, n, p, q, r \in \mathbb{R}$ (not all zeros)

$$m\, iY^+ + n\, iY^- + p\, i\Lambda Q + q\, yQ + r\, i\nabla Q = 0$$

then $m\,Y^+ + n\,Y^- + p\,\Lambda Q - q\,iyQ + r\,\nabla Q = 0$. We apply $\mathcal{L}$ to both sides of the equality $(L_+(\Lambda Q) = -2Q, L_-(xQ) = -2\nabla Q, L_+(\nabla Q) = 0)$ and get

$$me_0 Y^- - ne_0 Y^- - 2piQ - 2q\nabla Q = 0$$

so $m = n = p = q = 0$ as $Y^+, Y^-, iQ, \nabla Q$ are linearly independent thus $r = 0$, a contradiction. Therefore, $\det N \neq 0$ and with $|z^{in}| \gg 1$, we have that $\Omega$ is invertible around 0 and

$$||\Omega^{-1}|| \leq ||\mathrm{Gram}(iY^+, iY^-, i\Lambda Q, i\nabla Q, yQ)|| + 2$$

Therefore, for any $a^{in} \in [-(s^{in})^{-\frac{3}{2}}, (s^{in})^{-\frac{3}{2}}]$, we can choose

$$\mathfrak{b} = \Omega^{-1}((0, a^{in}, 0, 0, 0)), \qquad ||\mathfrak{b}|| \leq ||\Omega^{-1}||\,|a^{in}|$$

to conclude the lemma. $\qquad\square$

In fact, the coefficients $\mathfrak{b}_1$, $\mathfrak{b}_2$, $\mathfrak{b}_3$ can be determined explicitly from $\mathfrak{b}^+, \mathfrak{b}^-$ as follows

$$\mathfrak{b}_1 = \frac{\mathfrak{b}^+ \langle iY^+, i\Lambda Q \rangle + \mathfrak{b}^+ \langle e^{i\Gamma_0(\cdot)} iY^+ (\cdot + z^{in}\vec{e}_1), i\Lambda Q \rangle + \mathfrak{b}^- \langle iY^-, i\Lambda Q \rangle + \mathfrak{b}^- \langle e^{i\Gamma_0(\cdot)} iY^- (\cdot + z^{in}\vec{e}_1), i\Lambda Q \rangle}{\|\Lambda Q\|_{L^2}^2 + \langle e^{i\Gamma_0(\cdot)} i\Lambda Q (\cdot + z^{in}\vec{e}_1), i\Lambda Q \rangle}$$

$$(6.5)$$

$$\mathfrak{b}_2 = \frac{\mathfrak{b}^+ \langle iY^+, i\nabla Q \rangle + \mathfrak{b}^+ \langle e^{i\Gamma_0(\cdot)} iY^+ (\cdot + z^{in}\vec{e}_1), i\nabla Q \rangle + \mathfrak{b}^- \langle iY^-, i\nabla Q \rangle + \mathfrak{b}^- \langle e^{i\Gamma_0(\cdot)} iY^- (\cdot + z^{in}\vec{e}_1), i\nabla Q \rangle}{\|\nabla Q\|_{L^2}^2 - \langle e^{i\Gamma_0(\cdot)} [i\nabla Q] (\cdot + z^{in}\vec{e}_1), i\nabla Q \rangle}$$

$$(6.6)$$

$$\mathfrak{b}_3 = \frac{\mathfrak{b}^+ \langle iY^+, yQ \rangle + \mathfrak{b}^+ \langle e^{i\Gamma_0(\cdot)} iY^+ (\cdot + z^{in}\vec{e}_1), yQ \rangle + \mathfrak{b}^- \langle iY^-, yQ \rangle + \mathfrak{b}^- \langle e^{i\Gamma_0(\cdot)} iY^- (\cdot + z^{in}\vec{e}_1), yQ \rangle}{\|yQ\|_{L^2}^2 - \langle e^{i\Gamma_0(\cdot)} [yQ] (\cdot + z^{in}\vec{e}_1), yQ \rangle}$$

$$(6.7)$$

where $\Gamma_0(y) = -\frac{1}{2} iv^{in} \cdot (y + z^{in}\vec{e}_1) - \frac{1}{2} iv^{in} \cdot y$. This specific choice is made in order that initially, we have the following orthogonality conditions

$$\langle \eta_1(s^{in}), i\Lambda \rangle = \langle \eta_1(s^{in}), yQ \rangle = 0 \qquad (6.8)$$

and $\langle \eta_1(s^{in}), i\nabla Q \rangle = 0$. We recall the decomposition of $u(t)$: there exists a $\mathcal{C}^1$ function

$$\vec{q}(t) = (\lambda, z, \gamma, v) : [s_0, s^{in}] \to (0, +\infty) \times \mathbb{R}^d \times \mathbb{R} \times \mathbb{R}^d \text{ as}$$

such that we can modulate $u(t)$ on $[s_0, s^{in}]$ as

$$u(t(s), x) = \frac{e^{i\gamma(s)}}{\lambda(s)} (\mathbf{P} + \epsilon)(s, y)$$

and $\langle \eta_1(s), i\Lambda \rangle = \langle \eta_1(s), yQ \rangle = 0$. In here we obtain only two orthogonality conditions as the initial data satisfies only two (6.8). The proof of uniform estimates will remain the same except for some modifications that we will clarify immediately. Denote

$$a^{\pm}(s) = \langle \eta_1(s), iY^{\pm} \rangle, \qquad (6.9)$$

Lemma 11 allows us to establish a one-to-one mapping between the choice of $(\mathfrak{b}^+, \mathfrak{b}^-)$ and the constraints $a^+(s^{in}) = 0, a^-(s^{in}) = a^{in}$ for any choice of $a^{in}$. We now define the maximal time interval $[S(a^{in}), s^{in}]$ on which (3.21) holds and

$$|a^{\pm}(s)| \leq s^{-\frac{3}{2}} \qquad (6.10)$$

for all $s \in [S(a^{in}), s^{in}]$. We will prove that there exists a choice of

$$a^{in} \in [-(s^{in})^{-\frac{3}{2}}, (s^{in})^{-\frac{3}{2}}]$$

and $z^{in}$ such that $S(a^{in}) = s_0$. The first thing changed is that $\epsilon(s^{in})$ may not be zero, but we still have $\epsilon(s^{in}) \lesssim ||\mathfrak{b}|| \lesssim (s^{in})^{-\frac{3}{2}}$. This is enough to conclude that $|\mathbf{W}(s, \epsilon(s))| \lesssim C^* s^{-2}$ from the fact $\left|\frac{d}{ds}\mathbf{W}(s, \epsilon(s))\right| \lesssim C^* s^{-3}$. Next, from $\langle \eta_1(s^{in}), i\nabla Q \rangle = 0$, we deduce that

$$|M(u^{in}) - M(\mathbf{P}^{in})| \lesssim (C^*)^2 (s^{in})^{-2}$$

thus we still get $|\langle \eta_1, i\nabla Q \rangle| \lesssim (C^*)^2 s^{-2}$ from the fact $\left|\frac{d}{ds}M(\mathbf{P})\right| \lesssim (C^*)^2 s^{-3}$. The second thing which need to be modified is the coercivity of $\mathbf{W}$. By (1.17)

$$\mathbf{W}(s, \epsilon(s)) \gtrsim \|\epsilon(s)\|^2_{H^1} + O(s^{-3})$$

the process in Section 3 is still valid as long as we have (6.10). We claim the following preliminary estimates on the parameters $a^{\pm}(s)$.

**Lemma 12.** *For all $s \in [S(a^{in}), s^{in}]$,*

$$\left|\frac{da^{\pm}}{ds}(s) \mp e_0 a^{\pm}(s)\right| \lesssim ||\epsilon||^2_{H^1} \tag{6.11}$$

*Proof of Lemma 12.* Applying the inequality (3.28) with $A = -\operatorname{Im} Y^+, B = \operatorname{Re} Y^+$ and using the equation of $Y^{\pm}$ (1.16)

$$\left|\frac{d}{ds}\langle \eta_1, i\operatorname{Re} Y^+ - \operatorname{Im} Y^+ \rangle - \big[\langle \eta_1, -iL_-(\operatorname{Im} Y^+) - L_+(\operatorname{Re} Y^+)\rangle\right.$$

$$\left. - \langle \vec{m}_1 \cdot \vec{\mathbf{M}}Q, -i\operatorname{Im} Y^+ - \operatorname{Re} Y^+ \rangle\big]\right| \lesssim (C^*)^2 s^{-2} + s^{-1}|\vec{m}_1| \tag{6.12}$$

so we get

$$\left|\frac{d}{ds}\langle \eta_1, iY^+ \rangle - \langle \eta_1, i\mathcal{L}(Y^+)\rangle\right| \lesssim (C^*)^2 s^{-2} + s^{-1}|\vec{m}_1| + |\langle \vec{m}_1 \cdot \vec{\mathbf{M}}Q, Y^+ \rangle|.$$

This implies $\left|\frac{da^+}{ds}(s) - e_0 a^+(s)\right| \lesssim ||\epsilon||^2_{H^1}$. In the same way, we also obtain

$$\left|\frac{da^-}{ds}(s) + e_0 a^-(s)\right| \lesssim ||\epsilon||^2_{H^1}$$

as desired.                                                                                 $\square$

By the same arguments in Section 3, we improve all estimates in the bootstrap bounds except those of $a^{\pm}(s)$ and $z(s)$. It seems to us that the reasoning to close the bootstrap bound of $z(s)$ still works, in fact, it is, however we will control $a^{\pm}(s)$ through a suitable value of $a^{in}$ also by a topological argument so we have to choose $(z^{in}, a^{in})$ in the same time.

**Lemma 13** (Control of $a^+(s)$). *For all $a^{in} \in [-(s^{in})^{-\frac{3}{2}}, (s^{in})^{-\frac{3}{2}}]$, the following inequality holds for all $s \in [S(a^{in}), s^{in}]$*

$$|a^+(s)| \leq \frac{1}{2}s^{-\frac{3}{2}}. \tag{6.13}$$

*Proof of Lemma 13.* It follows (3.21), (6.11) and $a^+(s^{in}) = 0$ that for all $s \in [S(a^{in}), s^{in}]$

$$|a^+(s)| \lesssim (C^*)^2 e^{e_0 s} \int_s^{s^{in}} e^{-e_0 \tau} \tau^{-2} d\tau$$

$$= \frac{(C^*)^2}{e_0} e^{e_0 s} [e^{-e_0 s} s^{-2} - e^{-e_0 s^{in}} (s^{in})^{-2}] - 2 \frac{(C^*)^2}{e_0} e^{e_0 s} \int_s^{s^{in}} e^{-e_0 \tau} \tau^{-3} d\tau$$

$$\leq \frac{(C^*)^2}{e_0} s^{-2} \leq \frac{1}{2} s^{-\frac{3}{2}}$$

for $s_0$ to be large enough. $\qed$

**Lemma 14** (Control of $a^-(s)$ and closing the parameter $z$)**.** *There exist $z^{in}$ and $a^{in} \in [-(s^{in})^{-\frac{3}{2}}, (s^{in})^{-\frac{3}{2}}]$ such that $S(a^{in}) = s_0$.*

*Proof of Lemma 14.* We argue by contradiction. Consider $\zeta(s)$, $\xi(s)$ as defined in (3.45) and

$$\mathcal{N}(s) = s^3 (a^-(s))^2.$$

Suppose for all $(\zeta^\sharp, a^\sharp) \in \mathbb{D} = [-1, 1] \times [-1, 1]$, the choice of

$$\zeta^{in} = s^{in} + \zeta^\sharp s^{in} \log^{-\frac{1}{2}}(s^{in}), \quad a^{in} = a^\sharp (s^{in})^{-\frac{3}{2}}$$

gives us $S(a^{in}) = S(\zeta^\sharp, a^\sharp) \in (s_0, s^{in})$. Recall that

$$\dot{\xi}(s) = 2(\zeta(s) - s)(\dot{\zeta}(s) - 1)s^{-2} \log(s) - (\zeta(s) - s)^2 (2s^{-3} \log(s) - s^{-3}). \tag{6.14}$$

On the other hand, for $s \in (S(\zeta^\sharp, a^\sharp), s^{in}]$, then by (3.21) and (6.11), we have

$$\dot{\mathcal{N}}(s) = s^3 (3s^{-1} a^-(s) + 2 \frac{da^-}{ds}(s)) a^-(s)$$

$$= s^3 (3s^{-1} - 2e_0)(a^-(s))^2 + O\left(||\epsilon||_{H^1}^2 s^3 |a^-(s)|\right).$$

Due to the bound on $||\epsilon||_{H^1}^2$, we obtain

$$\dot{\mathcal{N}}(s) \leq s^3 (3s^{-1} - 2e_0)(a^-(s))^2 + C(C^*)^2 s^{-\frac{1}{2}} \sqrt{\mathcal{N}(s)}$$

then for $s_0$ large enough ($\frac{3}{s_0} < \frac{1}{2} e_0$ and $C(C^*)^2 s_0^{-\frac{1}{2}} < \frac{1}{2} e_0$), the estimate becomes

$$\dot{\mathcal{N}}(s) \leq -\frac{3}{2} e_0 \mathcal{N}(s) + C(C^*)^2 s^{-\frac{1}{2}} \sqrt{\mathcal{N}(s)}. \tag{6.15}$$

Denote

$$\Psi_1(s) = (\zeta(s) - s)(s)^{-1} \log^{\frac{1}{2}}(s),$$

$$\Psi_2(s) = a^-(s)(s)^{\frac{3}{2}}.$$

From the definition of $S(a^{in})$ and the continuity of flow, at the limit $S(\zeta^\sharp, a^\sharp)$, we have one of the following situation

$$\Psi_1(S(\zeta^\sharp, a^\sharp)) = \pm 1, \qquad \Psi_2 \in [-1, 1] \tag{6.16}$$

or

$$\Psi_2(S(\zeta^\sharp, a^\sharp)) = \pm 1, \qquad \Psi_1 \in [-1, 1]. \tag{6.17}$$

Remark that in the first case, we have

$$\dot{\xi}(S(\zeta^\sharp, a^\sharp)) < -(S(\zeta^\sharp, a^\sharp))^{-1} < 0$$

and in the second case we have $\mathcal{N}(S(\zeta^\sharp, a^\sharp)) = 1$

$$\dot{\mathcal{N}}(S(\zeta^\sharp, a^\sharp)) \leq -e_0 < 0.$$

A consequence of the above transversality property is the continuity of the map $(\zeta^\sharp, a^\sharp) \mapsto S((\zeta^\sharp, a^\sharp))$ thus the following map

$$\Psi \quad : \quad \mathbb{D} \quad \to \quad \partial\mathbb{D}$$
$$(\zeta^\sharp, a^\sharp) \quad \mapsto \quad (\Psi_1(S(\zeta^\sharp, a^\sharp)), \Psi_2(S(\zeta^\sharp, a^\sharp)))$$

is also continuous where $\partial\mathbb{D}$ is the boundary of $\mathbb{D}$. Note that if $a^\sharp = \pm 1$, then from (6.15), $\dot{\mathcal{N}}(s^{in}) < 0$, we have $S(\zeta^\sharp, a^\sharp) = s^{in}$ and if $\zeta^\sharp = \pm 1$, then from (6.14), $\dot{\xi}(s^{in}) < 0$, we also have $S(\zeta^\sharp, a^\sharp) = s^{in}$. Thus $\Psi(\zeta^\sharp, a^\sharp) = (\zeta^\sharp, a^\sharp)$ for all $(\zeta^\sharp, a^\sharp) \in \partial\mathbb{D}$, which means that the restriction of $\Psi$ to the boundary of $\mathbb{D}$ is the identity. But the existence of such a map contradicts the Brouwer fixed point theorem. In conclusion, there exists a final data $(z^{in}, a^{in})$ such that $S(a^{in}) = s_0$.                                                                 □

Finally, we still have the strong compactness result as in Lemma 8

$$u_n(t_0) \rightharpoonup u_0 \text{ weakly in } H^1(\mathbb{R}^d)$$

$$u_n(t_0) \to u_0 \text{ in } H^\sigma(\mathbb{R}^d), \text{ for } 0 \leq \sigma < 1$$

then we also consider $u$ the solution of (NLS) corresponding to $u_0$, by local well-posedness and continuous dependence (in [3]) for $L^2$ super-critical of (NLS), we have for all $t \in [t_0, +\infty)$,

$$u_n(t) \to u(t) \quad \text{in } H^\sigma(\mathbb{R}^d), \quad s_c \leq \sigma < 1$$

where $s_c$ is the critical exponent $s_c = \frac{d}{2} - \frac{2}{p-1} < 1$. Thus we can pass to the limit the decomposition $(\vec{q}, \epsilon)$ and get

$$\left\| u(t) - e^{i\gamma(t)} \sum_{k=1}^{2} Q\left(x + \frac{(-1)^k}{2} \delta x(t)\right) \right\|_{H^1} \lesssim t^{-1}. \tag{6.18}$$

**Acknowledgements.** This paper has been prepared as a part of my Master and PhD under the supervision of Y. Martel. I would like to thank my advisor for his constant support and many enlightening discussions. I also want to thank P. Raphaël for suggesting this work and helpful comments.

## References

[1] S. Agmon. Lectures on exponential decay of solutions of second-order elliptic equations: bounds on eigenfunctions of $N$-body Schrödinger operators. *Mathematical Notes (29)*, Princeton University Press, Princeton, 1982.

[2] T. Cazenave. Semilinear Schrödinger equations. *Courant Lecture Notes in Mathematics* , New York University, New York, 2003.

[3] T. Cazenave and F. B. Weissler. The Cauchy problem for the critical nonlinear Schrödinger equation in $H^s$. *Nonlinear Anal.*, 14(10):807–836, 1990.

[4] R. Côte, Y. Martel and F. Merle. Construction of multi-soliton solutions for the $L^2$-supercritical gKdV and NLS equations. *Rev. Mat. Iberoam.* 27 (2011), no. 1, 273–302.

[5] V. Combet. Multi-existence of multi-solitons for the supercritical nonlinear Schrödinger equation in one dimension. *Discrete Cont. Dyn. Sys.*, 34 (2014), 1961 – 1993.

[6] T. Duyckaerts and F. Merle. Dynamic of threshold solutions for energy-critical NLS. *Geom. Funct. Anal.* 18 (2009), no. 6, 1787–1840.

[7] T. Duyckaerts and S. Roudenko. Threshold solutions for the focusing 3D cubic Schrödinger equation. *Rev. Mat. Iberoam*, 26 (2010), 1 – 56.

[8] S.I. Ei and T. Ohta. Equation of motion for interacting pulses. *Physical Review E.*, 50 (1994), 4672–4678.

[9] J. Ginibre and G. Velo. On a class of nonlinear Schrödinger equations. I. The Cauchy problem, general case. *J. Funct. Anal.*, 32 (1979), 1 – 32.

[10] M. Grillakis, J. Shatah and W. A. Strauss. Stability theory of solitary waves in the presence of symmetry. *J. Funct. Anal.*, 197 (1987), 74 – 160.

[11] M. Grillakis. Analysis of the linearization around a critical point of an infinite dimensional Hamiltonian system. *Comm. Pure Appl. Math.*, 43 (1990), 299 – 333.

[12] M. A. Herrero and J. J. L. Velázquez. Flat blow-up in one-dimensional semilinear heat equations. *Differential Integral Equations* 5 (1992), 973–997.

[13] J. Krieger, Y. Martel and P. Raphaël. Two-soliton solutions to the three-dimensional gravitational Hartree equation. *Comm. Pure Appl. Math.* 62 (2009), no. 11, 1501–1550.

[14] E. Olmedilla. Multiple pole solutions of the nonlinear Schrödinger equation. *Physica D.* 25 (1987), 330–346.

[15] Y. Martel and F. Merle. Multi-solitary waves for nonlinear Schrödinger equations. *Annales de l'IHP (C) Non Linear Analysis*, 23 (2006), 849–864.

[16] Y. Martel and F. Merle. Description of two soliton collison for the quartic gKdV equation. *Ann. of Math.*, 174 (2011), 757–857.

[17] Y. Martel and F. Merle. Inelastic interaction of nearly equal solitons for the quartic gKdV equation. *Invent. Math.* 183 (2011), no. 3, 563–648.

[18] Y. Martel, F. Merle and T.-P. Tsai. Stability in $H^1$ of the sum of $K$ solitary waves for some nonlinear Schrödinger equations. *Duke Math. J.* 133 (2006), 405–466.

[19] Y. Martel and P. Raphael. Strongly interacting blow up bubbles for the mass critical NLS. *Preprint,* arXiv:1512.00900.

[20] F. Merle. Construction of solutions with exactly $k$ blow-up points for the Schrödinger equation with critical nonlinearity. *Comm. Math. Phys.* 129 (1990), no. 2, 223–240.

[21] F. Merle and P. Raphaël. On universality of blow-up profile for $L^2$ critical nonlinear Schrödinger equation. *Invent. Math.*, 156(3):565–672, 2004.

[22] F. Merle and P. Raphaël. The blow-up dynamic and upper bound on the blow-up rate for critical nonlinear Schrödinger equation. *Ann. of Math. (2)*, 161(1):157–222, 2005.

[23] P. Raphaël. Stability and blow up for the nonlinear Schr ödinger equation. *Lecture notes for the Clay summer school on evolution equations*, ETH, Zurich (2008).

[24] P. Raphaël and J. Szeftel. Existence and uniqueness of minimal blow-up solutions to an inhomogeneous mass critical NLS. *J. Amer. Math. Soc.*, 24(2):471–546, 2011.

[25] M. I. Weinstein. Modulational stability of ground states of nonlinear Schrödinger equations. *SIAM J. Math. Anal.*, 16 (1985), 472–491.

[26] M. I. Weinstein. Lyapunov stability of ground states of nonlinear dispersive evolution equations. *Comm. Pure Appl. Math.*, 39 (1986), 51–68.

[27] T. Zakharov and A.B. Shabat. Exact theory of two-dimensional self-focusing and one-dimensional self-modulation of waves in nonlinear media. *Sov. Phys. JETP* 34 (1972), 62–69.

CMLS, École Polytechnique, CNRS, Université Paris-Saclay, 91128 Palaiseau, France
*E-mail address*: `tien-vinh.nguyen@polytechnique.edu`